\documentclass[11pt,leqno]{article}

\usepackage{amsmath,amsthm}
\usepackage {latexsym}
\usepackage{amssymb}

\newcommand{\les}{\lesssim}

\newcommand{\bea}{\begin{eqnarray}}
\newcommand{\eea}{\end{eqnarray}}
\newcommand{\beq}{\begin{equation}}
\newcommand{\ee}{\end{equation}}

\newcommand{\R}{{\mathbb R}}

\newcommand{\si}{\sigma}
\renewcommand{\b}{\beta}

\def\pa{\partial}

\def\nn{\nonumber}

\newtheorem{theorem}{Theorem}[section]
\newtheorem{lemma}[theorem]{Lemma}
\newtheorem{defi}[theorem]{Definition}
\newtheorem{cor}[theorem]{Corollary}
\newtheorem{prop}[theorem]{Proposition}
\newtheorem{proposition}[theorem]{Proposition}

\theoremstyle{remark}
\newtheorem{remark}[theorem]{Remark}
\def\ve{\varepsilon}

\def\a{\alpha}
\def\ga{\gamma}
\def\de{\delta}
\def\Si{\Sigma}
\def\Ga{\Gamma}

\def\il{\int\limits}
\def\bm{\left( \begin{array}{cc}}
\def\endm{\end{array}\right)}
\newcommand{\eq}{\end{equation}}

\def\nab{\nabla}

\def\tr{\text{tr}}
\def\a{\alpha}
\def\b{\beta}
\def\ga{\gamma}
\def\de{\delta}
\def\Box{\square}
\def\Boxr{\widetilde{\square}}
\def\pa{\partial}
\def\pab{\bar\pa}

\def \rectangle#1#2{\hbox{\vrule\vbox to #2
              {\hrule\hbox to #1{\hfil}\vfil\hrule}\vrule}}
\def\sq{\,\,\rectangle{7pt}{7 pt}\,\,}
\def\Lb{\underline{L}}

\def\pas{\text{$\pa\mkern -9.0mu$\slash\,}}

\numberwithin{equation}{section}

\begin{document}

\title {Global existence for the Einstein vacuum equations\\
in wave coordinates }
\author {Hans Lindblad \thanks
{Part of this work was done while H.L. was a Member of the
Institute for  Advanced Study, Princeton, supported by the NSF grant
DMS-0111298 to the Institute. H.L. was also partially supported
by the NSF Grant DMS-0200226. } \,\,and Igor Rodnianski
\thanks{Part of this work was done while I.R. was a
Clay Mathematics Institute Long-Term
Prize Fellow.  His work was also partially supported by the NSF grant DMS--01007791.}\\
University of California at San Diego and  Princeton University}
\maketitle
\begin{abstract}
We prove global stability of Minkowski space for the Einstein vacuum
equations in harmonic (wave) coordinate gauge for the set of
restricted data coinciding with the Schwarzschild solution in the
neighborhood of space-like infinity. The result contradicts previous
beliefs that wave coordinates are "unstable in the large" and
provides an alternative approach to the stability problem originally
solved ( for unrestricted data, in a different gauge and with a
precise description of the asymptotic behavior at null infinity) by
D. Christodoulou and S. Klainerman. \vskip .3pc \noindent Using the
wave coordinate gauge we recast the Einstein equations as a system
of quasilinear wave equations and, in absence of the classical null
condition, establish a small data global existence result. In our
previous work we introduced the notion of a {\it  weak null condition}
and showed that the
Einstein equations in harmonic coordinates satisfy this condition.The
result of this paper relies on this observation and combines it with
the vector field method based on the symmetries of the standard
Minkowski space.

\vskip .5pc
\noindent
In a forthcoming paper we will address the question of stability of Minkowski
space for the Einstein vacuum equations in wave coordinates for all "small"
asymptotically flat data and the case of the Einstein equations coupled to a scalar
field.
\end{abstract}

\section {Introduction}
The focus of this paper is the question of global existence and
stability for the Einstein vacuum equations in "harmonic" (wave
coordinate) gauge. The Einstein equations determine a 4-d manifold
${\cal M}$ with a Lorentzian metric $g$ with vanishing Ricci
curvature
$$
R_{\mu\nu}=0.
$$
We consider the initial value problem: for a given a 3-d manifold $\Si$,
with a Riemannian metric $g_0$, and a symmetric two-tensor $k_0$, we
want to find a 4-d manifold ${\cal M}$, with a Lorentzian metric
$g$ satisfying the Einstein equations, and an imbedding
$\Sigma\subset M$ such that $g_0$ is the restriction of $g$ to
$\Sigma$ and $k_0$ is the second fundamental form of $\Sigma$ in
${\cal M}$. The initial value problem is overdetermined which
imposes compatibility conditions on the initial data: the
constraint equations
$$
R_0- {k_0}_{j}^i\, {k_0}^{j}_i + {k_0}_{i}^i \,{k_0}_j^j=0, \qquad
\nab^j {k_0}_{ij} - \nab_i \,{k_0}_j^j=0, \qquad \forall
i=1,...,3.
$$
Here $R_0$ is the scalar curvature of $g_0$ and $\nab$ is
covariant differentiation with respect to $g_0$. The Einstein
equations are invariant under diffeomorphisms. To have a working
formulation one needs to eliminate this freedom by fixing a gauge
condition or a system of coordinates.

 While the Einstein equations are
 independent of the choice of a coordinate system, the existence of
a special or preferred system of coordinates has been a subject of
debate \cite{Fo}. Historically, the first special coordinates were
the {\it harmonic} coordinates (also referred to as  {\it wave}
coordinates in current terminology). These obey the equation
$\Box_g\, x^\mu=0$, $\mu=0,1,2,3$, where $\Box_g=\nabla_\alpha
\nabla^\alpha$ is the geometric wave operator. Relative to the wave
coordinates a Lorentzian metric $g$ satisfies the {\it wave}
coordinate condition if:\footnote {We shall use below the standard
convention of summing over repeated indices and the notation
$\pa_\a=\pa/\pa x^\a$}
\begin{equation}\label{eq:wave-coord}
g^{\a\b}\pa_{\b} g_{\a\mu}=\frac 12 g^{\a\b}\pa_{\mu} g_{\a\b},
\qquad \forall \mu=0,..,3.
\end{equation}
In this system of coordinates, the vacuum Einstein equations
take the form of a system of quasilinear wave equations
\begin{equation}\label{eq:RE1}
g^{\a\b} \pa_\alpha\pa_\beta \, g_{\mu\nu} = {\cal N}_{\mu\nu} (g, \pa g), \qquad \forall \mu, \nu=0,..,3
\end{equation}
with a nonlinearity ${\cal N}(u,v)$ depending quadratically on $v$.
In this particular {\it gauge} Choquet-Bruhat \cite{CB1} was able to
establish the existence of a globally hyperbolic
development\footnote{For the definitions of global hyperbolicity and
maximal Cauchy development see \cite{H-E}, \cite{Wa} } of the
Einstein vacuum equations starting with an arbitrary set of initial
data prescribed on a 3-d space-like hypersurface and satisfying the
constraint equations. While the result of Choquet-Bruhat and a later
result of Choquet-Bruhat and Geroch \cite{CB-G}, establishing the
existence of a maximal Cauchy development, constructs solutions for
any given initial data set, it does not provide any information
about the geodesic completeness of the obtained solution. In the
language of the evolution equations these results only show the
existence of "local in time" solutions.

The global results have proved to be by far more resistant. The
outstanding global problem, which for a long time remained open, and
was finally ingeniously solved  by Christodoulou and Klainerman
\cite{C-K}, was that of the stability of Minkowski space. In
simplified language, it is the problem of constructing a global
solution to the Einstein vacuum equations from the initial data,
which is close to the Minkowski metric $m_{\mu\nu}$, and
asymptotically approaching the Minkowski space. The initial data
$(\Si, g_0, k_0)$ for the problem of stability of Minkowski space is
asymptotically flat, i.e., the complement of a compact set in $\Si $
is diffeomorphic to the complement of a ball in ${\Bbb R}^3$, and
there exists a system of coordinates $(x_1, x_2, x_3)$ with
$r=\sqrt{x_1^2 + x_2^2 +x_3^2}$ such that for all sufficiently large
$r$ the metric\footnote{The stability result of \cite{C-K} was
proved for {\it strongly} asymptotically flat data ${g_0}_{ij} =
(1+2M/r) \delta_{ij} + o(r^{-3/2})$,
 $k_0 = o(r^{-5/2})$.}
${g_0}_{ij} = (1+2M/r) \delta_{ij} + o(r^{-1-\si})$, and the second
fundamental form $k_0 = o(r^{-2-\si})$ for some $\si>0$. Here $M$ is
the mass, which by the positive mass theorem is positive unless the
data is flat, see Schoen and Yau \cite{S-Y} and Witten \cite{Wi}. In
addition, the data is required to satisfy a global smallness
assumption, which makes sure that it is sufficiently close to the
data $({\Bbb R}^3, \de, 0)$ for the Minkowski space.

To understand some of the difficulties of the problem we recall
that a generic system of quasilinear equations
\begin{equation}\label{eq:Quas}
\Box \phi_I = \sum_{|\a|\le|\b|\le 2}
A_{I,\a\b}^{JK} \pa^\a \phi_J \pa^\b\phi_K
+ {\text{cubic terms}}
\end{equation}
allows solutions with smooth arbitrarily small initial data which
blow up in finite time\footnote{This is in particular true for a
{\it semilinear} equation $\Box \phi = (\pa_{t}\phi)^{2}$,
\cite{J1}.}.  The key to global existence for such equations was
the {\it null} condition found by Klainerman,
\cite{K2}.
The small data global existence result for the equations satisfying the
null condition was established in \cite{C1}, \cite{K2}.
The null condition manifests itself in special
algebraic cancellations in the coefficients $ A_{I,\a\b}^{JK}$ of
the quadratic terms of the equation.\footnote{E.g. $\Box
\phi=(\pa_t \phi)^2-|\nabla_x \phi|^2$ satisfies the
null-condition.} It can be shown however, that the Einstein vacuum
equations in wave coordinates do not satisfy the null condition.
Moreover, Choquet-Bruhat \cite{CB3} showed that even without
imposing a specific gauge the Einstein equations violate the null
condition.

These considerations led to the suggestion that the wave coordinates
are not suitable for proving stability of Minkowski space. In fact,
considering a second iterate of the equation \eqref{eq:RE1},
Choquet-Bruhat \cite{CB2}  argued that the Einstein vacuum equations
are not stable in wave coordinates near the Minkowski solution. All
these resulted in the belief that the wave coordinates are unstable
in the large in the sense that a possible finite time blow up of
solutions of the equation \eqref{eq:RE1} is due to a coordinate
singularity.

The global stability of Minkowski space had been proved by
Christodoulou and Klainerman \cite{C-K} who avoided the use of a
preferred system of coordinates and instead relied on the
invariant formulation of the Einstein equations with the choice of
maximal time foliation (or the double null foliation in the new
proof of Klainerman and Nicolo \cite{K-N1}) and utilizing Bianchi
identities for the curvature. The special structure of the
quadratic terms plays a crucial part in the generalized energy
estimates which form the backbone of the proof but the null
condition can not be pointed out precisely.

A semiglobal stability result was also obtained in the work of
Friedrich \cite{Fr}. He used the conformal method
to reduce the global problem to a local one. The approach is
invariant and the special structure is again exploited implicitly.

In this paper we revisit the problem of global stability of Minkowski
space in wave coordinates. More precisely, we consider the data\footnote{The
existence of such data is guaranteed by the results of Corvino and Chrusciel-Delay,
\cite{Co}, \cite{C-D}.}
$({\Bbb R}^3, g_0, k_0)$ with the metric $g_0$ coinciding with the
spatial part of the Schwarzschild metric
$g_S= (1+M/r) ^4 dx^2$
in the region $r>1>>M$, vanishing second fundamental form $k_0$ for $r>1$,
and satisfying a global smallness assumption on ${\Bbb R}^3$.
We prove that for this initial data the wave coordinate gauge is stable in
the large: the reduced Einstein equations \eqref{eq:RE1} has a global
solution $g$ defining a future causally geodesically complete space-time,
\cite{H-E}.
The metric $g$ in wave coordinates $x^\a,\, \a=0,..,3$
approaches the Minkowski metric $m$:\,\,
$\sup_{x\in {\Bbb R}^3} |g(t,x)-m|\to 0$ as $t\to \infty$.

The intuition behind this result is based on the observation that
the Einstein vacuum equations in wave coordinates \eqref{eq:RE1}
satisfy the {\it weak null condition}. This notion was introduced
in \cite{L-R} for general quasilinear systems \eqref{eq:Quas} and
requires that the corresponding effective asymptotic system
\begin{equation}\label{eq:asymp}
(\pa_t + \pa_r)(\pa_t-\pa_r) \Phi_I ={r}^{-1} \sum_{n\le m\le 2}
A_{I,nm}^{JK} (\pa_t-\pa_r)^n \Phi_J \,\, (\pa_t-\pa_r)^m
\Phi_K,\qquad \Phi_I\sim r\phi_I
\end{equation}
has global solutions for all small initial data.\footnote{For the
precise definition see section \ref{weaknullcond}.} Here,
$$
A_{I,nm}^{JK}(\omega) = \sum_{|\a|=n, |\b|=m} A_{I,\a\b}^{JK}
\hat\omega^\a \hat \omega^\b,\qquad \hat \omega=(-1,\omega),\,\,
\omega\in {\Bbb S}^2.
$$
The classical null condition states that $A_{I,nm}^{JK}(\omega)
\equiv 0$ and thus implies the weak null condition. The asymptotic
system \eqref{eq:asymp} arises as an approximation of
\eqref{eq:Quas} when one neglects the derivatives tangential to the
outgoing Minkowski light cones, known to have faster decay.
The asymptotic equation was introduced in \cite{H1} to
predict the time of a blow-up for scalar wave equations known to
blow up in finite time, and was used in \cite{L2} to find some other scalar
wave equations for which the known blow-up mechanism was not
present. Asymptotic systems played an important role in the
analysis of the blow-up mechanisms in \cite{A1}.

In \cite{L-R} we have shown that the asymptotic system generated by the
Einstein equations in wave coordinates \eqref{eq:RE1} has global solutions
for all data. In this paper we consider the full nonlinear system
\eqref{eq:RE1}.
We should note that although the asymptotic system provides
useful heuristics about the behavior of solutions, in particular the
$L^\infty$ decay of the first derivatives of various components of
the metric $g$, it is barely used in our proof of the small data global existence
result for the full nonlinear equation  \eqref{eq:RE1}.  While it is tempting to
put forward a conjecture that, parallel to the result for the classical null
condition \cite{C1}, \cite{K2}, the weak null condition guarantees the
global existence result for small initial data, we can only argue that all
known examples seem to confirm it.
A simple example of an equation satisfying the weak null condition, violating
the standard null condition and yet possessing global solutions for all data
is given by the system
\begin{equation}\label{eq:simple}
\square \phi=w\cdot\pa^2\phi+\pa \psi\cdot\pa \psi,\qquad\quad
\Box\psi=0,\qquad\Box w=0
\end{equation}
 Another example is provided by the equation $\Box \phi=\phi
\Delta \phi$. The proof of a small data global existence result
for this equation is quite involved, \cite{L2} (radial case),
\cite{A3}. As we show in this paper the Einstein equations
\eqref{eq:RE1} is yet another example. Interestingly enough, at
the level of an effective asymptotic system the Einstein equations
can be modelled by the system \eqref{eq:simple}.

 The asymptotic behavior of null components
 of the Riemann curvature tensor $R_{\a\b\ga\de}$ of metric $g$\-- the so called
 "peeling estimates"\--  was discussed in the works of Bondi, Sachs and Penrose
 and becomes important in the framework of asymptotically
 simple space-times (roughly speaking, space-times which can be conformally
 compactified), see also the paper of Christodoulou \cite{C2} for further discussion
 of such space-times. Global solutions obtained in the work \cite{C-K} were
 accompanied by very precise analysis of its asymptotic behavior although not
 entirely consistent with peeling estimates. However, global solutions obtained
 by Klainerman-Nicolo \cite{K-N1} in the exterior\footnote{Outside of the domain
 of dependence of a compact set} stability of Minkowski space
 were shown to possess peeling estimates for special initial data, \cite{K-N2}.

 Our work is less precise about the asymptotic behavior and  is focused more on
 developing a technically relatively simple approach allowing us to prove stability
 of Minkowski space in a physically interesting wave coordinate gauge
In particular, we rely only on the standard Killing and conformal
Killing vector fields of Minkowski space and do not construct almost
Killing and conformal Killing vector fields adapted to the geometry
of null cones of the solution $g$.

Our proof is based on generalized energy estimates combined with
decay estimates. The  generalized energy estimates are used with
Minkowski vector fields $\{\pa_\a, \Omega_{\a\b}=x_\a\pa_\b-x_\b
\pa_\a, S=x^\a\pa_\a\}$. For the equations satisfying the standard
null condition uniform in time bounds on the generalized energies,
combined with global Sobolev (Klainerman-Sobolev) inequalities, are sufficient to infer
small data global existence. In our case however the generalized
energies slowly grow in time (at the rate of $t^\ve$) and need to
be complemented by independent, not following from the global
Sobolev inequalities, decay estimates. We derive the latter by
direct integration of the equation along the characteristics. It
is at this point that the intuition from the effective asymptotic
system is most useful. We show that all components of the metric
with exception of one decay at the rate of $t^{-1}$. The remaining
component however decays only as $t^{-1+\ve}$. Somewhat
surprisingly, the glue that holds together such weak decay
estimates and the generalized energy estimates is the wave coordinate
condition \eqref{eq:wave-coord}.

In this paper we only prove the result for a restricted set of data
coinciding with the Schwarzschild data outside of the ball of radius
one.\footnote{Since the initial metric is always of the form
$g_{ij}= (1+4M/r) \de_{ij} + o(r^{-1})$ with $M>0$, data coinciding
with the Schwarzschild outside of a compact set is the closest analogue
 of compactly supported or rapidly
decaying data usually considered in small data global existence results
for nonlinear wave equations.} This allows us to somewhat sidestep the problem of a
long range effect of a gravitational field. Due to the inward
bending of the light rays, solution arising from initial data
coinciding with the Schwarzschild data outside of the ball of radius
one will be equal to the Schwarzschild solution in the exterior of
the Minkowski cone $r=t+1$.

In our subsequent work we hope to be able to prove the stability of Minkowski
space in wave coordinates for general data.
In addition we hope to show that our method can be also used to treat the problem
of small data global existence for the Einstein
equations coupled to a scalar field.

{\bf Acknowledgments}:\,
The authors would like to thank  Demetrios Christodoulou and Sergiu Klainerman for
their inspiration and encouragement. We particularly benefited from
Sergiu Klainerman's suggestion to pursue first the problem with restricted data.
We would also like to thank Mihalis Dafermos
and Vince Moncrief  for stimulating discussions and useful suggestions.

\section{The main results and the strategy of the proof}
We now formulate the main results of our paper.
Our first result is global existence for the Einstein vacuum
equations in wave coordinates.
\begin{theorem}
Consider  the reduced Einstein vacuum equations \footnote{In what
follows we shall introduce the reduced wave operator
$\Boxr_g=g^{\a\b} \pa^2_{\a\b}$ and note that in wave coordinates
$\Boxr_g = \Box_g$, where $\Box_g \phi=|g|^{-1/2} \pa_\a
\big(g^{\a\b} |g|^{1/2}\pa_\b\phi\big)$ is the geometric wave
operator}
 \beq \label{eq:einst1}
\Boxr_g h_{\mu\nu}= g^{\a\b} \pa^2_{\a\b} h_{\mu\nu} = F_{\mu\nu}
(h)(\pa h, \pa h), \qquad \forall \mu, \nu=0,...,3,
 \eq
where $g_{\mu\nu}=m_{\mu\nu} + h_{\mu\nu}$ and the nonlinear
term $F$ is as in Lemma \ref{Einstwavecquad}.
We assume that the initial data $(g, \pa_t g)|_{t=0}=(g_0, g_1)$ are smooth,
the Lorentzian metric is of the form
$$
g_0= - a^2 dt^2 + {g_0}_{ij} dx^i dx^j
$$
 and

1) \,\, obey the wave coordinate condition
 \beq\label{eq:wavec10}
 g^{\a\a'}\pa_\a g_{\a'\mu} = \frac 12
g^{\a\a'} \pa_\mu g_{\a\a'}, \qquad \forall \mu=0,...,3,
 \eq

2)\,\, satisfy  the  constraint equations
$$
R_0 - |k_0|^2 + (\tr k_0)^2=0, \qquad
\nab^j {k_0}_{ij} - \nab_i \tr k_0=0, \quad \forall i=1,...,3,
$$
where $R_0$ is the scalar curvature of the metric ${g_0}_{ij}$, and
the second fundamental form $(k_0)_{ij}=-1/2 a^{-1} {g_{1}}_{ij}$.

3)\,\, we assume that the metric
$(g_0)_{ij}$ coincides with the spatial part of the Schwarzschild
metric $g_s$ (in wave coordinates):
$$
(g_0)_{ij}= \frac {r+2M}{r-2M} dr^2 + (r+ 2M)^2 (d\theta^2 + \sin^2\theta
\,d\phi^2), \qquad r>1
$$
and $g_1=0$ for $r>1$. Moreover, we assume that
the lapse function $a^2(r) = (r-2M)/(r+ 2M)$ for $r>1$ and
$a(r)=1$ for $r\le 1/2$

4)\,\, the data $(h_0, h_1)=(g_0-m, g_1)$ verify the
smallness condition
\begin{equation} \varepsilon=\sqrt{E_N(0)}
+M< \varepsilon_0,
\end{equation}
where $N\geq 10$ and
\begin{equation}
E_N(t)=\sup_{0\leq \tau\leq t}\sum_{|I|\leq N} \|\pa Z^I h(\tau,\cdot)\|_{L^2}^2
\end{equation}
Here $Z^I$ is a product of $|I|$ vector fields of the form
$\pa_i$, $x_i\pa_j-x_j\pa_i$, $t\pa_i+x_i\pa_t$ and
$t\pa_t+x^i\pa_i$.
Then there exists a unique global smooth solution $g$ with the property that
for some constant $C_N$,
\begin{align}
&E_N(t)\leq 16 \varepsilon^2 (1+t)^{2 C_N\varepsilon},\label{eq:strongenergyest0}\\
&\|g_{\mu\nu}(t) -m_{\mu\nu} \|_{L^\infty_x}\le C_N \varepsilon (1+t)^{-1+C_N\varepsilon}\nonumber.
\end{align}
\end{theorem}
\begin{remark}
The existence of data satisfying the assumptions of the theorem follows from
the work of \cite{Co}, \cite{C-D}, as argued in section 4.
\end{remark}
A corollary of the above result is the global stability of Minkowski space for
a restricted set of initial data.
\begin{theorem}
Let $({\Bbb R}^3, g_0, k_0)$ be the initial data set for the Einstein vacuum
equations $R_{\mu\nu}=0$. Assume that relative to some system of coordinates
$(x_1, x_2, x_3)$
the metric $g_0$ coincides with the spatial part of the
Schwarzschild metric $g_S$ outside the ball of radius one,
$$
g_0= (1+\frac {M}r)^4 dx^2, \qquad r>1,
$$
while the second fundamental form $k_0$ vanishes for $r>1$.
In addition, we assume that relative to that system of coordinates
$g_0$, $M$ and $k_0$ satisfy the smallness condition
$$
\sum_{0\le |I|\le N} \|\pa_x^I (g_0-\de)\|_{L^2(B_1)} +
\sum_{0\le |I|\le N-1}\|\pa_x^I k_0\|_{L^2(B_1)} +  M<\epsilon.
$$
Then there exists a future causally geodesically
complete\footnote{For the definition see \cite{H-E} and section 16
of this paper.} solution $g$ together with a global system of wave
coordinates with the property that the curvature tensor of $g$
relative to these coordinates decays to zero along any future
directed causal geodesic.
\end{theorem}
\noindent
 We now outline the strategy of
the proof.
\begin{remark}\label{rem:distinct}
Throughout the paper we shall use the notation
$A\les B$ for the inequality $A\le C B$ with some large
{\it universal} constant $C$. In our estimates we will make no
distinction between the tensors $h_{\a\b} = g_{\a\b} - m_{\a\b}$
and $H_{\a\b}=m_{\a\a'} m_{\b\b'} (g^{\a\b} - m^{\a\b})$,
since $H=-h + O(h^2)$ and the terms quadratic in $h$ are lower order.
\end{remark}
{\bf The continuity argument}
 For the proof we let $\delta$ be any fixed number
$0<\delta<1/2$. Let $g$ be a local smooth solution of the reduced
Einstein equations \eqref{eq:einst1}.
We start with the weak estimate
\beq
\label{eq:weakenergyest0}
E_N(t)\leq 64\varepsilon^2(1+t)^{2\delta}
\eq
By assumptions of the Theorem the estimate \eqref{eq:weakenergyest0}
holds for $t=0$. Let $[0,T]$ be
the largest time interval on which \eqref{eq:weakenergyest0} still holds.
 We shall show that if $\varepsilon>0$ is sufficiently
small then on the interval $[0,T]$ the inequality
\eqref{eq:weakenergyest0} implies the same inequality with the
constant $64$ replaced by $16$.  It will then follow that the
solution and the energy estimate \eqref{eq:weakenergyest0} can be
extended to a larger time interval $[0,T']$ thus contradicting the
maximality of $T$. This will imply that $T=\infty$ and the
solution is global. We will in fact prove  that for a sufficiently
small $\ve$ the stronger estimate \eqref{eq:strongenergyest0}
holds true on the interval $[0,T]$.

The global Sobolev inequality of Proposition
\ref{globalsobolev} and the
weak energy estimate \eqref{eq:weakenergyest0} imply the pointwise
decay estimates:
\begin{equation}
\sum_{|I|\leq N-2}|\pa Z^I h(t,x)|\leq
\frac{C\varepsilon(1+t)^\delta}{(1+t+r)(1+|t-r|)^{1/2}} ,
\qquad r=|x|
\end{equation}
From the assumption  that the constant  $\delta<1/2$
we derive the following weak decay estimates
\begin{equation} \label{eq:weakdecayest0}
|\pa Z^I h(t,x)|
\leq C\varepsilon(1+t+r)^{-1/2-\gamma}(1+|t-r|)^{-1/2-\gamma},
\qquad |I|\leq N-2
\end{equation}
with some fixed constant $\gamma>0$.
The weak decay estimates \eqref{eq:weakdecayest0}
will lead to
much stronger decay estimates in Theorem
\ref{decayeinst}.
In turn, using the stronger decay estimates in Theorem
\ref{decayeinst} we will be able to obtain
stronger energy estimates in Theorem
\ref{energyest}, i.e. \eqref{eq:strongenergyest0}.
These in particular will enable us to show that  the
estimate \eqref{eq:weakenergyest0} holds globally in time and  conclude the proof.
We remark that in the course of the proof all constants
will be independent of $\varepsilon>0$ but they will depend
on a lower bound for $\gamma>0$ (and hence on an upper bound for $\delta<1/2$).

As described above, the proof is a direct consequence of three
results. First is the global Sobolev inequality of
Proposition \ref{globalsobolev}, introduced by S. Klainerman \cite{K1}, giving decay
estimates in terms of energy estimates for the
generators of the Lorentz group. The second ingredient is the improved
decay estimates in Theorem \ref{decayeinst}. The final component  is
the energy
estimates in Theorem \ref{energyest} which rely on the improved decay
estimates.

{\bf Weak decay estimates.}
 As pointed out above we may start by
assuming the weak decay estimate \eqref{eq:weakdecayest0}.
Furthermore, since the solution $g=m+h$ coincides with the  Schwarzschild
solution of  mass $M\le \varepsilon$ in the region  $r\geq t+1$, we have \beq\label{eq:boundcond00}
 |Z^I h(t,x)|\les
\ve(1+r+t)^{-1},\qquad\text{when}\quad |x|=t+1 \eq
Hence
integrating \eqref{eq:weakdecayest0} from the light cone, where
\eqref{eq:boundcond00} holds, we get \beq\label{eq:weakdecay22}
 |Z^I h(t,x)|\les \ve(1+r+t)^{-1/2-\gamma}(1+|t-r|)^{1/2-\gamma},
 \eq
 Since the vector fields span the tangent space of the outgoing
light cones $r-t=q$ we infer, with $\pab$ denoting  the derivatives
tangential to the cones, that
 \beq\label{eq:weakdecaytangential}
 |\pab Z^I h|\les
\ve(1+r+t)^{-3/2-\gamma}(1+|t-r|)^{1/2-\gamma},
 \eq
This means that, close to the light cone $t=r$,  derivatives tangential to the forward light cones decay
 quite a bit better  than the expected decay rate from \eqref{eq:weakdecayest0} for a
generic derivative.

{\bf Wave coordinate condition.}
As we shall see below certain
components of the tensor $h$  decay faster than
others. This can be seen upon introduction of a null frame of vector fields
$L=\pa_t+\pa_r$, $\Lb=\pa_r-\pa_t$ and $S_1, S_2$:  two orthonormal vectors
tangential to the sphere of radius $r$  in $\bold{R}^3$. The first improved
estimates come from the wave coordinate condition
\eqref{eq:wavec10}. Writing
$g_{\alpha\beta}=m_{\alpha\beta}+h_{\alpha\beta}$ we obtain
from \eqref{eq:wavec10} that
$$
m^{\alpha\beta}\pa_\alpha h_{\beta\mu}=\pa_\mu  m^{\alpha\beta}
h_{\alpha\beta}+O( h\,\pa h)
$$
In particular, contracting with a vector field $T\in{\cal
T}=\{L,S_1,S_2\}$ and
using that for any symmetric 2-tensor $k$,
$m^{\alpha\beta}k_{\alpha\beta}=-k_{L\Lb}+
\delta^{AB} k_{AB}$, implies that we can express the
transversal derivative $\pa_{\Lb}$ of certain components of $h$
in terms of the tangential derivatives that decay better and a
quadratic term
$$
|(\pa h)_{LT}|\leq |\pab h|+ |h|\,|\pa h|\les \ve
(1+t+r)^{-1-2\gamma},\qquad |h_{LT}|\les \ve (1+|t-r|)(1+t+r)^{-1}
$$
Even though the estimate above does not give a better decay rate for all
components of $h$ it gives the decay exactly for those components which, as
it turns out, control the geometry, i.e., they lead to stronger energy and decay
estimates.

The above estimates will be sufficient to obtain improved estimates for the lowest
order energy of $h$. However, in order to get estimates for the energy of
$Z^I h$ we commute the vector fields $Z$ through
the equation for $h$. This generates additional  commutator terms.
The main commutator terms are controlled with the help of the
following additional estimate from the wave coordinate condition:
 \beq\label{eq:wavecdecay}
 |(\pa h)_{LT}|+|(\pa Z h)_{LL}|\leq \ve (1+t+r)^{-1-2\gamma},
 \qquad |h_{LT}|+|(Z h)_{LL}|\leq \ve (1+|t-r|)(1+t+r)^{-1}
 \eq
 We now describe derivation of  the stronger decay and energy estimates.

 {\bf Stronger decay estimates.}
  We rely on the following decay
estimate for  the wave equation on a curved background
\footnote{Recall that the reduced wave operator $\Boxr_g= g^{\a\b}
\pa^2_{\a\b}$.}:
\begin{multline}\label{eq:strongdecaywaveeq}
\|(1+t+r)\pa \phi(t,\cdot)\|_{L^\infty} \leq  C\int_0^t (1+\tau)\|
\Boxr_g\phi(\tau,\cdot)\|_{L^\infty}\,d\tau \\ +
 C\sup_{0\leq
\tau\leq t}\! \sum_{|I|\leq 1}\! \|Z^I\!
\phi(\tau,\cdot)\|_{L^\infty} + C\int_0^t\sum_{|I|\leq 2}
(1+\tau)^{-1} \| Z^I \phi(\tau,\cdot)\|_{L^\infty}\, d\tau
\end{multline}
The estimate \eqref{eq:strongdecaywaveeq} will be applied to the
components of the tensor $h$. The term $Z^I h$ on the right hand
side of the estimate will be controlled with the help of the weak
decay estimates, and   thus the decay rate of  $h$ will be
determined in terms of decay of $\Boxr_g h$. The
$L^\infty-L^\infty$ estimate \eqref{eq:strongdecaywaveeq} does not
rely on the fundamental solution as does the more standard
$L^1-L^\infty$ type estimate. This estimate was used \cite{L1} in
the constant coefficient case and here we establish it in the
variable coefficient case only under the assumption of the weak
decay of all of the components of the metric $g$ and the stronger
decay of the components of $g$ controlled by  the wave coordinate
condition. This analysis is by itself very interesting but we will
not go into it here and just refer the reader to the following
sections.

 We now analyze the inhomogeneous term in the equation for $h_{\mu\nu}$.
 The tensor
$h_{\mu\nu}= g_{\mu\nu}- m_{\mu\nu}$ verifies the reduced Einstein
equations of the form:
\begin{align}
&\Boxr_g h_{\mu\nu} = F_{\mu\nu} (h) (\pa h, \pa h),\nn\\
& F_{\mu\nu}(h)(\pa h,\pa h)=P(\pa_\mu h,\pa_\nu h)+Q_{\mu\nu}(\pa
h,\pa h)
+G_{\mu\nu}(h)(\pa h,\pa h)\label{eq:F},\\
& P(\pa_\mu h,\pa_\nu h)= \frac 14 \pa_{\mu} \tr h\, \pa_{\nu}\tr
h-\frac 12 \pa_{\mu} h^{\a\b}\pa_{\nu} h_{\a\b},\label{eq:P}
\end{align}
Here $Q_{\mu\nu}$ are linear combinations of the standard null-forms and
$G_{\mu\nu}(h)(\pa h,\pa h)$ is a quadratic form in $\pa h$ with
coefficients as smooth functions of $h$ vanishing at
$h=0$. The weak decay estimates imply that the last two
terms decay fast
\beq
 \label{eq:nullformtan0}
|Q_{\mu\nu}(\pa h,\pa h)|+|G_{\mu\nu}(h)(\pa h,\pa h)|\les
|\overline{\pa} h |\,|\pa h| + |h||\pa h|^2\les \varepsilon^2
(1+r+t)^{-2-2\ga} (1+|t-r|)^{-2\ga}
 \eq
 The problematic term is $P(\pa_\mu h,\pa_\nu h)$ since a priori the weak decay estimates only give  the decay rate  of
 $\varepsilon^2 (1+r+t)^{-1-2\ga} (1+|t-r|)^{-1-2\ga}$, which is not
 sufficient in the wave zone $t\approx r$.
 The crucial improvement comes as a result of a decomposition
 of the tensor  $P(\pa_\mu h,\pa_\nu h)$ with respect to a null
 frame  $\{L, \Lb, S_1, S_2\}$. Let $T\in {\cal T} =\{L, S_1, S_2\}$ be
 any of the vectors generating the tangent space to the forward
 Minkowski light cones and $U \in {\cal U} =\{L, \Lb, S_1, S_2\}$
 denote any of the null frame vectors.
 Define, for an arbitrary symmetric two tensor $k$,
 $|k|_{\cal TU}
 =\sum_{T\in {\cal T}, U\in {\cal U}}=|T^\mu U^\mu k_{\mu\nu}|$.
It then follows that
 \beq
|P(\pa h,\pa h)|_{\cal TU} \les |\overline{\pa} h|\,|\pa h|
\les \varepsilon^2
(1+r+t)^{-2-2\ga} (1+|t-r|)^{-2\ga}
 \eq
 On the other hand, the absolute value of the tensor $P(\pa h,\pa h)$ obeys the
 estimate
  \beq
  |P(\pa h,\pa h)|\les |\pa h|_{\cal TU}^2
+|\pa h|_{LL}|\pa h| \eq
We now decompose the system of equations for $h$ with respect  to the
null-frame
 \begin{align}
|\Boxr_g h|_{\cal TU}&\les \varepsilon^2
(1+r+t)^{-2-2\ga} (1+|t-r|)^{-2\ga} ,\label{eq:one}\\
|\Boxr_g h|_{\cal UU}&\les |\pa h|_{\cal TU}^2+ \varepsilon^2
(1+r+t)^{-2-2\ga} (1+|t-r|)^{-2\ga}\label{eq:two},
 \end{align}
where in the last inequality we also used the improved decay
estimate for $\pa h_{LL}$ obtained from the wave coordinate
condition. The result is a system of equations where the
components $\Boxr_g h_{TU}$ have very good decay properties, while
 $\Boxr_g h_{UU}$ for the remaining non tangential component
 depends,  to the highest order, only on the components $h_{TU}$ satisfying a better equation.
An additional subtlety in  the above analysis is the fact that
contraction with the null frame does not commute with $\Boxr_g$
(or even with $\square$). However, the decay estimate
\eqref{eq:strongdecaywaveeq} for the wave equation only uses the
principal radial part of $\square$:\, $\pa_t^2 - r^{-2}
\pa_r^2-2\, r^{-1} \pa_r $, which respects the null frame. This
analysis results in the improved decay estimates
\begin{equation}\label{eq:decay-first}
|\pa h|_{TU}\le C\ve (1+t)^{-1}, \qquad
|\pa h|\le C\ve (1+t)^{-1} \ln (2+t)
\end{equation}

{\bf The energy estimates.} We rely on the following energy
estimate for the wave equation, which holds under the assumption
that  the above decay estimates hold for the background metric
$g$:  for any  $\gamma>0$ \beq\label{eq:basicenergy}
 \int_{\Si_{T}} |\pa\phi|^{2} + \int_{0}^{T} \int_{\Si_{\tau}}
\frac {\gamma \,\, |\pab\phi|^{2}}{(1+|t-r|)^{1+2\gamma}}\leq
 8\int_{\Si_{0}} |\pa \phi|^{2} +  C\varepsilon
\int_0^T\int_{\Si_{t}} \frac {|\pa\phi|^{2}}{1+t}
 +16\int_0^T  \int_{\Sigma_t} |\Boxr_g\phi||\pa_t\phi|
\eq
 This implies that the energy of a solution of the homogeneous
 wave equation $\Boxr_g \phi=0$ grows but
 at the rate of at most $(1+t)^{C\ve}$.
 The presence of an additional space-time integral containing tangential derivatives
 on the right and side of  \eqref{eq:basicenergy} is crucial for our analysis.
 This type of estimate in the constant coefficient case basically follows
 by averaging of the energy estimates on light cones used e.g. in \cite{S1}.
 We also note that the energy estimates with space-time quantities involving
 special derivatives of a solution were also considered and used in the work of Alinhac, see
 e.g. \cite{A2}, \cite{A3}). In our work we use the space-time integral
 with derivatives spanning the tangent space to outgoing light cones and weights
 dependent on the distance to the Minkowski light cone $r=t+1$.
 We emphasize that the
 energy estimate \eqref{eq:basicenergy} is proved  only under the assumption of the
 weak decay of all components of the background metric $g$ together  with
 the strong decay of the components controlled from the wave coordinate condition.

It is worth noting that a combination of the energy estimates of the type 
\eqref{eq:basicenergy} and Klainerman-Sobolev inequalities would also yield a very simple
proof of the small data global existence result for semilinear equations 
$\Box \phi = Q(\pa\phi,\pa\phi)$ obeying the standard null condition. 
This fact appears to be previously unknown.

The energy estimate \eqref{eq:basicenergy} will be applied simultaneously
to all components of the tensor $h$.
As in the equations \eqref{eq:one}, \eqref{eq:two} the inhomogeneous term obeys the following estimate:
$$
|\Boxr_g h|\les \ve (1+r+t)^{-3/2-\ga} (1+|t-r|)^{1/2-\ga} |\pa
h|+ \ve (1+t)^{-1} |\pa h|,
$$
where in the last inequality we used the improved decay estimate for
the $|\pa h|_{\cal TU}$ components.
The energy estimate \eqref{eq:firstenergy} will then imply that
$$
E_0(t)\le 16\ve^2 (1+t)^{C\ve}.
$$

{\bf Higher order energy estimates.}
 In addition to the energy
estimates for the components of the tensor $h$ we need estimates
for the higher vector field derivatives of $h$: $Z^I h$ with
Minkowski vector fields $Z=\{\pa_\a, \Omega_{\a\b}, S\}$. To
obtain these estimates we apply $Z^I$ to the equation $\Boxr_{g}
h_{\mu\nu} = F_{\mu\nu}$ for $h$. Applying vector fields to the
nonlinear terms $F_{\mu\nu}$ yields similar nonlinear terms for
higher derivatives and these are can be dealt with using the
estimates already described above. We must note however that this
is where the additional space-time integral involving the
tangential derivatives on the left hand side of the energy
estimate \eqref{eq:firstenergy} becomes crucial. Consider for
example the term $\pa h\cdot\pab Z^I h$ generated by one of the
null forms in $F_{\mu\nu}$. We estimate its contribution, with the
help of the weak decay estimates, to the energy estimate as
follows:
$$
|\pa h|\, |\pab Z^I h||\pa_t Z^I h|\leq
\frac{C\ve|\pa_t Z^I h|}{(1+t)^{1/2+\gamma}}
\frac{|\pab Z^I h|}{(1+|t-r|)^{1/2+\gamma}}
\leq \frac{C \ve|\pa_t Z^I h|^2}{(1+t)^{1+2\gamma}}
+\frac{C\ve |\pab Z^I h|^2}{(1+|t-r|)^{1+2\gamma}}
$$
The integral of the first term is easily controlled by the
energy on time slices times an integrable factor in time.
The space time integral of the second term is in fact part of the energy
\eqref{eq:basicenergy}, and if we choose $\ve$ sufficiently small this
term can be absorbed by the space time integral on the left.
The idea with the space-time integral is that one can use the extra
decay in $|t-r|$ when one does not have full decay in $t$.

The more serious problem in higher order energy estimates lies
however in the commutators between $Z^I$ and the principal part
$\Boxr_g=g^{\alpha\beta}\pa_\alpha\pa_\beta$.

{\bf The commutators.} Writing
$g^{\alpha\beta}=m^{\alpha\beta}+H^{\alpha\beta}$ with $H^{\a\b} =
-m^{\a\a'}m^{\b\b'} h_{\a'\b'} + O(h^2)$, we show the following
commutator estimate\footnote{This commutator estimate applies to
the vector fields $Z=\{\pa_\a, \Omega_{\a\b}\}$. For the scaling
vector field $Z=S=x^\a\pa_\a$ the commutator expression should
have the form $\Boxr_g S - (S+2) \Boxr_g$ .}
 \beq\label{eq:commutatorestimate}
\big| [Z,\Boxr_g]\phi \big| \leq C\Big(\frac{|H|+|Z
H|}{1+t+r}+\frac{|Z H|_{LL}+|H|_{L\cal
T}}{1+|\,t-r|}\Big)\sum_{|I|\leq 1}|\pa Z^I \phi|\leq
\frac{C\ve}{1+t+r}\sum_{|I|\leq 1}|\pa Z^I \phi|
 \eq
 by the weak decay assumptions \eqref{eq:weakdecay22} and the improved decay
 from the wave coordinate condition
 \eqref{eq:wavecdecay}.  We should note that for a generic quasilinear wave
 equation commutators with Minkowski vector fields $Z$ give rise to
 uncontrollable error terms. In the special case of the equation
 $\Box\phi =\phi \Delta \phi$
 this problem can be overcome by modifying the vector fields $Z$, \cite{A3}.
 In our case it is the wave coordinate gauge that provides additional cancellations.

  This commutator estimate applied to $\phi=h_{\alpha\beta}$
 together with the analysis in the previous section now gives estimates for the
 energy  $E_1$ as well as for the stronger decay estimates for the
second derivatives of $h$, \eqref{eq:strongdecayind} with $|J|=1$.
This commutator
 will also show up as a top order term $[Z,\Boxr_g] \cdot Z^{I-1} h_{\alpha\beta}$
 in the energy estimate for $Z^I h$
 and the resulting term can be dealt with in the same way.

 The other top order term generated by the commutators $[Z^I,\Boxr_g]\phi$
 is of the form $(Z^I H^{\alpha\beta} ) \pa_\alpha\pa_\beta $.
 We first apply the pointwise estimate
 $$
 |(Z^I H^{\alpha\beta} ) \pa_\alpha\pa_\beta \phi|\leq
 C\Big(\frac{|Z^I H|}{1+t+r}
 +\frac{ |Z^I H|_{LL}}{1+|\, t-r|}\Big)
 \sum_{|K|\leq 1} |\pa Z^K \phi|
 $$
 To deal with its contribution to the energy estimate
 we use the Poincare estimate with a boundary term
\beq \label{eq:Poincstrategy} \il_{R^{3}} \frac {|Z^I
H|_{LL}^{2}\, dx}{(1+|t-r|)^{2+2\si}} \le C\il_{S_{(t+1)}} |Z^I
H|_{LL}^{2} \, dS+ C\il_{R^{3}} \frac {|\pa_{r} Z^I H|_{LL}^{2}\,
dx}{ (1+|t-r|)^{2\si}},\qquad \si>-1/2, \quad \si\neq 1/2 \ee
together with the fact that $h$ is Schwarzschild outside the cone
$r=t+1$, because of the inward bending of the Schwarzschild light
cones, and hence there $|Z^I h|\leq C\ve /(1+t)$. The way
coordinate condition implies that $|\pa Z^I H|_{LL}$ can be
estimated by $|\pab Z^I H|$ and lower order terms. The term
involving $|\pab Z^I H|$ is then controlled by the space-time
integral on the left hand side.

One can use a similar but more trivial argument
for decay estimates, i.e.
$$
|Z^I H|_{LL}\leq  {|Z^I H_{LL}|}_{r=t+1}+(1+|t-r|){|\pa_r Z^I H_{LL}|}_{L^\infty}
$$

{\bf The lower order terms.} So far we have only discussed the top
order terms, but there will also be several lower order terms
(relative to $|I|=k+1$) to deal with. These are typically of the
form \beq\label{eq:lowerorderterms} |\pa Z^J h|\, |\pa Z^K
h|\qquad \text{or}\qquad |Z^J h|\, |\pa^2 Z^{K-1} h|\leq
C\frac{|Z^J h|}{1+|t-r|}\, |\pa Z^{K} h| \eq with $|J|,
|K|<|I|=k+1$ The lower order terms are dealt with using induction.
We describe the induction argument for the decay estimates. From
this it will be clear how it also proceeds for the energy
estimates. We will inductively assume that we have the bounds:
\beq\label{eq:strongdecayind} |\pa Z^J h|+| Z^J h|(1+|t-r|)^{-1}
\leq C_k t^{-1+C_k\ve},\qquad |J|\leq k \eq The terms in
\eqref{eq:lowerorderterms} can then be estimated by $C_k^2 \ve^2
t^{-2+2C_k\ve}$. Including the top order terms using
\eqref{eq:commutatorestimate} applied to $\phi=Z^{I-1} h$, and
using \eqref{eq:strongdecaywaveeq} applied to $\Boxr_g Z^I h$ we
get an inequality of the form \beq\label{eq:gronwall} M(t)\leq
\int_0^t \frac{C\ve \, M(s)}{1+s}\, + \frac{C\ve
^2}{(1+s)^{1-C\ve}} \, ds \eq where $M(t)=(1+t)\| \pa Z^I
h(t,\cdot)\|_{L^\infty}$. The Gronwall's inequality then  gives
the bound $M(t)\leq C(1+t)^{2C\ve}$.

\section {The Einstein equations in wave coordinates}\label{section:einstwavec}

For a Lorentzian metric  $g_{\mu\nu}$, where $\mu,\nu=0,...,3$
we denote
\beq
\label{eq:christof}
\Gamma_{\mu\,\,\, \nu}^{\,\,\, \lambda} =\frac{1}{2}
g^{\lambda\delta} \big( \pa_\mu g_{\delta\nu} +
\pa_\nu g_{\delta\mu}-\pa_\delta g_{\mu\nu}\big),
\eq
the Christoffel symbols of $g$ and
\beq
\label{eq:ctensor}
R_{\mu\,\,\, \nu\delta }^{\,\,\, \lambda}=
\pa_\delta \Gamma_{\mu\,\,\, \nu}^{\,\,\, \lambda}
-\pa_\nu\Gamma_{\mu\,\,\, \delta}^{\,\,\, \lambda}
+\Gamma_{\rho\,\,\, \delta}^{\,\,\, \lambda} \Gamma_{\mu\,\,\, \nu}^{\,\,\, \rho}
-\Gamma_{\rho\,\,\, \nu}^{\,\,\, \lambda} \Gamma_{\mu\,\,\, \delta}^{\,\,\, \rho}
\eq
its Riemann curvature tensor with
$R_{\mu\nu}=R_{\mu\,\,\, \nu\alpha}^{\,\,\, \alpha}$,
the Ricci tensor.

We consider the metric $g$ satisfying the Einstein vacuum equations
\beq
\label{eq:Einst}
R_{\mu\nu}=0.
\eq
We impose the wave coordinate condition:
\beq
\label{eq:wavecord}
 \Ga^\lambda:=g^{\alpha\beta}\, \Gamma_{\alpha\,\,\, \beta}^{\,\,\, \lambda}=0
\eq
It follows that assuming \eqref{eq:wavecord} we have that
the reduced wave operator $\Boxr_g=g^{\a\b}$
\beq
\label{eq:boxg}
\Boxr_g={\square}_g=\frac{1}{\sqrt{|g|}}\pa_\alpha
 g^{\alpha\beta}\sqrt{|g|}\pa_\beta
\eq
The following lemma provides the description of the Einstein vacuum
equations in wave coordinates as a system of quasilinear wave equations
for $g_{\mu\nu}$.
\begin{lemma} \label{Einstwavec}
Let metric $g$ satisfy the Einstein vacuum equations \eqref{eq:Einst} together with
the wave coordinate condition \eqref{eq:wavecord}.
Then $g_{\mu\nu}$ solves the following system of {\it reduced} Einstein equations:
\beq
\label{eq:Einstred} {\Boxr}_g  g_{\mu\nu}=\widetilde{P}(\pa_\mu
g,\pa_\nu g) +\widetilde{Q}_{\mu\nu}(\pa g,\pa g) \eq where
\begin{align}
&\widetilde{P}(\pa_\mu g,\pa_\nu g)=
\frac{1}{4}  g^{\alpha\alpha^\prime}\pa_\mu
g_{\alpha\alpha^\prime} \,  g^{\beta\beta^\prime}\pa_\nu
g_{\beta\beta^\prime}- \frac{1}{2} g^{\alpha\alpha^\prime}g^{\beta\beta^\prime} \pa_\mu
g_{\alpha\beta}\, \pa_\nu g_{\alpha^\prime\beta^\prime}
\label{eq:deftildeM} \\
&
\widetilde{Q}_{\mu\nu}(\pa g,\pa g)= \pa_{\alpha} g_{\beta\mu}\,
 \,g^{\alpha\alpha^\prime}g^{\beta\beta^\prime}
 \pa_{\alpha^\prime} g_{\beta^\prime\nu}
-g^{\alpha\alpha^\prime}g^{\beta\beta^\prime} \big(\pa_{\alpha}
g_{\beta\mu}\,\,\pa_{\beta^\prime} g_{\alpha^\prime \nu}
-\pa_{\beta^\prime} g_{\beta\mu}\,\,\pa_{\alpha}
g_{\alpha^\prime\nu}\big)
\label{eq:tildenullform}\\
&\qquad +g^{\a\a'}g^{\b\b'}\big (\pa_\mu g_{\a'\b'}
\pa_\a g_{\b\nu}- \pa_\a g_{\a'\b'}
\pa_\mu g_{\b\nu}\big )  + g^{\a\a'}g^{\b\b'}\big (\pa_\nu g_{\a'\b'} \pa_\a g_{\b\mu} - \pa_\a
g_{\a'\b'} \pa_\nu g_{\b\mu}\big )\nn\\
&\qquad +\frac 12 g^{\a\a'}g^{\b\b'}\big (\pa_{\b'} g_{\a\a'}
\pa_\mu g_{\b\nu} - \pa_{\mu} g_{\a\a'}
\pa_{\b'} g_{\b\nu} \big ) +\frac 12 g^{\a\a'}g^{\b\b'}
\big (\pa_{\b'} g_{\a\a'} \pa_\nu g_{\b\mu} -
\pa_{\nu} g_{\a\a'} \pa_{\b'} g_{\b\mu} \big )\nn
\end{align}
Furthermore, the wave coordinate condition  \eqref{eq:wavecord}  reads
\beq
\label{eq:wavec2}
g^{\alpha\beta}\pa_\alpha g_{\beta\mu}
=\frac{1}{2}g^{\alpha\beta}\pa_\mu g_{\alpha\beta},\qquad\text{or}\qquad
\pa_\alpha g^{\alpha\nu}=\frac{1}{2} g_{\alpha\beta} g^{\nu\mu}
 \pa_\mu g^{\alpha\beta}
\eq
\end{lemma}
\begin{proof}
The proof of \eqref{eq:wavec2} is immediate.

We now observe that
$$
\pa_\alpha g_{\beta\mu}=\Gamma_{\alpha\beta\mu}+\Gamma_{\alpha\mu\beta},
\qquad\quad\text{where}\qquad
\Gamma_{\mu\alpha\nu}=g_{\alpha\lambda}
\Gamma_{\mu\,\,\, \nu}^{\,\,\, \lambda}.
$$
It follows that $g_{\alpha\lambda} \pa_\beta \Gamma_{\mu\,\,\, \nu}^{\,\,\, \lambda}
=\pa_\beta \Gamma_{\mu\alpha\nu}-
(\Gamma_{\beta\alpha\lambda}+\Gamma_{\beta\lambda\alpha})
\Gamma_{\mu\,\,\, \nu}^{\,\,\, \lambda}$ so also using that
$\Gamma_{\alpha\lambda\beta}=\Gamma_{\beta\lambda\alpha}$ we obtain
\beq
\label{eq:Rexpr}
R_{\mu\alpha\nu\beta}=g_{\alpha\lambda}
R_{\mu\,\,\, \nu\beta }^{\,\,\, \lambda}=
\pa_\beta \Gamma_{\mu\alpha \nu}
-\pa_\nu\Gamma_{\mu\alpha\beta}
+\Gamma_{\nu\lambda\alpha} \Gamma_{\mu\,\,\, \beta}^{\,\,\, \lambda}
-\Gamma_{\alpha\lambda\beta} \Gamma_{\mu\,\,\, \nu}^{\,\,\, \lambda}
\eq
 It follows from \eqref{eq:wavec2} that
 \beq
g^{\alpha\beta}\Big(\pa_\mu \pa_\alpha g_{\beta\nu}-\frac{1}{2}
\pa_\mu\pa_\nu g_{\alpha\beta}\Big)=- \pa_\mu
g^{\alpha\beta}\Big(\pa_\alpha g_{\beta\nu}-\frac{1}{2}\pa_\nu
g_{\alpha\beta}\Big)=
g^{\alpha\alpha^\prime}g^{\beta\beta^\prime}\,\pa_\mu
g_{\alpha^\prime\beta^\prime}\Big(\pa_\alpha
g_{\beta\nu}-\frac{1}{2}\pa_\nu g_{\alpha\beta}\Big)
 \eq
 and hence
\begin{multline}
\label{eq:calc} g^{\alpha\beta} \big(\pa_\alpha
\Gamma_{\mu\beta\nu} -\pa_\nu \Gamma_{\mu\beta\alpha}\big)=
\frac{g^{\alpha\beta}}{2}  \big( \pa_\alpha\pa_\mu g_{\beta\nu} +
\pa_\alpha\pa_\nu g_{\beta\mu}-\pa_\alpha\pa_\beta g_{\mu\nu}\big)
-\frac{g^{\alpha\beta}}{2} \big( \pa_\nu\pa_\mu g_{\beta\alpha} +
\pa_\nu\pa_\alpha g_{\beta\mu}-\pa_\nu\pa_\beta g_{\mu\alpha}\big)\\
=-\frac{g^{\alpha\beta}}{2}\pa_\alpha\pa_\beta g_{\mu\nu} +
\frac{g^{\alpha\beta}}{2}\big(\pa_\alpha\pa_\mu g_{\beta\nu}
+\pa_\nu\pa_\beta g_{\mu\alpha} -\pa_\nu\pa_\mu
g_{\beta\alpha}\big)\\
=-\frac{1}{2}\, g^{\alpha\beta}\pa_\alpha\pa_\beta g_{\mu\nu}
+\frac{1}{2}g^{\alpha\alpha^\prime}g^{\beta\beta^\prime}\,\big
(\pa_\mu g_{\a^\prime\b^\prime} \pa_\a g_{\b\nu} +\pa_\nu
g_{\a^\prime\b^\prime} \pa_\a g_{\b\mu} -\pa_\nu
g_{\a^\prime\b^\prime} \pa_\mu g_{\a\b}\big) .
\end{multline}
Here by \eqref{eq:wavec2} we can write
\begin{multline}\label{eq:metricproduct}
g^{\alpha\alpha^\prime}g^{\beta\beta^\prime}\,\pa_\mu
g_{\a^\prime\b^\prime} \pa_\a g_{\b\nu}=
g^{\alpha\alpha^\prime}g^{\beta\beta^\prime}\,\pa_\a
g_{\a^\prime\b^\prime} \pa_\mu g_{\b\nu}+
g^{\alpha\alpha^\prime}g^{\beta\beta^\prime}\,\big(\pa_\mu
g_{\a^\prime\b^\prime} \pa_\a g_{\b\nu}- \pa_\a
g_{\a^\prime\b^\prime} \pa_\mu g_{\b\nu}\big)\\
= \frac{1}{2}
g^{\alpha\alpha^\prime}g^{\beta\beta^\prime}\,\pa_{\beta^\prime}
g_{\a^\prime\alpha} \,\pa_\mu g_{\b\nu}+
g^{\alpha\alpha^\prime}g^{\beta\beta^\prime}\,\big(\pa_\mu
g_{\a^\prime\b^\prime} \,\pa_\a g_{\b\nu}- \pa_\a
g_{\a^\prime\b^\prime}\,\pa_\mu g_{\b\nu}\big)\\
= \frac{1}{2}
g^{\alpha\alpha^\prime}g^{\beta\beta^\prime}\,\pa_{\mu}
g_{\a^\prime\alpha} \,\pa_{\beta^\prime}g_{\b\nu}+
g^{\alpha\alpha^\prime}g^{\beta\beta^\prime}\Big(\frac{1}{2}\big(
\pa_{\beta^\prime} g_{\a^\prime\alpha} \,\pa_\mu g_{\b\nu}
-\pa_{\mu} g_{\a^\prime\alpha} \,\pa_{\beta^\prime}g_{\b\nu}\big)
+ \big(\pa_\mu g_{\a^\prime\b^\prime}\, \pa_\a g_{\b\nu}- \pa_\a
g_{\a^\prime\b^\prime}\, \pa_\mu g_{\b\nu}\big)\Big)\\
=\frac{1}{4}
g^{\alpha\alpha^\prime}g^{\beta\beta^\prime}\,\pa_{\mu}
g_{\a^\prime\alpha} \,\pa_{\nu}g_{\b\beta^\prime}+
g^{\alpha\alpha^\prime}g^{\beta\beta^\prime}\Big(\frac{1}{2}\big(
\pa_{\beta^\prime} g_{\a^\prime\alpha} \,\pa_\mu g_{\b\nu}
-\pa_{\mu} g_{\a^\prime\alpha} \,\pa_{\beta^\prime}g_{\b\nu}\big)
+ \big(\pa_\mu g_{\a^\prime\b^\prime}\, \pa_\a g_{\b\nu}- \pa_\a
g_{\a^\prime\b^\prime}\, \pa_\mu g_{\b\nu}\big)\Big)
\end{multline}
Hence by  \eqref{eq:metricproduct}
 and \eqref{eq:metricproduct} with $\mu$ and $\nu$ interchanged
 we get
\begin{multline}\label{eq:metricproducttwo}
\frac {1}{2}g^{\alpha\alpha^\prime}g^{\beta\beta^\prime}\, \Big
(\pa_\mu g_{\alpha^\prime\beta^\prime}\, \pa_\a g_{\b\nu} +\pa_\nu
g_{\alpha^\prime\beta^\prime} \,\pa_\a g_{\b\mu} -\pa_\nu
g_{\alpha^\prime\beta^\prime}\, \pa_\mu g_{\a\b}\Big)=
g^{\alpha\alpha^\prime}g^{\beta\beta^\prime}\, \big(\frac{1}{4}
\pa_{\mu} g_{\a^\prime\alpha} \,\pa_{\nu}g_{\b\beta^\prime}
 -\frac{1}{2}\pa_\nu
g_{\alpha^\prime\beta^\prime}\, \pa_\mu g_{\a\b}\big)\\
+\frac{1}{2}g^{\alpha\alpha^\prime}g^{\beta\beta^\prime}\,\Big(
\big(\pa_\mu g_{\a^\prime\b^\prime}\, \pa_\a g_{\b\nu}- \pa_\a
g_{\a^\prime\b^\prime}\, \pa_\mu g_{\b\nu}\big) +\big(\pa_\nu
g_{\a^\prime\b^\prime}\, \pa_\a g_{\b\mu}- \pa_\a
g_{\a^\prime\b^\prime}\, \pa_\nu g_{\b\mu}\big) \Big)\\
\frac{1}{4}g^{\alpha\alpha^\prime}g^{\beta\beta^\prime} \Big(\big(
\pa_{\beta^\prime} g_{\a^\prime\alpha} \,\pa_\mu g_{\b\nu}
-\pa_{\mu} g_{\a^\prime\alpha} \,\pa_{\beta^\prime}g_{\b\nu}\big)
+\big( \pa_{\beta^\prime} g_{\a^\prime\alpha} \,\pa_\nu g_{\b\mu}
-\pa_{\nu} g_{\a^\prime\alpha}
\,\pa_{\beta^\prime}g_{\b\mu}\big)\Big)
\end{multline}
 On the other hand
\begin{multline}\label{eq:gammaproduct}
\Gamma_{\nu\alpha\beta} \Gamma_{\mu}^{\,\,\, \alpha\beta}
=\frac{1}{4} \big( \pa_\nu g_{\beta\alpha} +
\pa_\beta g_{\alpha\nu}-\pa_\alpha g_{\beta\nu}\big)
g^{\alpha\alpha^\prime}g^{\beta\beta^\prime}
 \big( \pa_{\mu} g_{\beta^\prime\alpha^\prime} +
\pa_{\beta^\prime} g_{\alpha^\prime\mu}-\pa_{\alpha^\prime}
 g_{\beta^\prime\mu}\big)\\
 =\frac{1}{4}
\pa_\nu g_{\alpha\beta}\,\,
g^{\alpha\alpha^\prime}g^{\beta\beta^\prime}\,
 \pa_\mu g_{\alpha^\prime\beta^\prime}
+\frac{1}{2} \pa_{\alpha} g_{\beta\mu}\,
\,g^{\alpha\alpha^\prime}g^{\beta\beta^\prime}
 \pa_{\alpha^\prime} g_{\beta^\prime\nu}
 -\frac{1}{2} \pa_{\alpha} g_{\beta\mu}\,
\,g^{\alpha\alpha^\prime}g^{\beta\beta^\prime} \pa_{\beta^\prime}
g_{\alpha^\prime\nu}\\
=g^{\alpha\alpha^\prime}g^{\beta\beta^\prime}\,\Big( \frac{1}{4}
\pa_\nu g_{\alpha\beta}\,\,
 \pa_\mu g_{\alpha^\prime\beta^\prime}
+\frac{1}{2} \pa_{\alpha} g_{\beta\mu}\, \,
 \pa_{\alpha^\prime} g_{\beta^\prime\nu}
 -\frac{1}{2} \pa_{\beta^\prime} g_{\beta\mu}\,
\, \pa_{\alpha} g_{\alpha^\prime\nu}\Big)\\
 -\frac{1}{2}g^{\alpha\alpha^\prime}g^{\beta\beta^\prime}\,
\Big(\pa_{\alpha} g_{\beta\mu}\, \, \pa_{\beta^\prime}
g_{\alpha^\prime\nu}-\pa_{\beta^\prime} g_{\beta\mu}\, \,
\pa_{\alpha} g_{\alpha^\prime\nu}\Big)\\
=g^{\alpha\alpha^\prime}g^{\beta\beta^\prime}\,\Big( \frac{1}{4}
\pa_\nu g_{\alpha\beta}\,\,
 \pa_\mu g_{\alpha^\prime\beta^\prime}
 -\frac{1}{8} \pa_{\mu} g_{\beta\beta^\prime}\,
\, \pa_{\nu} g_{\alpha\alpha^\prime} +\frac{1}{2} \pa_{\alpha}
g_{\beta\mu}\, \,
 \pa_{\alpha^\prime} g_{\beta^\prime\nu}
 \Big)\\
 -\frac{1}{2}g^{\alpha\alpha^\prime}g^{\beta\beta^\prime}\,
\Big(\pa_{\alpha} g_{\beta\mu}\, \, \pa_{\beta^\prime}
g_{\alpha^\prime\nu}-\pa_{\beta^\prime} g_{\beta\mu}\, \,
\pa_{\alpha} g_{\alpha^\prime\nu}\Big)
\end{multline}
where the last inequality follows from \eqref{eq:wavec2}.

Taking the trace of \eqref{eq:Rexpr} and using \eqref{eq:calc},
\eqref{eq:wavecord} we obtain \beq \label{eq:Rexpr1}
R_{\mu\nu}=-\frac{1}{2} g^{\alpha\beta}\pa_\alpha\pa_\beta
g_{\mu\nu} +\Gamma_{\nu\alpha\beta} \Gamma_{\mu}^{\,\,\,
\alpha\beta} +\frac
{1}{2}g^{\alpha\alpha^\prime}g^{\beta\beta^\prime}\, \Big (\pa_\mu
g_{\alpha^\prime\beta^\prime}\, \pa_\a g_{\b\nu} +\pa_\nu
g_{\alpha^\prime\beta^\prime} \,\pa_\a g_{\b\mu} -\pa_\nu
g_{\alpha^\prime\beta^\prime}\, \pa_\mu g_{\a\b}\Big), \eq
 Using \eqref{eq:gammaproduct} and \eqref{eq:metricproducttwo}
 we get
 \begin{multline}
R_{\mu\nu}=-\frac{1}{2} g^{\alpha\beta}\pa_\alpha\pa_\beta
g_{\mu\nu}+g^{\alpha\alpha^\prime}g^{\beta\beta^\prime}\,\Big(
-\frac{1}{4} \pa_\nu g_{\alpha\beta}\,\,
 \pa_\mu g_{\alpha^\prime\beta^\prime}
 +\frac{1}{8} \pa_{\mu} g_{\beta\beta^\prime}\,
\, \pa_{\nu} g_{\alpha\alpha^\prime}
 \Big)\\
 +\frac{1}{2}g^{\alpha\alpha^\prime}g^{\beta\beta^\prime}
 \pa_{\alpha} g_{\beta\mu}\, \,
 \pa_{\alpha^\prime} g_{\beta^\prime\nu}
-\frac{1}{2}g^{\alpha\alpha^\prime}g^{\beta\beta^\prime}\,
\Big(\pa_{\alpha} g_{\beta\mu}\, \, \pa_{\beta^\prime}
g_{\alpha^\prime\nu}-\pa_{\beta^\prime} g_{\beta\mu}\, \,
\pa_{\alpha} g_{\alpha^\prime\nu}\Big)\\
+\frac{1}{2}g^{\alpha\alpha^\prime}g^{\beta\beta^\prime}\,\Big(
\big(\pa_\mu g_{\a^\prime\b^\prime}\, \pa_\a g_{\b\nu}- \pa_\a
g_{\a^\prime\b^\prime}\, \pa_\mu g_{\b\nu}\big) +\big(\pa_\nu
g_{\a^\prime\b^\prime}\, \pa_\a g_{\b\mu}- \pa_\a
g_{\a^\prime\b^\prime}\, \pa_\nu g_{\b\mu}\big) \Big)\\
\frac{1}{4}g^{\alpha\alpha^\prime}g^{\beta\beta^\prime} \Big(\big(
\pa_{\beta^\prime} g_{\a^\prime\alpha} \,\pa_\mu g_{\b\nu}
-\pa_{\mu} g_{\a^\prime\alpha} \,\pa_{\beta^\prime}g_{\b\nu}\big)
+\big( \pa_{\beta^\prime} g_{\a^\prime\alpha} \,\pa_\nu g_{\b\mu}
-\pa_{\nu} g_{\a^\prime\alpha}
\,\pa_{\beta^\prime}g_{\b\mu}\big)\Big)
 \end{multline}
The result now follows.
\end{proof}
Let $m$ denote the standard Minkowski metric
$$
m_{00}=-1,\qquad  m_{ii}=1,\quad \text{if}\quad i=1,...,3,
\quad\text{and}\quad m_{\mu\nu}=0,\quad \text{if}\quad \mu\neq \nu,
$$
Define a 2-tensor $h$ from the decomposition
$$
g_{\mu\nu}=m_{\mu\nu}+h_{\mu\nu}.
$$
Let $m^{\mu\nu}$ be the inverse of $m_{\mu\nu}$.
Then for small $h$
$$
H^{\mu\nu}=g^{\mu\nu}-m^{\mu\nu}=-h^{\mu\nu}+O^{\mu\nu}(h^2),
\qquad\text{where}\qquad
h^{\mu\nu}=m^{\mu\mu^\prime}m^{\nu\nu^\prime}h_{\mu^\prime\nu^\prime}
$$
and $O^{\mu\nu}(h^2)$ vanishes to second order at $h=0$.

As a consequence of Lemma \ref{Einstwavec} we get:
\begin{lemma} \label{Einstwavecquad}
If Einstein's equation's \eqref{eq:Einst} and the wave coordinate
condition \eqref{eq:wavecord} hold then \beq \label{Eired}
{\Boxr}_g h_{\mu\nu}=F_{\mu\nu}(h)(\pa h,\pa h) \eq where
$F_{\mu\nu}(h)(\pa h,\pa h)$ is a quadratic form in $\pa h$ with
coefficients that are smooth functions of $h$. More precisely,
\beq\label{eq:defF} F_{\mu\nu}(h)(\pa h,\pa h)=P(\pa_\mu h,\pa_\nu
h)+Q_{\mu\nu}(\pa h,\pa h) +G_{\mu\nu}(h)(\pa h,\pa h) \eq where
\beq\label{eq:defP} P(\pa_\mu h,\pa_\nu h)=\frac{1}{4}
m^{\alpha\alpha^\prime}\pa_\mu h_{\alpha\alpha^\prime} \,
m^{\beta\beta^\prime}\pa_\nu h_{\beta\beta^\prime}  -\frac{1}{2}
m^{\alpha\alpha^\prime}m^{\beta\beta^\prime} \pa_\mu
h_{\alpha\beta}\, \pa_\nu h_{\alpha^\prime\beta^\prime} \eq and
\begin{multline} \label{eq:nullform} Q_{\mu\nu}(\pa h,\pa h)= \pa_{\alpha}
h_{\beta\mu}\, \, m^{\alpha\alpha^\prime}m^{\beta\beta^\prime}
 \pa_{\alpha^\prime} h_{\beta^\prime\nu}
-m^{\alpha\alpha^\prime}m^{\beta\beta^\prime} \big(\pa_{\alpha}
h_{\beta\mu}\,\,\pa_{\beta^\prime} h_{\alpha^\prime \nu}
-\pa_{\beta^\prime} h_{\beta\mu}\,\,\pa_{\alpha}
h_{\alpha^\prime\nu}\big)\\
+m^{\a\a'}m^{\b\b'}\big (\pa_\mu h_{\a'\b'}\, \pa_\a h_{\b\nu}-
\pa_\a h_{\a'\b'} \,\pa_\mu h_{\b\nu}\big )  +
m^{\a\a'}m^{\b\b'}\big (\pa_\nu h_{\a'\b'} \,\pa_\a gh_{\b\mu} -
\pa_\a
h_{\a'\b'}\, \pa_\nu h_{\b\mu}\big )\\
+\frac 12 m^{\a\a'}m^{\b\b'}\big (\pa_{\b'} h_{\a\a'}\, \pa_\mu
h_{\b\nu} - \pa_{\mu} h_{\a\a'}\, \pa_{\b'} h_{\b\nu} \big )
+\frac 12 m^{\a\a'}m^{\b\b'} \big (\pa_{\b'} h_{\a\a'}\, \pa_\nu
h_{\b\mu} - \pa_{\nu} h_{\a\a'} \,\pa_{\b'} h_{\b\mu} \big )\nn
 \end{multline}
 is a null form and
$G_{\mu\nu}(h)(\pa h,\pa h)$ is a quadratic form in $\pa h$ with
coefficients smoothly dependent on $h$ and vanishing when $h$
vanishes: $G_{\mu\nu}(0)(\pa h,\pa h)=0$.

Furthermore
\beq
\label{eq:wavec}
m^{\alpha\beta}\pa_\alpha h_{\beta\mu}
=\frac{1}{2}m^{\alpha\beta}\pa_\mu h_{\alpha\beta}+G_\mu(h)(\pa h),\qquad
\text{or}\qquad
\pa_\alpha H^{\alpha\nu}=\frac{1}{2}g_{\alpha\beta}
 \big(m^{\nu\mu}+H^{\nu\mu}\big) \pa_\mu H^{\alpha\beta}
\eq
where $G_\mu(h)(\pa h)$ is a linear function of $\pa h$ with coefficients
that are smooth functions of $h$ and that vanishes when $h$ vanishes:
$G_\mu(0)(\pa h)=0$.
\end{lemma}
Observe that the terms in \eqref{eq:defP} do not satisfy the
classical null condition. However  the trace $m^{\mu\nu}
h_{\mu\nu}$ satisfies a nonlinear wave equation with semilinear
terms obeying the  the null condition:
$$
g^{\alpha\beta}\pa_\alpha\pa_\beta m^{\mu\nu} h_{\mu\nu}
=Q(\pa h,\pa h)+G(h)(\pa h,\pa h).
$$

\section{The initial data }
In this section we discuss
the initial data for which the results of our paper apply.
We shall consider the asymptotically flat data, satisfying  a
global smallness condition, with the property that it coincides with
the Schwarzschild data  outside
the ball of radius one.

We start by showing the existence of such data. Let
$(g_{0},k_{0})$ be asymptotically flat initial data for the
Einstein equations consisting of the Riemannian metric $g_{0}$ and
a second fundamental form $k_{0}$. The initial data for the vacuum
Einstein satisfy the constraint equations
\begin{align}
&R_0-(\tr k_0)^{2} + |k_0|^{2}=0,\label{eq:enconstr}\\
& \nabla^{j} {k_0}_{ij}-\nabla_{i} \tr k_0 =0\label{eq:momconstr}
\end{align}
We restrict our attention to the time-symmetric case $R_0=k_0=0$.
Then, if $(g_0,k_0)$ is sufficiently close to the Minkowski data and $g_0$ satisfies
the parity condition $g_0(x)=g_0(-x)$,
by the results of Corvino \cite{Co} and Chrusciel-Delay \cite{C-D} one can
construct a new set of initial data $(g,k)$ with the properties that
\begin{itemize}
\item The initial data $(g,k)$ coincides with $(g_{0}, k_{0})$ on
the ball of radius 1/2. \item $(g,k)$ is exactly the Schwarzschild
data $(g_{S}^x,0)$ of mass $M$ outside $B_1$, the ball of radius
one.
\end{itemize}
At this point we specify the smallness conditions:
\begin{equation}\label{eq:smallness}
M\le \epsilon, \qquad
\sum_{0\le |I|\le N} \Big (\|\pa_x^I (g-\de)\|_{L^2(B_1)} +
\sum_{0\le |J|\le N-1}\|\pa^{J}_x k\|_{L^2(B_1)}\Big )\le \epsilon
\end{equation}
for some sufficiently large integer $N$. Here $\pa_x^I$  denotes the derivative
$\pa_{x_1}^{I_1}...\pa_{x_n}^{I_n}$, where
 $(I_1,..,I_n)$ is an arbitrary multi-index with the
 property that  $I_1+..+I_n= |I|$.

We have two expressions for the Schwarzschild metric
in isotropic and wave coordinates:
\begin{align}
&g_{S} = -\frac{(1- M/{r})^{2}}{(1+M/{r})^{2}}dt^{2} +
(1+\frac {M}{r})^{4} dx^{2}, \label{eq:Siso}\\
&g_{s} = - \frac {r-2M}{r+2M} dt^{2} + \frac {r+2M}{r-2M} dr^{2}
+ (r+2M)^{2} (d\theta^{2} + \sin^{2}\theta d\phi^{2})
\label{eq:Swave}
\end{align}
The expressions $g_S^x$ and $g_s^x$ will denote the spatial parts
the Schwarzschild metric in respective coordinates.
Observe that
\beq
\label{eq:ASch}
g_{s} = m + \frac {4M}{r} (dt^{2} + dx^{2}) + O (r^{-2}),
\eq
We now find the coordinate change transforming the metric $g_{S}$
into $g_{s}$. Set
\beq
\label{eq:change}
t= \tau,\qquad
r=\rho  + \frac{M^{2}}\rho
\ee
In the coordinates $\tau, \rho$ the metric $g_{s}$ takes the
form $g_{S}$. This change of coordinates is one-to-one for the
values $\rho > M$. Since the mass $M<< 1$ we can  define
the change of coordinates $r=\Phi (\rho)$, where $\Phi$ coincides
with the map \eqref{eq:change} for $\rho >1$ and the identity
transformation for $\rho\le 1/2$.
Thus we have constructed the initial data $(g,k)$ such that
\begin{itemize}
\item The initial data $(g,k)$ coincides (in new coordinates) with
$(g_{0}, k_{0})$ on the ball of radius 1/2. \item $(g,k)$ is
exactly the Schwarzschild data $(g^x_{s},0)$ outside the ball of
radius one. \item Moreover, the new data still obeys the smallness
condition \eqref{eq:smallness}.
\end{itemize}
The constructed metric is already in wave coordinates on its
Schwarzschild part. We now describe the procedure which produces
the initial data $(g, \pa_{t}g)$ associated with $(g,k)$ and
satisfying the wave coordinate condition.

Recall that {\it {a priori}} we are only given the spatial part of
the metric $g_{ij}$ together with a second fundamental form
$k_{ij}$. We now define the full space-time metric $g_{\a\b}$ on the
Cauchy hypersurface $\Si_{0}$ as follows: \beq \label{eq:sptimet}
g_{0i}=0, \qquad g_{00} = -a(r), \ee where the function
\begin{align*}
&a(r)= \frac {r-2M}{r+2M},\qquad {\text{for}}\,\,\,r>1\\
&a(r) = 1,\qquad {\text{for}}\,\,\,r\le \frac 12
\end{align*}
Thus defined metric coincides with the full Schwarzschild
metric $g_s$ for $r>1$.
We further define
\beq
\label{eq:timek}
\pa_{t} g_{ij} = -2 a k_{ij}
\ee
It remains to determine $\pa_{t} g_{0\a}$. We find it by satisfying
the wave coordinate condition
$$
g^{\b\mu} \pa_{\mu} g_{\a\b} = \frac 12 g^{\mu\nu}\pa_{\a} g_{\mu\nu}
$$
Setting $\a=0$ we obtain
$$
\frac 12 g^{00} \pa_{t} g_{00} = - g^{\b i}\pa_{i} g_{0\b} + \frac 12
g^{{ij}}\pa_{t} g_{ij}
$$
This defines $\pa_{t} g_{00}$. On the other hand setting $\a=i$ we
obtain
$$
g^{00}\pa_{t} g_{0i} = - g^{\b j}\pa_{j} g_{i\b} +
\frac 12 g^{\mu\nu}\pa_{i} g_{\mu\nu}
$$
This determines $\pa_{t} g_{0i}$.
Observe that since the metric $g$ coincides with the Schwarzschild
metric $g_{s}$, already satisfying the wave coordinate condition,
outside the ball of radius one, we have that on that set the initial data takes
the form $(g_{s},0)$.
Hence we constructed the initial data $(g,\pa_{t}g)$ with the
properties that
\begin{itemize}
\item The initial data $(g,\pa_{t}g)$ corresponds to the initial
data $(g,k)$ prescribed originally. \item $(g,\pa_{t}g)$ is
exactly the Schwarzschild data $(g_{s},0)$ outside the ball of
radius one. \item The initial data verifies the wave coordinate
condition. \item The initial data satisfies the smallness
condition
\begin{equation}\label{eq:small-data}
\sum_{0\le |I|\le N} \Big (\|\pa_x^I (g-m)\|_{L^2(B_1)} +
\sum_{0\le |J|\le N-1}\|\pa^{J}_x \pa_t g\|_{L^2(B_1)}\Big )\le \epsilon
\end{equation}
\end{itemize}
Now with the initial data $(g,\pa_{t}g)$ we solve the reduced Einstein
equations \eqref{eq:Einstred}. It follows from the proof of Lemma
\ref{Einstwavec} that, in the notation $\Gamma^{\lambda} = g^{\a\b}
\Gamma^{\lambda}_{\a\b}$, the reduced Einstein equations can be written in
the form:
\beq
\label{eq:reduE}
{\bf {R}}_{\a\b} - \frac 12 (D_{\a} \Gamma_{\b} + D_{\b}\Gamma_{\a})-
\Gamma_\si N^\si_{\a\b}(g,\pa g)=0.
\ee
Here $D$ denotes a covariant derivative with respect to the
space-time metric $g$ and $ N^\si_{\a\b}$ are some given functions
depending on $g$ and $\pa g$.
Observe that the initial data $(g,\pa_{t}g)$ were chosen in such
a way that the wave coordinate condition $\Gamma^{\lambda}=0$ is
satisfied on the initial hypersurface $\Si_{0}$. We now argue that
this condition is propagated, i.e, the solution of the reduced Einstein
equations \eqref{eq:reduE}  obeys $\Gamma^{\lambda}=0$ on {\it{any}}
hypersuface $\Si_{t}$. We would have thus
shown that a solution of the reduced Einstein equations is, in fact,
a solution of the vacuum Einstein equations.

To prove that $\Gamma^{\lambda} =0$ we differentiate \eqref{eq:reduE}
and use the contracted Bianchi identity
$D^{\b} {\bf {R}}_{\a\b} = \frac 12 D_{\a} {\bf{R}}$
\begin{align*}
0&= 2(D^{\b}{\bf{R}}_{\a\b} - \frac 12 D_{\a} {\bf{R}}) =
D^{\b} D_{\a} \Gamma_{\b} + D^{\b}D_{\b} \Gamma_{\a} -
D_{\a} D^{\b} \Gamma_{\b}- 2D^\b (\Gamma_\si N^\si_{\a\b})-
D_\a(\Gamma_\si N^{\si\b}_{\b})\\ &= D^{\b}D_{\b} \Gamma_{\a} +
{\bf{R}}_{\a\ga}\Gamma^{\ga}-2(D_\b \Ga_\si ) N^{\si\b}_{\b}
-(D_\a\Ga_\si)  N^{\si\b}_{\b} - 2 \Ga_\si (D_\b N^{\si\b}_{\a})-
\Ga_\si  (D_\a N^{\si\b}_{\b})
\end{align*}
Therefore, $\Gamma^{\lambda}$ satisfies a covariant wave equation, on the
background determined by the constructed metric $g$, with
the initial condition $\Gamma^{\lambda}=0$. It remains to show that
$D_{t}\Gamma^{\lambda}=0$ on $\Si_{0}$ and the conclusion that
$\Gamma^{\lambda}\equiv 0$ will follows by the uniqueness result for
wave equation.

We recall that the initial data $(g,k)$ verifies the constraint
equations \eqref{eq:enconstr}, \eqref{eq:momconstr}, which imply
 that
on $\Si_{0}$
$$
{\bf{R}}_{TT} + \frac 12 {\bf{R}}=0,\qquad
{\bf{R}}_{Ti}=0,
$$
where $T= - (g_{00})^{-1} \pa_t$ is the unit future oriented normal to
$\Si_0$.
Therefore returning to \eqref{eq:reduE} we obtain that
\begin{align*}
&0={\bf{R}}_{00} + \frac 12 {\bf{R}}=
-(g_{00})^{-1} D_{t} \Gamma_{0} +  D^{i}\Gamma_{i},\\
&0={\bf{R}}_{0i}= \frac 12 D_{t}\Gamma_{i} + \frac 12 D_{i}\Gamma_{0}
\end{align*}
This finishes the proof that $\Gamma^{\lambda}\equiv 0$.

We also know that the time-independent Schwarzschild metric
$g_{s}$ is a solution of the Einstein vacuum equation
${\bf{R}}_{\a\b}=0$. Moreover, since $g_{s}$ satisfies the wave
coordinate condition it also verifies the reduced Einstein
equations \eqref{eq:reduE}. Since the initial data $(g,\pa_{t}g)=
(g_{s},0)$ outside the ball of radius two, constructed solution
will coincide with the Schwarzschild solution in the exterior of
the null cone developed from the sphere of radius one in
$\Si_{0}$.

We end the discussion of the initial data by comparing the light
cones of Minkowski and Schwarzschild spaces in the wave coordinates
of the   Schwarzschild space.
\begin{lemma}
For an arbitrary $R>2M$ the forward null cone of the metric $g_s$, intersecting the time slice $t=0$
along the sphere of radius $R$, is contained in the interior of the
Minkowski cone $t-r = R$.
\end{lemma}
\begin{proof}
The null cone intersecting the time slice $t=0$ along the
sphere of radius $R$ can be realized as the level hypersurface
$u=0$ of the optical function $u$ solving the eikonal equation
$$
g_{s}^{\a\b} \pa_{\a} u \,\,\pa_{\b} u =0
$$
with the initial condition that $u=0$ on the sphere of radius $R$
at time $t=0$. Because of the spherical symmetry of the
Schwarzschild metric $g_{s}$ and the initial condition we look for
a spherically symmetric solution $u=u(t,r)$. The eikonal equation
then reads
$$
\frac {r+2M}{r-2M} (\pa_{t} u )^{2} = \frac {r-2M}{r+2M} (\pa_{r} u )^{2}
$$
Let $t=\ga(r)$ be a null geodesic, originating from some point on the sphere of radius
$R$ at $t=0$, such that  $u(\ga(r),r)=0$. Then
$$
\pa_{t} u \dot \ga(r) + \pa_{r} u  =0
$$
Substituting this into the eikonal equation we obtain that
$$
\Big (\frac {r+2M}{r-2M}\Big )^{2}= |\dot\ga (r)|^{2}
$$
Taking the square root and integrating we obtain that
$$
\ga (r) = \ga (R) \pm \Big (r-R + 4M \ln \frac {r-2M}{R-2M}\Big )
$$
Thus the null geodesics are described the curves
$$
t=\pm \Big (r-R + 4M \ln \frac {r-2M}{R-2M}\Big )
$$
In particular, the forward null cone is contained in the
interior of the set $t\ge r-R$.
\end{proof}

\section {The null-frame and null-forms}\label{section:null}
Below we introduce a standard Minkowski null-frame used throughout the paper.
At each point $(t,x)$ we introduce a pair of null vectors
$(L,\Lb)$
$$
L^0=1, \quad L^i=x^i/|x|, \,\,\, \, i=1,2,3,\quad\text{and}\quad
\underline {L}^0=1,\quad \underline {L}^i=-x^i/|x|, \,\,\, \,
i=1,2,3,
$$
Adding two orthonormal tangent to the sphere $S^{2}$ vectors
$S_1,S_2$ which are orthogonal to $\omega$ defines a null
frame $(L,\Lb, S_1,S_2)$.
\begin{remark} Since $S^2$ does not admit a global
orthonormal frame $S_1,S_2$ we could alternatively introduce
a global frame induced by the projections of
the coordinate vector fields $e_i$.

Let $P$ be the orthogonal projection of a vector field in $\bold{R}^3$
along $\omega=x/|x|$ onto the tangent space of the sphere;
$
PV=V-\langle V,\omega\rangle \omega.
$
For $i\!=\! 1,2,3$ denote the projection of $\pa_i$ by
\begin{equation}\label{eq:deftan}
\pas_i\!=A_i^j\pa_j=\pa_i-\omega_i\omega^j\pa_j,\qquad\text{where}\qquad
A_i^j=(P e_i)^j=\delta_i^j-\omega_i\omega^j,\qquad i=1,2,3,
\end{equation}
where $e_i$ is the usual orthonormal basis in $\bold{R}^3$, and
the sums are over $j=1,2,3$ only. Let $\pab_0\!=\!
L^\alpha\pa_\alpha$ and $\pab_i=\pas_i$, for $i=1,2,3$. Then
a linear combination of the derivatives 
$\{\pab_0,...,\pab_3\}$ spans the tangent space of the forward
light cone.
\end{remark}

In what follows $A, B$ will denote any of the vectorfields $S_1, S_2$. 
We will use the summation conventions
$$
X^A A^\alpha = X^\beta  {S_1}_\beta S_1^\alpha + 
X^\beta {S_2}_\beta S_2^\alpha,
\qquad X_A Y_A =  X^\alpha Y^\beta {S_1}_\alpha {S_1}_\beta  +
  X^\alpha Y^\beta {S_2}_\alpha {S_2}_\beta. 
$$
Obvious generalizations of the above 
conventions will be used for higher order tensors.
 
We record the following null frame decomposition of a vector field $X=X^\a\pa_\a$:
$X^\alpha=X^L L^\alpha+
X^{\underline{L}}\underline{L}^\alpha+X^A A^\alpha$.
Relative to a
null frame the Minkowski metric $m$ has the following form
$$
m_{LL}=m_{\Lb\Lb}=m_{LA}= m_{\Lb A}=0, \qquad m_{L\Lb}=m_{\Lb L}=-2,
 \qquad m_{AB}=\delta_{AB},
$$
i.e.
$m_{\alpha\beta} X^\alpha Y^\beta=-2(X^L Y^{\underline{L}}+X^{\underline{L}}Y^L)
+X^A Y^A$.
Recall that we raise and lower indices of any tensor relative
to the Minkowski metric $m$, i.e.,  $X_\alpha=m_{\alpha\beta}X^\beta$.
We define $X_Y=m_{\alpha\beta} X^\alpha Y^\beta=X_\alpha Y^\alpha$.
Then $X_Y=X^L Y_L+X^{\underline{L}} Y_{\underline{L}}+X^A Y_A$.
It is useful to remember the following rule:
$$
X^{L}= -\frac 12 X_{\Lb},\quad X^{\Lb}=-\frac 12 X_{L},\quad X^A=X_A.
$$
Then
$$
m^{LL}=m^{\Lb\Lb}=m^{LA}= m^{\Lb A}=0, \qquad m^{L\Lb}=m^{\Lb L}=-1/2,
 \qquad m^{AB}=\delta^{AB}
$$
i.e. $m^{\alpha\beta} X_\alpha Y_\beta=-\frac{1}{2}\big(X_L Y_{\underline{L}}
+X_{\underline{L}} Y_L\big)+X_A Y_A$.

\begin{defi}\label{def:nullcoord}
Denote $q=r-t$ and $s=t+r$ the null coordinates of the Minkowski metric $m$
and $\pa_{q} =\frac 12 (\pa_{r}-\pa_{t})$
and $\pa_{s}=\frac 12 (\pa_{t}+\pa_{r})$, the corresponding null vector fields
\end{defi}
Let
$k_{XY}=k_{\alpha\beta} X^\alpha Y^\beta$.
Then
\beq\label{eq:trace}
\operatorname{tr}\, k=
m^{\alpha\beta} k_{\alpha\beta}
=-\frac{1}{2}\big(k_{L\underline{L}}+k_{\underline{L}L}\big)+\overline{\operatorname{tr}}\, k
\eq
where
\beq\label{eq:tantrace}
\overline{\operatorname{tr}}\,k=\delta^{AB}k_{AB}
=\overline{\delta}^{ij}k_{ij},\qquad\text{and}\qquad
\overline{\delta}^{ij}=\delta^{ij}-\omega^i\omega^j
\eq
where the sum is over $i,j=1,2,3$ only.

If $k$ and $p$ are symmetric it follows that
\begin{multline}\label{eq:traceofproduct}
p_{\alpha\beta} k^{\alpha\beta}=
m^{\alpha\alpha^\prime} m^{\beta\beta^\prime}p_{\alpha\beta}
k_{\alpha^\prime\beta^\prime}\\
=\frac{1}{4}\big(p _{LL} k_{\underline{L}\underline{L}}
+p_{\underline{L}\underline{L}} k_{{L}{L}}
+ 2 p_{L\underline{L}} k_{\underline{L}L}\big)
-\delta^{AB}\big(p_{A L}k_{B\underline{L}}
+p_{A \underline{L}}k_{B{L}}\big)
+\delta^{AB}\delta^{A^\prime B^\prime}p_{AA^\prime}k_{BB^\prime}\\
=\frac{1}{4}\big(p _{LL} k_{\underline{L}\underline{L}}
+p_{\underline{L}\underline{L}} k_{{L}{L}}
+ 2 p_{L\underline{L}} k_{\underline{L}L}\big)
-\overline{\delta}^{ij}\big(p_{i L}k_{j\underline{L}}
+p_{i\underline{L}}k_{j{L}}\big)
+\overline{\delta}^{ij}\overline{\delta}^{i^\prime j^\prime}p_{ii^\prime}k_{jj^\prime}
\end{multline}
\begin{lemma}\label{nullframeP}
With $P(p,k)$ given by \eqref{eq:defP} we have for symmetric 2-tensors $p$ and $k$:
\begin{multline}\label{eq:nullframeP}
P(p,k)=\frac{1}{4} m^{\alpha\beta} p_{\alpha\beta}m^{\alpha\beta}
k_{\alpha\beta}-\frac{1}{2}m^{\alpha\alpha^\prime}
m^{\beta\beta^\prime}p_{\alpha\beta} k_{\alpha^\prime\beta^\prime}
\\
=-\frac{1}{8}\big(p_{LL} k_{\underline{L}\underline{L}}
+p_{\underline{L}\underline{L}} k_{{L}{L}}\big)
-\frac{1}{4}\delta^{AB}\delta^{A^\prime
B^\prime}\Big(2p_{AA^\prime}k_{BB^\prime}-
p_{AB} k_{A^\prime B^\prime}\Big)\\
+\frac{1}{4}\delta^{AB}\big(2p_{A L}k_{B\underline{L}} +2p_{A
\underline{L}}k_{B{L}}- p_{AB}
k_{L\underline{L}}-p_{L\underline{L}}k_{AB} \big)
\end{multline}
i.e. at least one of the factors contains only tangential components.
\end{lemma}

Furthermore
$$
p^{\alpha\beta}\pa_\alpha =p^{L\beta}\pa_L+p^{\underline{L}\beta}
\pa_{\underline{L}}
+p^{A\beta}\pa_A
=-\frac{1}{2}p_{\underline{L}}^{\,\,\,\beta}\pa_L
-\frac{1}{2} p_L^{\,\,\,\beta}\pa_{\underline{L}}
+p^{A\beta}\pa_A
$$

We introduce the following notation. Let ${\cal T}=\{L,S_1,S_2\}$,
${\cal U}=\{\Lb,L,S_1,S_2\}$, ${\cal L}=\{L\}$ and ${\cal S}=\{S_1,S_2\}$.
For any two of these families ${\cal V}$ and ${\cal W}$ and an
arbitrary two-tensor $p$ we denote
\begin{align}
|p|_{{\cal V W}}
&=\sum_{V\in{\cal V},W\in{\cal W} , }
|p_{\beta\gamma} V^\beta W^\gamma|,\label{eq:norm1}\\
|\pa p|_{{\cal V W}}
&=\sum_{U\in {\cal U},V\in{\cal V},W\in{\cal W} , }
|(\pa p)_{\alpha \beta\gamma}U^\alpha V^\beta W^\gamma|,
\label{eq:norm2}\\
|\overline{\pa} p|_{{\cal V W}}
&=\sum_{T\in {\cal T},V\in{\cal V},W\in{\cal W} , }
|(\pa p)_{\alpha \beta\gamma}T^\alpha V^\beta W^\gamma|
\label{eq:norm3}
\end{align}

Let $Q$ denote a null form, i.e. $Q_{\alpha\beta}(\pa \phi,\pa \psi)
= \pa_\alpha  \phi\, \pa_\beta \psi-\pa_\beta \phi\, \pa_\alpha \psi$
if $\alpha\neq \beta$ and $Q_{0}(\pa \phi,\pa \psi)=m^{\alpha\beta} \pa_\alpha  \phi\, \pa_\beta \psi$.

\begin{lemma}\label{tander} If $P$ is as in Lemma
\ref{nullframeP} then \beq\label{eq:tanP} |P(p,k)|\les |p\,|_{\cal
TU}|k|_{\cal TU} +|p\,|_{LL}|k|+|p\,|\,|k|_{LL} \eq If $Q(\pa
\phi,\pa \psi)$ is a null form then \beq\label{eq:nullformtan}
|Q(\pa \phi,\pa \psi)|\les |\overline{\pa} \phi| |\pa \phi| + |\pa
\phi||\overline{\pa} \psi| \eq

Furthermore
\beq \label{eq:dernullframeZ0}
|k^{\alpha\beta} \pa_\alpha\phi \, \pa_\beta \phi|\les
\big(|k|_{LL}|\pa\phi|^2+|k|\,|\pab \phi||\pa\phi|\big)
\eq
\beq \label{eq:dernullframeZ1}
|L_\alpha k^{\alpha\beta} \pa_\beta\phi |\les
\big(|k|_{LL}|\pa\phi|+|k|\,|\pab \phi|\big)
\eq
\beq \label{eq:dernullframeZ2}
|(\pa_\alpha k^{\alpha\beta}) \pa_\beta\phi |\les
\big(|\pa k|_{LL}+|\overline{\pa} k|\big)\,|\pa\phi|
+|\pa k|\, |\overline{\pa} \phi|
\eq
\end{lemma}
\begin{proof}
The proof of \eqref{eq:nullformtan} for the null form $Q_{0}$
follows directly from \eqref{eq:trace}. To prove the claim for the
null forms $Q_{\a\b}$ use that \beq\label{eq:decomalpha} \pa_{i} =
L_{i}(\pa_s+\pa_{q})+\bar\pa_i,\quad i=1,2,3,\qquad
\pa_0=L_0(\pa_s-\pa_q)
\end{equation}
Therefore,
$$
|Q_{\a\b}(\pa\phi,\pa\psi)|= |\pa_{\a}\phi\pa_{\b}\psi-\pa_{\b}\phi
\pa_{\a}\psi|\leq C|\pab\phi|\,| \pa\psi|
+ C|\pa\phi|\,|\pab\psi|
$$

The estimates
\eqref{eq:dernullframeZ0}-\eqref{eq:dernullframeZ2} follow from
\eqref{eq:traceofproduct}.
\end{proof}

\begin{lemma}\label{tander2}
If $k^{\alpha\beta}$ is a symmetric tensor
and $\phi$ a function then
\beq \label{eq:dernullframeZ}
|k^{\alpha\beta} \pa_\alpha \pa_\beta \phi|\les
\big(|k|_{LL}|\pa^2\phi|+|k|\,|\pab\pa \phi|\big)
\eq
Also,
with $\overline{\operatorname{tr} } \,k=\delta^{AB} k_{AB}=(\delta^{ij}-\omega^i\omega^j)k_{ij}$ we have
\beq \label{eq:nuZ}
\big|k^{\alpha\beta} \pa_\alpha\pa_\beta \phi-
k_{LL}\pa_{q}^2\phi - 2 k_{L\Lb}
\pa_s \pa_{q} \phi -r^{-1}\,\overline{\operatorname{tr} } \,k\,  \pa_q \phi\big|\les
|k|_{L{\cal T}}|\pab\pa\phi|+|k|\,\big(|\pab^2 \phi| + r^{-1}
|\bar\pa\phi|\big).
\eq
\end{lemma}
\begin{proof}
The estimate \eqref{eq:dernullframeZ} follow from
\eqref{eq:traceofproduct}.
We have
\beq
\pa_i\omega_j=r^{-1}(\delta_{ij}-\omega_i\omega_j)=
r^{-1}\overline{\delta}_{ij}
\eq
Furthermore $\pa_i=\bar\pa_i+\omega_i\pa_r$,
where $\pa_r=\omega^j\pa_j$ so $[\bar\pa_i,\pa_r]=
(\bar\pa_i\omega^k)\pa_k$ and
\begin{multline}
\pa_i\pa_j=(\bar\pa_i+\omega_i\pa_r)(\bar\pa_j+\omega_j\pa_r)\\=\bar\pa_i\bar\pa_j
+\omega_i\,\omega^k\,\bar\pa_j \pa_k+
\omega_j\omega^k\,\bar\pa_i\pa_k
+\omega_i\omega_j\pa_r^2
+(\bar\pa_i\omega_j)\pa_r+\omega_j(\bar\pa_i\omega^k)\pa_k \\
=\bar\pa_i\bar\pa_j
+\omega_i\,\bar\pa_j \pa_r+
\omega_j\,\bar\pa_i\pa_r
+\omega_i\omega_j\pa_r^2
+r^{-1}\overline{\delta}_{ij}\pa_r
-r^{-1}\omega_i\bar\pa_j
\end{multline}
Furthermore
\beq
\pa_0\pa_i=\pa_t(\bar\pa_i+\omega_i\pa_r)=
\bar\pa_i\pa_t+\omega_i\pa_t\pa_r
\eq
Hence
\begin{multline}
k^{\alpha\beta}\pa_\alpha\pa_\beta
=k^{00}\pa_t^2+2k^{0i}\omega_i\pa_t\pa_r+k^{ij}\omega_i\omega_j\pa_r^2+r^{-1}\overline{\operatorname{tr}}\, k\,
\pa_r\\
+k^{ij}\bar\pa_i\bar\pa_j-r^{-1} k^{ij}\omega_i
\bar\pa_j+2k^{0 j}\bar\pa_j\pa_t
+2k^{ij}\omega_i\bar\pa_j\pa_r
\end{multline}
If we substitute $\pa_t=\pa_s-\pa_q$, $\pa_r=\pa_s+\pa_q$
and identify
\beq
k_{LL}=k^{00}-2k^{0i}\omega_i+k^{ij}\omega_i\omega_j,
\qquad
k_{L\underline{L}}=-k^{00}+k^{ij}\omega_i\omega_j,
\qquad
k_{\underline{L}\underline{L}}=k^{00}+2k^{0i}\omega_i+k^{ij}\omega_i\omega_j
\eq
 and
\beq
-k^{0j}+k^{ij}\omega_j=k_{0}^{\,\,\, j}
+k_{i}^{\,\,\, j}\omega^i=k_{L}^{\,\,\, j},\qquad
k^{0j}+k^{ij}\omega_j=-k_{0}^{\,\,\, j}
+k_{i}^{\,\,\, j}\omega^i=k_{\underline{L}}^{\,\,\, j}
\eq
 we get
\begin{multline}
k^{\alpha\beta}\pa_\alpha\pa_\beta
=k_{LL}\pa_q^2+2k_{L\underline{L}}\pa_s\pa_q
+k_{\underline{L}\underline{L}}\pa_s^2
+r^{-1}\overline{\operatorname{tr}}\, k\,
\pa_q\\
+k^{ij}\bar\pa_i\bar\pa_j+r^{-1}\overline{\operatorname{tr}}\, k\,\pa_s-r^{-1} k^{ij}\omega_i
\bar\pa_j+2k_{L}^{\,\,\,j}\bar\pa_j\pa_q
+2k_{\underline{L}}^{\,\,\, j}\bar\pa_j\pa_s
\end{multline}
Finally, we can also write
\beq
2k_{L}^{\,\,\,j}\bar\pa_j\pa_q=
k_{L}^{\,\,\,j}\bar\pa_j(\omega^k\pa_k-\pa_t)
=k_{L}^{\,\,\,j}\omega^k\bar\pa_j\pa_k-
k_{L}^{\,\,\,j}\bar\pa_j\pa_t+r^{-1}k_{L}^{\,\,\,j}
\bar\pa_j,
\eq
 since $(\bar\pa_j\omega^k)\pa_k=r^{-1}\bar\pa_j$.
The inequality \eqref{eq:nuZ} now follows.
\end{proof}

\begin{cor}\label{cor:BoxH}
Let $\phi$ be a solution of the reduced wave equation $\Boxr_g\phi
=F$ with a metric $g$ such that $H^{\a\b}=g^{\a\b}-m^{\a\b}$
satisfies the condition that $|H^{L\Lb}|<\frac 14$. Then \beq
\label{eq:phiwave} \Big|\Big(4\pa_s -
\frac{H_{LL}}{2g^{L\Lb}}\pa_q-\frac{\overline{\operatorname{tr}}\, H
+H_{L\Lb}}{2g^{L\underline{L}}\,\, r} \Big) \pa_q
(r\phi)+\frac{rF}{2g^{L\underline{L}}}\Big|\les r|\triangle_\omega
\phi|+ |H|_{L{\cal T}}\,\, r\, |\pab\pa\phi|+ |H|\,\big(r\,|\pab^2
\phi| + |\bar\pa\phi|+r^{-1}|\phi|\big) \eq where
$\triangle_\omega=\bar\triangle=\delta^{ij}\bar\pa_i \bar\pa_j$.
\end{cor}
\begin{proof}
Define the new metric
$$
\tilde g^{\a\b} = \frac {g^{\a\b}}{-2 g^{L\Lb}}.
$$
The equation $g^{\a\b}\pa_{\a}\pa_{\b}\phi = F$ then takes the form
$$
\tilde g^{\a\b} \pa_{\a}\pa_{\b}\phi = \frac {F}{-2g^{L\Lb}},
$$
which can also be written as
$$
\Box \phi + (\tilde g^{\a\b}-m^{\a\b}) \pa_{\a}\pa_{\b}\phi = \frac {F}{-2g^{L\Lb}}
$$
Let $k^{\a\b}$ be the tensor $k^{\a\b} = (\tilde g^{\a\b}-m^{\a\b})$
Observe that
\begin{align*}
k^{\a\b} = (-2g^{L\Lb})^{-1} \big(g^{\a\b} + 2 m^{\a\b}
g^{L\Lb}\big ) & = (-2g^{L\Lb})^{-1} \big ( H^{\a\b} + m^{\a\b}
(2g^{L\Lb} + 1)\big )\\ &= (-2g^{L\Lb})^{-1} \big ( H^{\a\b} + 2 m^{\a\b}
H^{L\Lb} \big )
\end{align*}
Thus,
\begin{equation}
\label{eq:Lbdisappear}
k_{L\Lb}=0,\qquad \quad k_{L{\cal T}}= (-2g^{L\Lb})^{-1} H_{L{\cal T}},
\qquad \quad \overline{\operatorname{tr}}\, k \,
=(-2g^{L\Lb})^{-1} \big(\overline{\operatorname{tr}}\,H \,+H_{L\Lb}\big)
\end{equation}
Moreover, $|k|\les |H|$, since $g^{L\Lb}= H^{L\Lb} - \frac 12 $
and by the assumptions of the Corollary $|H^{L\Lb}|<\frac 14$.

Now using \eqref{eq:nuZ} of Lemma \ref{tander2}, with the condition that
$k_{L\Lb}=0$, together with the
decomposition
$$
\Box \phi = -\pa_t^2 \phi + \triangle \phi = \frac 1r
(\pa_t+\pa_r)(\pa_r-\pa_t) r \phi + \triangle_\omega\phi= \frac 4r
\pa_s \pa_q r \phi + \triangle_\omega\phi.
$$
we find that the identity
$\Box \phi + k^{\a\b} \pa_{\a}\pa_{\b}\phi = ({-2g^{L\Lb}})^{-1} F$
leads to the inequality
$$
\big|4 \pa_{s}\pa_{q} r \phi +r k_{LL}\pa_{q}^2\phi
+\overline{\operatorname{tr}}\, k
\,\pa_q\phi+(2g^{L\underline{L}})^{-1} r F\big|\les
r|\triangle_\omega \phi| + r |k|_{L{\cal T}}|\pab\pa\phi|+
|k|\,\big(r\, |\pab^2 \phi| +|\bar\pa\phi|\big)
$$
Finally, identity \eqref{eq:Lbdisappear} and a crude  estimate $|k|\les
|H|$ yield the desired result.
\end{proof}

\section{The weak null condition and
asymptotic expansion of Einstein's equation's in wave coordinates}
\label{weaknullcond}

Let us now first describe the weak null condition. The results of this section
appear in \cite{L-R}.
Consider the Cauchy problem for a system of nonlinear wave equations in three space dimensions:
\beq
\label{eq:quas}
-\sq u_i\,\,=F_i(u,u^\prime, u^{\prime\prime}), \qquad i=1,...,N,\qquad
u=(u_1,...,u_N),
\eq
where $ -\square=-\pa_t^2+\sum_{j=1}^3\pa_{x^j}^2$.
We assume that $F$ is a function of $u$ and its derivatives
of the form
\beq
\label{eq:nonform}
\,\,\,
F_i(u,u^\prime, u^{\prime\prime})={a_{i \alpha\beta}^{jk} \partial^\alpha u_j
\, \partial^\beta u_k}+G_i(u,u^\prime,u^{\prime\prime}),
\eq
where  $G_i(u,u^\prime,u^{\prime\prime})$ vanishes to third order as
$(u,u^\prime,u^{\prime\prime})\!\to\! 0$ and
$a_{i\, \alpha\beta}^{jk}\!=\!0$ unless $|\alpha|\!\leq \! |\beta|\!\leq 2\!$ and
$|\beta|\!\geq\! 1$.
Here we used the summation convention over repeated indices.
We assume that the initial data
\beq
\label{eq:indata}
u(0,x)=\varepsilon u_0(x)\in C^\infty
,\quad  u_t(0,x)=\varepsilon u_1(x)\in C^\infty
\eq
is small and decays fast as $|x|\to \infty$.
We are going to determine conditions on the nonlinearity
such that the equation \eqref{eq:quas} is
compatible with
the asymptotic expansion as $|x|\to \infty$ and $|x|\sim t$
\beq
\label{eq:asymexp}
u(t,x)\sim  \varepsilon
U(q,s,\omega)/|x|,\quad \text{ where }\quad
q=|x|-t,\,\,\,s=\varepsilon \ln{|x|},\,\,\, \omega=x/|x|,
\eq
for all sufficiently small $\epsilon >0$.
The linear and some nonlinear wave equations allow for such an expansion
with $U$ independent of $s$ and the next term
decaying like $\varepsilon/|x|^2$, see \cite{H1,H2}.
Substituting \eqref{eq:asymexp} into \eqref{eq:quas}
and equating powers of order $\varepsilon^2/|x|^2$
we see that
\beq
\label{eq:assystem}
2\partial_s\partial_q U_i =A_{i\, mn}^{jk}(\omega)
(\partial_q^m U_j)(\partial_q^n U_k),\qquad\quad U\big|_{s=0}=F_0,
\eq
 where
\beq \label{eq:nonstr} A_{i,mn}^{jk}(\omega)=\sum_{|\alpha|=m,
|\beta|=n}{a_{i,\alpha\beta}^{jk}
\hat{\omega}^\alpha\hat{\omega}^\beta},\qquad \text{where }
\hat{\omega}=(-1,\omega)\quad\text{and } \hat{\omega}^\alpha=
\hat{\omega}_{\alpha_1}...\hat{\omega}_{\alpha_k} \eq
 In fact, $ \Box
u=-\ve^{-1}\pa_s\pa_q ( r u)+\text{angular derivatives} $ and $
\pa_\mu=\hat{\omega}_\mu\pa_q+\text{tangential derivatives} $.

One can show that \eqref{eq:quas}-\eqref{eq:indata} has a solution
as long as $\varepsilon\log{t}$ is bounded, provided that
$\varepsilon>0$ is sufficiently small and the solution of
\eqref{eq:assystem} exists up to that time,
\cite{J-K,H1,H2,L1,L2}. The only exception is the case
$A^{jk}_{i00}\neq 0$, which has shorter life span. In cases where
the solution of \eqref{eq:assystem} blows up it has been shown
that  solutions of  \eqref{eq:quas}-\eqref{eq:indata} also break
down in some finite time $T_\varepsilon\leq e^{C/\varepsilon}$,
\cite{J1,H1,A1}. John's example was \beq \label{eq:blowupex1}
\square u= u_t \,\triangle u \eq for which \eqref{eq:assystem}  is
the Burger's equation
$(2\pa_s-U_{\!q}\pa_q) U_{\!q}=0$, which is known to blow up.
The equation
\beq
\label{eq:blowupex2}
\square u= u_t^2
\eq
is another example where solutions blow up,
for which  \eqref{eq:assystem} is $\pa_s U_{\!q}=U_{\!q}^2$,
that also  blows up.

The {\it null condition} of \cite{K2} is equivalent to
\beq
\label{eq:ncond}
A_{i\, mn}^{jk}(\omega)=0\quad\text{for all } (i,j,k,m,n),\,\,\omega\in
{\Bbb S}^2.
\eq
The results of \cite{C1}, \cite{K2} assert
that \eqref{eq:quas}-\eqref{eq:indata} has  global solutions
for all sufficiently small initial data, provided that the null
condition is satisfies.
In this case the asymptotic equation
\eqref{eq:assystem} trivially can be solved globally.
Moreover, similar to the linear case, its solutions approach a
limit as $s\!\to\! \infty$
and the solutions of  \eqref{eq:quas}-\eqref{eq:indata}
decay like solutions of linear equations. A typical example of an
equation satisfying the null condition is
\beq
\label{eq:example0}
\Box u=u_t^2-|\nabla_x u|^2
\eq
There is however a more general class of nonlinearities for which solutions of
\eqref{eq:assystem} do not
blow up:

We say that
a system \eqref{eq:quas} satisfies the
{\it weak null condition} if the solutions of the corresponding asymptotic
system
\eqref{eq:assystem} exist for all $s$ and if the solutions together
with its derivatives
grow at most exponentially in $s$ for all initial data decaying
sufficiently fast in $q$.

Under the weak null condition assumption solutions of
\eqref{eq:assystem}
satisfy the equation \eqref{eq:quas} up to terms of order
$\varepsilon^2\!/|x|^{3-C\varepsilon}\!\!$, but
need only decay like $\varepsilon/ |x|^{1-C\varepsilon}\!\!$.
An example of the equation satisfying the weak null condition is given by
\beq
\label{eq:example1}
\Box u=u \triangle u
\eq
In \cite{L2} it was proven that \eqref{eq:example1} have small
global solutions in the spherically symmetric case and recently \cite{A3}
established this result without the symmetry assumption.  The
equation \eqref{eq:example1}
appears to be similar to  \eqref{eq:blowupex1} but a closer look
shows that the corresponding asymptotic equation:
\beq
\label{eq:ex1asym}
 (2\pa_s -U\pa_q )U_q=0
\eq
 has global solutions growing exponentially in $s$, see \cite{L2}.
The system
\beq
\label{eq:example2}
\square u=v_t^2 ,\qquad \square v=0
\eq
is another example that satisfies the weak null condition.
The equation \eqref{eq:example2} appears to resemble \eqref{eq:blowupex2}.
The system however decouples: \,$v$ satisfies a
linear homogeneous equation
and given $v$ we have a linear inhomogeneous equation for $u$, and
global existence follows.
The corresponding asymptotic system is
\beq
\label{eq:ex2asym}
\partial_s \partial_q U=(\partial_q V)^2,\qquad \pa_s\pa_q V=0
\eq
The solution of the second equation in \eqref{eq:ex2asym} is
independent of $s$: $V_q=V_q(q,\omega)$ and substituting this into the first
equation we see that
$U_q(s,q,\omega)=s V_q(q,\omega)^2$ so $\pa u$
only decays like $|x|^{-1}\ln{|x|}$.

We show below that the Einstein vacuum equations in wave coordinates
satisfy the weak null condition, i.e.
that the corresponding
asymptotic system \eqref{eq:assystem}
admits global solutions. In fact, each of the quadratic nonlinear terms in the Einstein
equations is
either of the type appearing in
\eqref{eq:example0}, \eqref{eq:example1} or \eqref{eq:example2}.

\begin{theorem}
Let $h$ be a symmetric 2-tensor and let \beq \label{eq:asymh}
h_{\mu\nu}(t,x)\sim  \varepsilon U_{\mu\nu}(s,q, \omega)/|x|,\quad
\text{ where }\quad q=|x|-t,\,\,\,s=\varepsilon \ln{|x|},\,\,\,
\omega=x/|x|. \eq be an asymptotic ansatz. Then the asymptotic
system for the the Einstein equations in wave coordinates
\eqref{Eired}, obtained by formally equating the terms with the
coefficients $\epsilon^{2}{|x|^{-2}}$, takes the following form:
\beq\label{eq:Einstasympt} \big(2\pa_s-U_{LL}\pa_q\big)\pa_q
U_{\mu\nu} =L_\mu L_\nu P(\pa_q U,\pa_q U),\qquad \forall \mu, \nu
=0,...,3 \eq where $U_{LL}=m^{\alpha\alpha^\prime}
m^{\beta\beta^\prime}U_{\alpha^\prime\beta^\prime} L_\alpha
L_\beta$ and $P(\pa_q U,\pa_q U)= \frac{1}{4}\pa_q\tr U \,\pa_q
\tr U-\frac{1}{2} \pa_q U_{\alpha\beta}\, \pa_q U^{\alpha\beta}$.
The asymptotic form of the wave coordinate condition
\eqref{eq:wavec} is \beq\label{eq:wavecasympt} 2\pa_q
U_{L\mu}=L_\mu \,\pa_q \operatorname{tr} U,\qquad \forall
\mu=0,...,3, \eq where $U_{L\mu}= m^{\alpha\alpha^\prime}
U_{\alpha^\prime\mu} L_\alpha$ and
 $\tr U=m^{\alpha\beta}U_{\alpha\beta}$
The solution of the system
\eqref{eq:Einstasympt}-\eqref{eq:wavecasympt} exists globally and,
thus, the Einstein vacuum equations \eqref{Eired} in wave
coordinates satisfies the weak null condition. Moreover, the
component $\pa_q U_{\Lb\Lb}$ grows at most as $s$ while the
remaining components are uniformly bounded.
\end{theorem}
The asymptotic form \eqref{eq:Einstasympt} follows by a direct
calculation from \eqref{Eired}.
Observe that the null form
$Q_{\mu\nu}(\pa h, \pa h)$
disappears after passage to the asymptotic system.

Next we note that \eqref{eq:wavecasympt} is preserved under the flow of
\eqref{eq:Einstasympt}.
Contracting \eqref{eq:Einstasympt} with $L^\mu L^\nu$ we obtain
$$
(2\pa_s-U_{LL}\pa_q)\pa_{q} U_{LL} =0,
$$
which can be solved globally. More generally, contracting
\eqref{eq:Einstasympt} with the vector fields $\{L, S_1, S_2\}$ we
obtain \beq\label{eq:Einstasympt2} (2\pa_s-U_{LL}\pa_q)\pa_q
U_{TU}=0,\qquad\text{if}\quad T\in\{L,S_1, S_2\}\quad\text{and}\quad
U\in\{L, \Lb, S_1, S_2\}, \eq which can be solved globally now that
$U_{LL}$ has been determined. Note that the components $\pa_{q}
U_{TU}$ are constant along the integral curves of the vector field
$2\pa_{s} -U_{LL}\pa_{q}$. The remaining unknown component
$U_{\Lb\Lb}$ can be determined by contracting the equation
\eqref{eq:Einstasympt} with the vector field $\Lb$. \beq
\label{eq:LbLb} (2\pa_s-U_{LL}\pa_q)\pa_q U_{\Lb\Lb}= 4 P(\pa_{q}
U, \pa_{q} U), \eq By Lemma \ref{nullframeP} the quantity $P(\pa_q
U, \pa_q U)$ does not contain the term $(\pa_q U_{\Lb\Lb})^2$.
Thus, the equation \eqref{eq:LbLb} can be solved globally and
produces solutions growing exponentially in $s$. A more precise
information can be obtained from the asymptotic form of the wave
coordinate condition \eqref{eq:wavecasympt}. For contracting it
with
 the null frame
$\{L, S_1, S_2\}$ we obtain $\pa_q U_{LT}=0$, if $T\in\{L,S_1, S_2\}$.
Therefore, \beq \label{eq:asymptnullframeP} P(\pa_q U,\pa_q U)=
-\frac{1}{4}\delta^{AB}\delta^{A^\prime B^\prime} \Big( 2 \pa_q
U_{AA^\prime} \,\pa_q U_{BB^\prime} -\pa_q U_{AB} \,\pa_q
U_{A^\prime B^\prime}\Big) -\frac{1}{2}\delta^{AB}\pa_q U_{AB} \,
\pa_q U_{L\Lb}, \eq
It follows from \eqref{eq:Einstasympt2} that $P$ is already
determined and is, in fact, constant along the characteristics of
the field $2\pa_{s} -U_{LL}\pa_{q}$. Therefore, integrating
\eqref{eq:LbLb} we infer that $\pa_q U_{\Lb\Lb}$ grows at most
like $s$.

\section {Vector fields and commutators}\label{section:commute}
Let $Z\in {\cal Z}$ be any of the vector fields
$$
\Omega_{\a\b} = -x_{\a}\pa_{\b} + x_{\b}\pa_{\a},\qquad
 S = t\pa_{t} + r\pa_{r},\qquad \pa_\alpha,
$$
where $x_0=-t$ and $x_i=x^i$, for $i\geq 1$.
Let $I=(\iota_1,...,\iota_k)$, where $|\iota_i|=1$, be an ordered
multiindex of length $|I|=k$ and let $Z^I=Z^{\iota_1}\cdot\cdot\cdot Z^{\iota_k}$ denote a product of $|I|$
such derivatives. With a slight abuse of notation we will also identify
the index set with vector fields, so $I=Z$ means the index $I$ corresponding to
the vector field $Z$. Furthermore, by a sum over $I_1+I_2=I$ we mean a sum
over all possible order preserving
partitions of the ordered multiindex $I$ into two
ordered multiindices $I_1$ and $I_2$, i.e.
if $I=(\iota_1,...,\iota_k)$, then
$I_1=(\iota_{i_1},...,\iota_{i_n})$
and $I_2=(\iota_{i_{n+1}},...,\iota_{i_k})$,
where $i_1,...,i_k$ is any reordering of the integers $1,...,k$ such that $i_1<...<i_n$ and $i_{n+1}<...<i_k$
and $i_1,...,i_k$.
 With this convention
Leibniz rule becomes $Z^I(fg)=\sum_{I_1+I_2=I} (Z^{I_1} f)(Z^{I_2}
g)$. We denote by $\bar\pa$ the tangential derivatives, i.e.,
$\bar\pa=\{\bar\pa_0,\bar\pa_1,\bar\pa_2,\bar\pa_3\}$ and note 
that the span of the tangential derivatives $\{\bar\pa_0,\bar\pa_1,\bar\pa_2,\bar\pa_3\}$ coincides with the 
linear span  of the vectorfields $\{\pa_{s},
\pa_{S_1},\pa_{S_2}\}$.
\begin{lemma}\label{lem:coord}
We have the following expressions for the coordinate vector
fields:
\begin{align}
&\pa_{t} = \frac {tS - x^{i}\Omega_{0i}}{t^{2}-r^{2}},
\label{eq:tZ}\\
&\pa_{r} = \omega^i\pa_i =\frac {t\omega^{i}\Omega_{0i} -
rS}{t^{2}-r^{2}},\label{eq:rZ}\\
&\pa_{i} = \frac {-x^{j}\Omega_{ij} + t \Omega_{0i} - x_{i} S}{t^{2}-r^{2}}=
-\frac {x_{i} S} {t^{2}-r^{2}} + \frac {x_{i}
x^{j}\Omega_{0j}}{t(t^{2}-r^{2})} + \frac {\Omega_{0i}}t
\label{eq:iZ}
\end{align}
In particular,
\begin{equation}
\label{eq:somega}
\pa_{s} =\frac{1}{2}\big(\pa_t+\pa_r\big)
= \frac {S + \omega^{i}\Omega_{0i}}{2(t+r)},\qquad
\pab_{i}= \pa_{i}- \omega_{i}\pa_{r}
=\frac{\omega^j\Omega_{ij}}{r} = \frac {-\omega_{i}\omega^{j}\Omega_{0j} +
\Omega_{0i}}t.
\end{equation}
\end{lemma}
\begin{lemma}\label{tanderZ} For any function $f$ we have the estimate
\beq \label{eq:tanZ}
(1+t+|q|)|\bar \pa  f|+(1+|q|)|\pa f|\les C \sum_{|I|=1}|Z^I f|,\qquad \quad
|\pa f|\les |\bar\pa f|+|\pa_q f|
\eq
where $|\bar\pa f|^2=
|\bar\pa_0 f|^2+|\bar\pa_1 f|^2+|\bar\pa_2 f|^2+|\bar\pa_3 f|^2$ and $\bar\pa_0=\pa_s$. Furthermore
\beq\label{eq:2tanZ}
|\bar\pa^2 f|\les \frac{C}{r}\sum_{|I|\leq 2}
\frac{|Z^I f|}{1+t+|q|},
\eq
where $|\bar\pa^2 f|^2=\sum_{\alpha,\beta=0,1,2,3}
|\bar\pa_\alpha\bar\pa_\beta f|^2$.

Moreover, if $k^{\alpha\beta}$ is a symmetric tensor then
\beq \label{eq:derframeZ}
|k^{\alpha\beta} \pa_\alpha\pa_\beta \phi|\leq
C\bigg (\frac{|k|}{1+t+|q|}+ \frac{|k|_{LL}}{1+|q|}\bigg )
\sum_{|I|\le 1}|\pa Z^I \phi|
\eq
\end{lemma}
\begin{proof} First we note that if $r+t\leq 1$ then
\eqref{eq:tanZ} holds since the usual derivatives $\pa_\alpha$ are
included in the sum on the right. The inequality for $|\bar\pa f|$
in \eqref{eq:tanZ} follows directly from \eqref{eq:somega}; one
has to divide into two cases $r\leq t$ and $r\geq t$ and use
two different expressions depending on the relative size of $r$ and $t$.
The inequality for $|\pa f|$ in \eqref{eq:tanZ} follows from
\eqref{eq:tZ} and the first identity in \eqref{eq:iZ}.

If $t+r<1$ then \eqref{eq:2tanZ} follows from  \eqref{eq:somega}, since $|\pa_i\omega_j|\leq C
r^{-1}$ and the sum on the right of \eqref{eq:2tanZ} contains the
usual derivatives. Since $|\Omega_{ij}\omega_k|\leq C$ and
$\Omega_{ij}r=\Omega_{ij}t=0$, for $1\leq i,j\leq 3$ it follows,
by applying $\bar\pa_i=r^{-1}\omega^j\Omega_{ij}$ to the
expressions in \eqref{eq:somega},  that
\beq \label{eq:spacetanZ}
|\bar\pa_i \bar\pa_\alpha f|\leq Cr^{-1} (t+r)^{-1} \sum_{|I|\leq
2}|Z^I f|. \eq
Once again we  distinguish  the cases
$r<t$ and $r>t$ and use different expressions for $\bar\pa_i$.
 With the notation
$\bar\pa_0=2\pa_s$ \eqref{eq:spacetanZ} holds also for $\alpha=0$.
Since $[\pa_s,\bar\pa_i]=0$ it only remains to prove
\eqref{eq:2tanZ} for $\pa_s^2$. Since $S\omega^j=0$,
$|\Omega_{0i}\omega^j|\leq C t r^{-1}$, $S(t+r)=2(t+r)$ and
$|\Omega_{0i}(t+r)| \leq C(t+r)$ ,\eqref{eq:2tanZ} follows also
for $\pa_s^2$.

The inequality \eqref{eq:derframeZ} follows from
Lemma \ref{tander2}, \eqref{eq:tanZ} and the commutator identity
$[Z,\pa_i]=a_i^j\pa_j$.
\end{proof}

\begin{lemma}\label{waveeqframe} Suppose $\Boxr_g \phi=F$.
Then \beq\label{eq:waveoper} \Big|\Big(4\pa_s
-\frac{H_{LL}}{2g^{L\Lb}}\pa_q -\frac{\overline{\operatorname{tr}}\,
H \,+H_{L\Lb}}{2g^{L\underline{L}}\,\, r} \Big) \pa_q
(r\phi)+\frac{rF}{2g^{L\underline{L}}}\Big| \les \Big(1+
\frac{r\,|H|_{\cal LT}}{1+|q|}+ |H|\Big) r^{-1} \sum_{|I|\le 2}|Z^I
\phi| \eq
\end{lemma}
\begin{proof} By Corollary \ref{cor:BoxH}
\begin{multline*}
\Big|\Big(4\pa_s -\frac{H_{LL}}{2g^{L\Lb}}\pa_q
-\frac{\overline{\operatorname{tr}}\, H \,+H_{L\Lb}}
 {2g^{L\underline{L}}\,\, r}\Big) \pa_q
(r\phi)+\frac{r F}{2g^{L\underline{L}}}\Big|
\\
\les r|\triangle_\omega \phi|+
r |H|_{L{\cal T}}|\pab\pa\phi|+ |H|\,\big( r\,|\pab^2 \phi| +
|\bar\pa\phi|+r^{-1}|\phi|\big)
\end{multline*}
where $\triangle_\omega=\delta^{ij}\bar\pa_i \bar\pa_j$. Here all
the the derivatives can be reexpressed in terms of the vector
fields $Z$ and $\pa_q$ using \ref{tanderZ}, yielding the
expression \eqref{eq:waveoper}. Note that
$$
|\bar\pa \pa\phi|
\les \frac{\sum_{|I|=1}|Z^I\pa \phi|}{1+t+|q|}
\les \frac{\sum_{|I|\leq 1} |\pa Z^I \phi|}{1+t+|q|}
\les \frac{\sum_{|I|\leq 2} |Z^I\phi|}{(1+|q|)(1+t+|q|)}.
$$
\end{proof}
\begin{lemma} \label{commut1}
Let $Z=Z^\mu\pa_\mu$ be any of the vector fields  above and let $c_\alpha^\mu$
be defined by
$$
[\pa_\alpha,Z]=c_\alpha^{\,\,\, \mu}\pa_\mu ,\qquad c_\alpha^{\,\,\, \mu}=\pa_\alpha Z^\mu
$$
Then $c_\alpha^\mu$ are constants and
$$
c_{LL}=c^{\Lb\Lb}=0.
$$
Furthermore
$$
[Z,\square]=-c_Z\square
$$
where $c_Z$ is either $0$ or $2$.

In addition, if $Q$ is a null form, then
\beq\label{eq:nullformcommut}
Z Q
(\pa \phi,\pa \psi)=Q(\pa \phi,\pa Z \psi)+Q(\pa Z \phi,\pa \psi)
+ \tilde Q (\pa \phi,\pa \psi)
\eq
for some null form $\tilde Q$ on the right hand-side.
\end{lemma}
\begin{proof}
Since $Z=Z^\alpha\pa_\alpha$ is a Killing or conformally Killing vector field we have
\beq\label{eq:confkill}
\pa_{\a} Z_{\beta} + \pa_{\beta} Z_{\alpha} = f m_{\a\b}
\eq
where $Z_\alpha=m_{\alpha\beta} Z^\beta$.
In fact, for the vector fields above, $f=0$ unless $Z=S$ in which case $f=2$.
In particular,
$$
L^{\alpha} L^{\beta}\pa_\alpha Z_\beta =0.
$$
If $c_{\alpha}^{\,\,\,\mu}$ is  as defined above and
$c_{\alpha\beta} = c_{\alpha}^{\,\,\,\mu}m_{\mu\beta}=\pa_\alpha Z_\beta$
the above simply means that
$c_{LL}=c^{\Lb\Lb}=0$. which proves the first part of the lemma.
To verify \eqref{eq:nullformcommut} we first consider the null form
$Q=Q_{\a\b}$
We have
\begin{align*}
Z Q_{\a\b}(\pa\phi,\pa\psi) &= Q_{\a\b}(\pa Z \phi,\pa\psi) +
 Q_{\a\b}(\pa\phi,\pa Z \psi) \\ &+ [Z,\pa_\a] \phi \pa_\b \psi -
\pa_\b \phi [Z,\pa_\a]\psi + [Z,\pa_\b] \phi \pa_\a \psi -
\pa_\a \phi [Z,\pa_\b]\psi\\ & =
 Q_{\a\b}(\pa Z \phi,\pa\psi) +
 Q_{\a\b}(\pa\phi,\pa Z \psi) - c_\a^\mu ( \pa_\mu \phi \pa_\b \psi
- \pa_\b\phi \pa_\mu \psi ) - c_\b^\mu ( \pa_\mu \phi \pa_\a \psi
- \pa_\a\phi \pa_\mu \psi)\\ & =
 Q_{\a\b}(\pa Z \phi,\pa\psi) +
 Q_{\a\b}(\pa\phi,\pa Z \psi) - c_\a^\mu Q_{\mu\b} (\pa\phi,\pa\psi)-
c_\b^\mu Q_{\mu\a} (\pa\phi,\pa\psi)
\end{align*}
The calculation for the null form $Q_0(\pa\phi,\pa\psi) =
m^{\a\b} \pa_\a\phi \pa_\b \psi$ proceeds as follows:
\begin{align*}
Z Q_0(\pa\phi,\pa\psi) &= Q_0(\pa Z \phi,\pa\psi) +
 Q_0 (\pa\phi,\pa Z \psi) + m^{\a\b} [Z,\pa_\a] \phi \pa_\b \psi  +
m^{\a\b} \pa_\a \phi [Z,\pa_\b] \psi\\ & =
 Q_0(\pa Z \phi,\pa\psi) +
 Q_0(\pa\phi,\pa Z \psi) + m^{\a\b} c_\a^\mu \pa_\mu \phi \pa_\b \psi
 + m^{\a\b} c_b^\mu \pa_\a\phi \pa_\mu \psi \\ & =
 Q_0(\pa Z \phi,\pa\psi) +
 Q_0(\pa\phi,\pa Z \psi) + f m^{\a\b} \pa_\a\phi\pa_\b\psi\\ &=
Q_0(\pa Z \phi,\pa\psi) +
 Q_0(\pa\phi,\pa Z \psi) + f Q_0 (\pa\phi,\pa\psi),
\end{align*}
where $f$ is a constant associated with a Killing (conf. Killing) vector field $Z$ via a relation
$c^{\a\b} + c^{\b\a} = fm^{\a\b}$.
\end{proof}

\begin{lemma} \label{commut2}
If $k^{\alpha\beta}$ is a symmetric tensor then
\beq
k^{\alpha\beta}[\pa_\alpha\pa_\beta, Z]=k_Z^{\alpha\beta}\pa_\alpha\pa_\beta,
\qquad\text{where}\qquad
k_Z^{\alpha\beta}=k^{\alpha\gamma} c_{\gamma}^{\,\,\, \beta}
+k^{\gamma\beta}c_{\gamma}^{\,\,\,\alpha},
\qquad c_{\alpha}^{\,\,\,\mu} =\pa_\alpha Z^\mu .
\eq
In particular $k^{\alpha\beta}_S=2k^{\alpha\beta}$ and
\beq \label{eq:kLT}
|k_Z|_{LL}\leq 2|k|_{L\cal T}.
\eq
In general
\beq\label{eq:kformula}
[k^{\alpha\beta}\pa_\alpha\pa_\beta, Z^I]=\sum_{I_1+I_2=I, \,|I_2|<|I|}
k^{I_1\alpha\beta}\pa_\alpha\pa_\beta Z^{I_2},
\eq
where
\beq
k^{J\alpha\beta}= \!\sum_{ |K|\leq |J|}c^{J\alpha\beta}_{K\mu\nu}
Z^K k^{\mu\nu}=
-Z^J k^{\alpha\beta}
-\!\sum_{K+Z=J} Z^K k_Z^{\alpha\beta}
+\!\!\sum_{ |K|\leq |J|-2} d^{J\alpha\beta}_{K\mu\nu} Z^K k^{\mu\nu}
\eq
for some constants
$c^{J\alpha\beta}_{M\mu\nu}$ and $d^{J\alpha\beta}_{M\mu\nu}$. Here the sum \eqref{eq:kformula} means the sum over all possible order preserving partitions of the ordered multiindex $I$
into two ordered multiindices $I_1$ and $I_2$.
\end{lemma}
\begin{proof}
First observe that since the vector fields $Z$ are linear in $t$
and $x$ we have
$$
[\pa^{2}_{\a\b}, Z] = [\pa_{\b}, Z] \pa_{\a} +
[\pa_{\a}, Z] \pa_{\b}=
c_{\beta}^{\,\,\gamma}\pa_\gamma\pa_\alpha
+c_{\alpha}^{\,\,\,\gamma} \pa_\gamma\pa_\beta.
$$
which proves the first statement and the second follows since
$c_{L}^{\,\,\Lb} =0$.

To prove \eqref{eq:kformula} we first write
$$
Z^I \big(k^{\alpha\beta}\pa_\alpha\pa_\beta \phi\big)
=\sum_{K+J=I} (Z^K k^{\alpha\beta}) Z^J
\big( \pa_\alpha\pa_\beta \phi\big)
$$
Then we observe that \beq\label{eq:Zcom2} Z^{J}
\pa_\alpha\pa_\beta\phi=\sum_{J_1+J_2=J,\,
J_1=(\iota_{11},...,\iota_{1n})}
\left[Z^{\iota_{11}},\left[Z^{\iota_{12}},\left[...,
\left[Z^{\iota_{1\,n-1}},[Z^{\iota_{1n}},
\pa^2_{\a\b}]\,\right]...\right] \,\right]\,\right] Z^{J_2}\phi,
\eq where the sum is over all order preserving partitions of the
ordered multiindex $J=(\iota_1,...,\iota_k)$ into two ordered
multiindices $J_1=(\iota_{11},...,\iota_{1n})$ and
$J_2=(\iota_{21},...,\iota_{2k})$. It therefore follows that
$$
k^{J\a\b} = - \sum_{K+L=J,\, L=(\iota_1,...,\iota_l)}
(Z^K k^{\a\b}) \left[Z^{\iota_1},\left[Z^{\iota_2},\left[...,
\left[Z^{\iota_{l-1}},[Z^{\iota_l}, \pa^2_{\a\b}]\,\right]...\right]
\,\right]\,\right]
$$
The desired representation follows after taking into account that
$$
(Z^K k^{\a\b}) [Z,\pa^2_{\a\b}] = - (Z^K k^{\a\b}_Z)\pa_\alpha\pa_\beta
$$
\end{proof}

\begin{cor}\label{commut3}
 Let $\Boxr_g=\Box+H^{\alpha\beta}\pa_\alpha\pa_\beta$.
 Then with $\hat{Z}=Z+c_Z$
\beq \label{eq:curvwaveeqcommut1} \Boxr_g Z \phi- \hat{Z}\Boxr_g
\phi= -({\hat Z} H^{\alpha\beta}+
H_{{Z}}^{\alpha\beta})\pa_\alpha\pa_\beta \phi, \eq As a
consequence, we have \beq \label{eq:curvwaveeqcommutest1}
\big|\Boxr_g Z \phi- \hat{Z}\Boxr_g \phi|\les \bigg (\frac {|Z
H|+|H|}{1+t+|q|} +\frac {|Z H|_{LL}+|H|_{L\cal T}}{1+|q|}\bigg )
\sum_{|I|\le 1}|\pa Z^I \phi| \eq In general \beq
\label{eq:curvwaveeqcommut2} \Boxr_g Z^I \phi- \hat{Z}^I\Boxr_g
\phi= -\sum_{I_1+I_2=I, \, |I_2|<|I|}\hat H^{{I}_1\alpha\beta}
\pa_\alpha\pa_\beta Z^{I_2}\phi, \eq where \beq \hat
H^{J\alpha\beta}= \!\sum_{ |M|\leq |J|}
{c}^{J\alpha\beta}_{M\mu\nu} \hat{Z}^M H^{\mu\nu}= -\hat Z^J
H^{\alpha\beta} - \!\sum_{M+Z=J } \hat Z^M H_{Z}^{\alpha\beta}
+\!\!\sum_{ |M|\leq |J|-2} {d}^{J\alpha\beta}_{M\mu\nu} \hat Z^M
H^{\mu\nu} \eq We have
\begin{align}
|\Boxr_g Z^I \phi|\les |\hat{Z}^I \Boxr_g \phi| +\frac 1{1+t+|q|}
\,\,\,\sum_{|K|\leq |I|,}\,\, \sum_{|J|+(|K|-1)_+\le |I|} \,\,\,
|Z^{J} H|\,\, {|\pa Z^{K} \phi|}
\label{eq:curvwaveeqcommutest2} \\
+ \frac 1{1+|q|}
 \sum_{|K|\leq |I|}\Big(\sum_{|J|+(|K|-1)_+\leq |I|} \!\!\!\!\!|Z^{J} H|_{LL}
+\!\!\!\!\!\sum_{|J^{\prime}|+(|K|-1)_+\leq |I|-1}\!\!\!\!\!|Z^{J^{\prime}}
H|_{L\cal T}
+\!\!\!\!\!\sum_{|J^{\prime\prime}|+(|K|-1)_+\leq |I|-2}\!\!\!\!\!
|Z^{J^{\prime\prime}} H|\Big) {|\pa Z^{K} \phi|}\nn
\end{align}
where $(|K|-1)_+=|K|-1$ if $|K|\geq 1$ and $(|K|-1)_+=0$ if $|K|=0$.
\end{cor}
\begin{proof}
First observe that
\begin{align*}
\hat Z \Boxr_g \phi &= (Z+ c_Z) \Box \phi +
(Z+c_Z)H^{\a\b}\pa^2_{\a\b} \phi
\\ &= \Box Z \phi + H^{\a\b} \pa^2_{\a\b} Z\phi +
(Z H^{\a\b} ) \pa^2_{\a\b} \phi + (H^{\a\b}_Z + c_Z
H^{\a\b})\pa^2_{\a\b}\phi \\ &= \Boxr_g Z\phi + (Z H^{\a\b} )
\pa^2_{\a\b} \phi + (H^{\a\b}_Z + c_Z H^{\a\b})\pa^2_{\a\b}\phi
\end{align*}
Recall now that the constant $c_Z$ is different from $0$ only in
the case of the scaling vector field $S$. Moreover, in that case
$$
H^{\a\b}_S + c_S H^{\a\b} =0
$$
The inequality \eqref{eq:curvwaveeqcommutest1} now follows from
\eqref{eq:curvwaveeqcommut1}, \eqref{eq:kLT}
and the estimate \eqref{eq:derframeZ}.
The general commutation formula \eqref{eq:curvwaveeqcommut2}
follows from the following
calculation, similar to the one in Lemma \ref{commut2}.
We have
$$
\hat Z^{I}\Boxr_{g}\phi = \hat Z^{I}\Box \phi + \hat Z^{I}
H^{\a\b}\pa^{2}_{\a\b}\phi = \Box Z^{I}\phi + \sum_{J+K=I} \hat
Z^{J} H^{\a\b} Z^{K}\pa^{2}_{\a\b}\phi
$$
If we now use \eqref{eq:Zcom2} we get
\eqref{eq:curvwaveeqcommut2} as in the proof of
Lemma \ref{commut2}. The inequality \eqref{eq:curvwaveeqcommutest2} now follows from
\eqref{eq:curvwaveeqcommut2}, \eqref{eq:kLT}
and the estimate \eqref{eq:derframeZ}.
\end{proof}

\section{Basic energy identities}\label{section:energywave}
We now establish basic energy identities for solutions of the
equation \beq \label{eq:quasinh} \Boxr_{g} \phi =F \eq We denote by
$\Si_{t}$ the hyper surfaces $t=$const, by $C_{t_{1}}^{t_{2}}(q)$
the forward light cones with a vertex at $(q,0)$ and truncated at
times $t_{1}, t_{2}$. We also denote by $K_{t_{1}}^{t_{2}}(q)$ the
interior of the light cone $C_{t_{1}}^{t_{2}}(q)$and by $B_{t,r}$
the ball of radius $r$ centered at $(t,0)$.
\begin{lemma}
Let $\phi$ be a solution of \eqref{eq:quasinh}. Then for
any $t_{1}\le t_{2}$ and an arbitrary $q\le t_{2}$
\begin{align}
\int_{\Si_{t_{2}}}\big (-g^{00}|\pa_{t}\phi|^{2}+
g^{ij}\pa_{i}\phi \pa_{j}\phi \big ) & =
\int_{\Si_{t_{1}}}\big (-g^{00}|\pa_{t}\phi|^{2}+
g^{ij}\pa_{i}\phi \pa_{j}\phi \big )\nn \\ &-
2\int_{t_{1}}^{t_{2}}\int_{\Si_{\tau}}
\Big ( \pa_{\a} g^{\a\b} \pa_{\b}\phi \pa_{t} \phi -
\frac 12 \pa_{t} g^{\a\b}\pa_{\a}\phi \pa_{\b}\phi +
F\pa_{t}\phi \Big ),\label{eq:eglobal}
\end{align}
and
\begin{align}
\int_{B_{t_{1}-q}}\big (-g^{00}|\pa_{t}\phi|^{2}+
g^{ij}\pa_{i}\phi \pa_{j}\phi \big ) & +
\int_{C_{t_{1}}^{t_{2}}(q)}|\pab\phi|^{2}  =
\int_{B_{t_{2}-q}}\big (-g^{00}|\pa_{t}\phi|^{2}+
g^{ij}\pa_{i}\phi \pa_{j}\phi \big )\nn \\ & +
2\int_{K_{t_{1}}^{t_{2}}(q)}
\Big ( \pa_{\a} g^{\a\b} \pa_{\b}\phi \pa_{t} \phi -
\frac 12 \pa_{t} g^{\a\b}\pa_{\a}\phi \pa_{\b}\phi +
F\pa_{t}\phi \Big )\nn \\ & +
2\int_{C_{t_{1}}^{t_{2}}(q)}
\big ( 2(g^{\a\b}-m^{\a\b}) L_\a \pa_{\b}\phi\pa_{t}\phi + (g^{\a\b}-m^{\a\b})
\pa_{\a}\phi\pa_{\b}\phi\big )\label{eq:elocal}
\end{align}
\end{lemma}
\begin{proof}
We multiply the equation \eqref{eq:quasinh} by $\pa_{t}\phi$
and integrate over the space-time slab between the hyper surfaces
$\Si_{t_{1}}$ and $\Si_{t_{2}}$. We have
\begin{align*}
-\int_{t_{1}}^{t_{2}}\int_{\Si_{\tau}} g^{\a\b}\pa^{2}_{\a\b}\phi
\pa_{t}\phi & = \int_{t_{1}}^{t_{2}}\int_{\Si_{\tau}} \big (g^{\a\b}\pa_{\b}\phi
\pa_{t}\pa_{\a}\phi + \pa_{\a} g^{\a\b}\pa_{\b}\phi
\pa_{t}\phi\big ) \\ &- \int_{\Si_{t_{2}}} g^{0\b}\pa_{\b}\phi
\pa_{t}\phi + \int_{\Si_{t_{1}}} g^{0\b}\pa_{\b}\phi
\pa_{t}\phi\\ & = \frac 12 \int_{\Si_{t_{2}}} \big (- g^{00}|\pa_{t}\phi|^{2}
+g^{ij}\pa_{i}\phi\pa_{j}\phi\big ) -
\frac 12 \int_{\Si_{t_{1}}} \big (- g^{00}|\pa_{t}\phi|^{2}
+g^{ij}\pa_{i}\phi\pa_{j}\phi\big )\\ & +
\int_{t_{1}}^{t_{2}}\int_{\Si_{\tau}} \big (\pa_{\a} g^{\a\b}\pa_{\b}\phi
\pa_{t}\phi -\frac 12 \pa_{t} g^{\a\b}\pa_{\a}\phi\pa_{\b}\phi \big )
\end{align*}
and the desired identity \eqref{eq:eglobal} follows.
Similarly, integrating over the region $K_{t_{1}}^{t_{2}}(q)$ we obtain
\begin{align*}
\int_{B_{t_{1}-q}}\big (-g^{00}|\pa_{t}\phi|^{2}+
g^{ij}\pa_{i}\phi \pa_{j}\phi \big ) & -
\int_{C_{t_{1}}^{t_{2}}(q)}\big ( 2g^{\a\b}L_\a \pa_{\b}\phi\pa_{t}\phi +
g^{\a\b}\pa_{\a}\phi\pa_{\b}\phi\big ) \\ & =
\int_{B_{t_{2}-q}}\big (-g^{00}|\pa_{t}\phi|^{2}+
g^{ij}\pa_{i}\phi \pa_{j}\phi \big )\\ & +
2\int_{K_{t_{1}}^{t_{2}}(q)}
\Big ( \pa_{\a} g^{\a\b} \pa_{\b}\phi \pa_{t} \phi -
\frac 12 \pa_{t} g^{\a\b}\pa_{\a}\phi \pa_{\b}\phi +
F\pa_{t}\phi \Big )
\end{align*}
Subtracting the Minkowski part from the metric $g$ in the
$C_{t_{1}}^{t_{2}}(q)$ integral leads to the identity
\eqref{eq:elocal}.
\end{proof}
\begin{cor}
\label{cor:Energy}
Let $\phi$ be a solution of the equation \eqref{eq:quasinh}
with a metric $g$ satisfying the condition that
\beq
\label{eq:esmall}
|H|\le \frac 14,\qquad H^{\a\b}=g^{\a\b}-m^{\a\b}.
\eq
Then for any $0<\gamma \le 1$
\begin{align}
 \int_{\Si_{t_{2}}} \big (|\pa_{t}\phi|^{2}
+|\nab_x\phi|^{2}\big ) & + \int_{t_{1}}^{t_{2}} \int_{\Si_{\tau}}
\frac {\gamma\,\,|\pab\phi|^{2}}{(1+|q|)^{1+2\gamma}}\le
 4 \int_{\Si_{t_{1}}} \big (|\pa_{t}\phi|^{2}
+|\nab_x \phi|^{2}\big ) \nn \\ &+ 8
\int_{t_{1}}^{t_{2}} \int_{\Si_{\tau}} \big |\pa_{\a} g^{\a\b}\pa_{\b}\phi
\pa_{t}\phi -\frac 12 \pa_{t} g^{\a\b}\pa_{\a}\phi\pa_{\b}\phi
+ F\pa_{t}\phi \big |
\label{eq:energy}\\ &+ 2\int_{t_{1}}^{t_{2}} \int_{\Si_{\tau}}
\frac {\gamma}{(1+|q|)^{1+2\ga}}\big |
(g^{\a\b}-m^{\a\b})\pa_{\a}\phi\pa_{\b}\phi +
2(g^{\Lb\b} - m^{\Lb\b})\pa_{\b}\phi\pa_{t}\phi \big |\nn
\end{align}
\end{cor}
\begin{proof} First we note that \eqref{eq:esmall} implies that
\begin{equation}\label{eq:esmallest}
\frac{3}{4}\big(|\pa_t\phi|^2+|\nabla_x\phi|^2\big)
\leq -g^{00}|\pa_t\phi|^2+g^{ij} \pa_i\phi\,\pa_j\phi
\leq \frac{5}{4}\big(|\pa_t\phi|^2+|\nabla_x\phi|^2\big)
\end{equation}
The inequalities \eqref{eq:elocal} and \eqref{eq:eglobal}imply that
\begin{align}
\int_{C_{t_{1}}^{t_{2}}(q)}|\pab\phi|^{2}  &\le
\int_{\Si_{t_2}}\big (-g^{00}|\pa_{t}\phi|^{2}+
g^{ij}\pa_{i}\phi \pa_{j}\phi \big )\label{eq:ezlocal} \\ & +
2\int_{K_{t_1}^{t_2}(q)}
\pa_{\a} g^{\a\b} \pa_{\b}\phi \pa_{t} \phi -
\frac 12 \pa_{t} g^{\a\b}\pa_{\a}\phi \pa_{\b}\phi +
F\pa_{t}\phi \nn \\ & +
2\int_{C_{t_{1}}^{t_{2}}(q)}
 2(g^{\a\b}-m^{\a\b}) L_\a \pa_{\b}\phi\pa_{t}\phi + (g^{\a\b}-m^{\a\b})
\pa_{\a}\phi\pa_{\b}\phi \\
&\le
\int_{\Si_{t_1}}\big (-g^{00}|\pa_{t}\phi|^{2}+
g^{ij}\pa_{i}\phi \pa_{j}\phi \big )\\ & +
2\int_{t_1}^{t_2}\int_{\Sigma_t}\Big|
\pa_{\a} g^{\a\b} \pa_{\b}\phi \pa_{t} \phi -
\frac 12 \pa_{t} g^{\a\b}\pa_{\a}\phi \pa_{\b}\phi +
F\pa_{t}\phi \Big|\nn \\ & +
2\int_{C_{t_{1}}^{t_{2}}(q)}
 \Big| 2(g^{\a\b}-m^{\a\b}) L_\a \pa_{\b}\phi\pa_{t}\phi + (g^{\a\b}-m^{\a\b})
\pa_{\a}\phi\pa_{\b}\phi\Big| \nn
\end{align}
We multiply the above inequality by an integrable factor
$\gamma(1+|q|)^{-1-2\ga}$ and integrate with
respect to $q$ in the interval $(-\infty, t_2]$ to obtain:
\begin{align}
\int_{t_1}^{t_2}\int_{\Sigma_t}
\frac{\gamma|\pab\phi|^{2}}{(1+|q|)^{1+2\gamma}}
&\leq \frac{5}{4}\int_{\Sigma_{t_1}}
\big(|\pa_t\phi|^2+|\nabla_x\phi|^2\big)\nn
\\ &+ 2
\int_{t_{1}}^{t_{2}} \int_{\Si_{\tau}} \big |\pa_{\a} g^{\a\b}\pa_{\b}\phi
\pa_{t}\phi -\frac 12 \pa_{t} g^{\a\b}\pa_{\a}\phi\pa_{\b}\phi
+ F\pa_{t}\phi \big |
\label{eq:energy}\\ &+ 2\int_{t_{1}}^{t_{2}} \int_{\Si_{\tau}}
\frac {\gamma}{(1+|q|)^{1+2\ga}}\big |
(g^{\a\b}-m^{\a\b})\pa_{\a}\phi\pa_{\b}\phi +
2(g^{\Lb\b} - m^{\Lb\b})\pa_{\b}\phi\pa_{t}\phi \big |\nn
\end{align}
where we also used \eqref{eq:esmallest}. On the other hand using \eqref{eq:esmallest} and \eqref{eq:eglobal} yields
\begin{align}
\int_{\Sigma_{t_1}}
\big(|\pa_t\phi|^2+|\nabla_x\phi|^2\big)&\leq
\frac{5}{3} \int_{\Sigma_{t_1}}
\big(|\pa_t\phi|^2+|\nabla_x\phi|^2\big)\\
&+\frac{8}{3}\int_{t_{1}}^{t_{2}} \int_{\Si_{\tau}} \big |\pa_{\a} g^{\a\b}\pa_{\b}\phi
\pa_{t}\phi -\frac 12 \pa_{t} g^{\a\b}\pa_{\a}\phi\pa_{\b}\phi
+ F\pa_{t}\phi \big |
\end{align}
and the corollary follows.
\end{proof}

\section{Poincar\'e and Klainerman-Sobolev inequalities}
We now state the following useful version of the Poincar\'e
inequality.
\begin{lemma}\label{le:Poinc}
Let $f$ be a smooth function.
Then for any $\ga >-1/2, \ga\ne 1/2$ and any positive $t$
\beq
\label{eq:Poinc}
\il_{R^{3}} \frac {|f(x)|^{2}\, dx}{(1+|t-r|)^{2+2\ga}}
\le C\il_{S_{(t+1)}} |f|^{2} \, dS+
C\il_{R^{3}} \frac {|\pa_{r} f(x)|^{2}\, dx}{(1+|t-r|)^{2\ga}}
\ee
provided that the left hand side is bounded.
Here $S_{(t+1)}$ is the sphere of radius $t+1$
and $r=|x|$.
\end{lemma}
\begin{proof}
Using polar coordinates $x=r\omega$ we write
$$
|f(r,\omega)|^{2} - |f(t+1,\omega)|^{2} =
-2 \int_{r}^{t+1} \pa_{r} f(\rho,\omega )\cdot f(\rho,\omega)
\,d\rho
$$
Hence
$$
|f(r,\omega)|^{2}r^2 \les  |f(t+1,\omega)|^{2}(t+1)^2
+2\int_{r}^{t+1} |\pa_{r} f(\rho,\omega )|\,
|f(\rho,\omega)|\, \rho^2\,d\rho,\qquad\text{if}\quad
r\leq t+1.
$$
Therefore multiplying  by $(1+|t-r|)^{-2-2\gamma}$ and
integrating with respect to $r$ from $0$ to $t+1$:
\begin{align*}
\int_{0}^{{t+1}}\!\!
\frac {|f(r,\omega)|^{2} \, r^2\, dr}{(1+|t-r|)^{2+2\ga}}
&\les \int_0^{t+1} \!
\frac{|f(t+1,\omega)|^{2}\, (t+1)^2\, dr}{(1+|t-r|)^{2+2\gamma}}
+ \int_{0}^{t+1}\int_r^{t+1}\frac{|\pa_r f(\rho,\omega)|
|f(\rho,\omega)|}{(1+|t-r|)^{2+2\gamma}}
\,\rho^2\,  d\rho  \, d r
\nn\\
&\les
|f(t+1,\omega)|^{2}\, (t+1)^2
+ \int_{0}^{t+1}\int_0^\rho\frac{|\pa_r f(\rho,\omega)|
|f(\rho,\omega)|}{(1+|t-r|)^{2+2\gamma}}
\, dr \, \rho^2\, d \rho
\nn\\
&\les
|f(t+1,\omega)|^{2}\, (t+1)^2
+ \int_{0}^{t+1}\frac{|\pa_r f(\rho,\omega)|
|f(\rho,\omega)|}{(1+|t-\rho|)^{1+2\gamma}}
 \, \rho^2\, d \rho
\nn\\
&\les  |f(t+1,\omega)|^{2}(t+1)^2
+\Big(\int_{0}^{t+1}
\frac{|\pa_r f(\rho,\omega)|^2\, \rho^2\,d\rho}
{(1+|t-\rho|)^{2\gamma}}\Big)^{1/2}
\Big(\int_{0}^{t+1}
\frac{|f(\rho,\omega)|^2\, \rho^2\,d\rho}
{(1+|t-\rho|)^{2+2\gamma}}\Big)^{1/2},
\end{align*}
where we first changed the order of integration and then
used Cauchy-Schwarz inequality.
It therefore follows that
$$
\int_{0}^{{t+1}}\!\!
\frac {|f(r,\omega)|^{2} \, r^2\, dr}{(1+|t-r|)^{2+2\ga}}
\les |f(t+1,\omega)|^{2}\, (t+1)^2
+\int_{0}^{t+1}
\frac{|\pa_r f(\rho,\omega)|^2\, \rho^2\,d\rho}
{(1+|t-\rho|)^{2\gamma}}
$$
and if we also integrate over the angular variables we get
$$
\il_{|x|\le (t+1)} \frac {|f(x)|^{2}\, dx}{(1+|t-r|)^{2+2\ga}}
\les \il_{S_{(t+1)}} |f|^{2} dS+
\il_{|x|\le (t+1)} \frac {|\pa_{r} f(x)|^{2}\, dx}{(1+|t-r|)^{2\ga}}
$$
On the other hand, if we instead integrate from $t+1$ to $2(t+1)$ we similarly obtain
\begin{align*}
\int_{t+1}^{{2(t+1)}}\!\!
\frac {|f(r,\omega)|^{2} \, r^2\, dr}{(1+|t-r|)^{2+2\ga}}
&\les \int_{t+1}^{2(t+1)} \!
\frac{|f(t+1,\omega)|^{2}\, (t+1)^2\, dr}{(1+|t-r|)^{2+2\gamma}}
+ \int_{t+1}^{2(t+1)}\!\!\int_{t+1}^{r}\frac{|\pa_r f(\rho,\omega)|
|f(\rho,\omega)|}{(1+|t-r|)^{2+2\gamma}}
\, \rho^2\,d\rho  \,  d r
\nn\\
&\les
|f(t+1,\omega)|^{2}\, (t+1)^2
+ \int_{t+1}^{2(t+1)}\int_\rho^{2(t+1)}\frac{|\pa_r f(\rho,\omega)|
|f(\rho,\omega)|}{(1+|t-r|)^{2+2\gamma}}
\, dr \, \rho^2\, d \rho
\nn\\
&\les
|f(t+1,\omega)|^{2}\, (t+1)^2
+ \int_{t+1}^{2(t+1)}\frac{|\pa_r f(\rho,\omega)|
|f(\rho,\omega)|}{(1+|t-\rho|)^{1+2\gamma}}
 \, \rho^2\, d \rho
\end{align*}
and as before it follows that
$$
\il_{(t+1)\le |x|\le 2(t+1)} \frac {|f(x)|^{2}\, dx}{(1+|t-r|)^{2+2\ga}}
\les \il_{S_{(t+1)}} |f|^{2} \, dS+
\il_{(t+1)\le |x|\le 2(t+1)} \frac {|\pa_{r} f(x)|^{2}\, dx}{(1+|t-r|)^{2\ga}}
$$
Finally, in the region $r\ge 2(t+1)$ the estimate \eqref{eq:Poinc}
would follow from the Hardy type inequality: \beq\label{eq:hardy}
\il_{|x|\ge (t+1)} \frac {|f(x)|^{2}\, dx}{|x|^{2+2\ga}} \le
\il_{|x|\ge (t+1)} \frac {|\pa_{r} f(x)|^{2}\, dx }{|x|^{2\ga}} +
(t+1)^{-1-2\ga} \il_{S_{(t+1)}} |f|^{2}\, dS, \eq that hold
provided the left hand side is bounded. One can for the proof
assume that $f$ ha compact support since we can choose a sequence
of compactly supported functions converging to a given function
$f$ in the norm defined by the right hand side as long as the norm
in the left of $f$ is bounded. \eqref{eq:hardy} for compactly
supported smooth functions can be easily seen from integrating the
identity
$$
\pa_r \left(\frac{r^2 f^2}{r^{1+2\gamma}}\right)
=\frac {2 r^{2}}{r^{1+2\ga}} f\cdot \pa_{r} f
 +(1-2\ga)\frac {r^{2}}{r^{2+2\ga}} f^{2},\qquad \gamma\neq -1/2
$$
from $r=t+1$ to $r=\infty$ and using Cauchy-Schwarz as above.
\end{proof}

We now state the global Sobolev inequality, which is
due to S. Klainerman \cite{K1}.
\begin{prop}\label{globalsobolev}
The following inequality holds for an arbitrary smooth
function $\phi$.
$$
|\phi(t,x)|(1+t+|t-r|)(1+|t-r|)^{1/2} \leq
C\sum_{|I|\leq 3} \|Z^I \phi(t,\cdot)\|_{L^2}.
$$
\end{prop}

\section
{Decay estimates for the wave equation on a curved space time  }
\label{section:decaywaveeq}
In this section we will derive some basic estimates for the scalar
wave equation on a curved background. The results will require some
weak assumptions on the metric $g$, which will be easily verified in
the case of a metric satisfying the reduced Einstein equations.

We consider the reduced
scalar wave equation: \beq\label{eq:waveeqscal} \Boxr_g \phi=F .\eq
The following result is a generalization of the lemma in \cite{L1}
to the variable coefficient case:

\begin{lemma} \label{decaywaveeq1} Suppose that $\phi$ satisfies the
reduced scalar wave equation \eqref{eq:waveeqscal} on a curved background
with a metric $g$. Suppose that
$H^{\alpha\beta}=g^{\alpha\beta}-m^{\alpha\beta}$ satisfies
\beq\label{eq:decaymetric1}
|H|\leq \frac{1}{4},
\qquad  \text{and}\qquad
|H|_{L{\cal T}}\leq  \frac{1}{4}\frac{|q|+1}{1+t+|x|}.
\eq
when $t/2\leq |x|\leq 2t$ and
\beq\label{eq:decaymetric2}
\int_0^{\infty}\!\! \|\,H(t,\cdot)\|_{L^\infty(D_t)}\frac{dt}{1+t}
\leq \frac{1}{4},\qquad\text{where}
\qquad D_t=\{x\in\bold{R}^3;\, t/2\leq |x|\leq 2t\}.
\eq
Then for any $t\ge 0$ and $x\in \R^3$, 
\begin{multline}\label{eq:decaywaveeq1}
(1+t+|x|)\,|\pa \phi(t,x)|
\leq  C\sup_{0\leq \tau\leq t}\!
\sum_{|I|\leq 1}\!
\|Z^I\! \phi(\tau,\cdot)\|_{L^\infty}\\
+ C\int_0^t
\Big((1+\tau)\| F(\tau,\cdot)\|_{L^\infty(D_\tau)}
+\sum_{|I|\leq 2}
(1+\tau)^{-1} \| Z^I \phi(\tau,\cdot)\|_{L^\infty(D_\tau)}\Big)\, d\tau
\end{multline}
\end{lemma}
\begin{proof} Since by Lemma \ref{tanderZ}
\beq\label{eq:awaycone} (1+|t-r|) |\pa\phi|+(1+t+r)|\bar\pa\phi|
\leq C\sum_{|I|=1}|Z^I\phi|,\qquad r=|x|, \eq the inequality
\eqref{eq:decaywaveeq1} holds when $r<t/2+1/2$ or $r>2t-1$.
Furthermore, since \beq\label{eq:qder} (1+r)|\pa_q\phi|\leq
C|\pa_q(r\phi)|+C|\phi|,\qquad r\geq 1 \eq it follows that
\beq\label{eq:qrder}
(1+t+r) |\pa\phi|\leq C\sum_{|I|\leq
1}|Z^I\phi|+C |\pa_q (r\phi)|
\eq Hence it suffices to prove that
$|\pa_q (r\phi)|$ is bounded by the right hand side of
\eqref{eq:decaywaveeq1} when $t/2+1/2<r<2t-1$. By Lemma
\ref{waveeqframe}
\beq\label{eq:waveoper5}
 \big|(4\pa_s -
\frac{H_{LL}}{2g^{L\Lb}}\pa_{q})\pa_{q}( r \phi)\big| \les \Big(1+
\frac{r\,|H|_{\cal LT}}{1+|q|}+ |H|\Big) r^{-1}\sum_{|I|\le 2} |Z^I
\phi|+ |H|\, r^{-1}\, |\pa_q (r\phi) | + r |F| \eq and using the
decay assumptions \eqref{eq:decaymetric1} and \eqref{eq:qrder} we
get \beq\label{eq:waveoper6} \big|(4\pa_s
-\frac{H_{LL}}{2g^{L\Lb}}\pa_{q})\pa_{q}( r \phi)\big| \les
\frac{|H|}{1+t} |\pa_q (r\phi) |+\sum_{|I|\le 2} \frac{Z^I
\phi|}{1+t} + C (t+1)|F|,\qquad\text{when}\quad t/2+1/2 \leq r\leq
2t-1 \eq Along an integral curve $(t,x(t))$ of the vector field
$\pa_s+H^{\Lb \Lb}(2g^{L\Lb})^{-1}\pa_q$, contained in the region
$t/2+1/2\leq |x|\leq 2t-1$, we have the following equation for
$\psi=\pa_q(r\phi)$: \beq \label{eq:qui} \Big|\frac{
d}{dt}\psi\Big|\leq \hat{h} |\psi|+ f \eq where $\hat{h}=C|H|/(1+t)$
and $f= C t |F|+C\sum_{|I|\le 2}|Z^I \phi|/(1+t)$. Hence multiplying
\eqref{eq:qui} with the integrating factor $e^{-\hat{H}}$, where
$\hat{H}=\int \hat{h}(s)\, ds $ we get \beq \Big|\frac{
d}{dt}\Big(\psi e^{-\hat{H}}\Big)\Big|\leq  f e^{-\hat{H}}. \eq If
we integrate backwards along an integral curve from any point
$(t,x)$ in the set $t/2+1/2\leq |x|\leq 2t-1$ until the first time
the curve intersects the boundary of the set at $(\tau,y)$,
$|y|=\tau/2+1/2$ or $|y|=2\,\tau-1$,  we obtain
$$
|\psi(t,x)|
\leq \exp\Big({\int_{\tau}^t \|\hat{h}(\sigma,\cdot)\|_{L^\infty}\, d\sigma}\Big)
|\psi(\tau,y)|
+\int_{\tau}^t \exp\Big({\int_{\tau^\prime}^t \|\hat{h}(\sigma,\cdot)\|_{L^\infty}\,
 d\sigma}\Big)
\|f(\tau^\prime,\cdot)\|_{L^\infty}\, d\tau^\prime,
$$
where the $L^\infty$ norms are taken only over the set $t/1+1/2\leq |x|\leq 2t-1$.
(Note that any integral curve has to intersect either of the two boundaries
$r=t/2+1/2$ or $r=2t-1$ since the slope of the curve $x(t)$
has to be close to $1$ when $H_{LL}$ is small.)
The lemma now follows from taking the supremum over $x$ in the set
$t/2+1/2\leq |x|\leq 2t-1$,
using that on the cones $|y|=\tau/2+1/2$ or $|y|=2\tau-1$
we have that $|\psi|\leq Cr|\pa_q\phi|+C |\phi|\leq
C\sum_{|I|\leq 1} |Z^I\phi|$,
by \eqref{eq:awaycone}, and using that by
\eqref{eq:decaymetric2}
$\int_{0}^t \|\hat{h}(\sigma,\cdot)\|_{L^\infty}\, d\sigma
\leq \frac{1}{4}.$
\end{proof}

For second order derivatives we have an estimate which gives
a slightly worse decay:
\begin{lemma} \label{decaywaveeq2}   Let $\phi$ be a solution of the
reduced scalar wave equation on a curved background
with a metric $g$. Assume that
$H^{\alpha\beta}=g^{\alpha\beta}-m^{\alpha\beta}$ satisfies
\beq\label{eq:decaymetric3}
\sum_{|I|\leq 1}|Z^I H|\leq \frac{\tilde{\varepsilon}}{4},\qquad \text{and}\qquad
\sum_{|I|\leq 1}|Z^I H|_{LL}+|H|_{L{\cal T}}
\leq\frac{\tilde{\varepsilon}}{4}\, \frac{|q|+1}{1+t+|x|}.\eq
when $t/2\leq |x|\leq 2t$
for some $\tilde{\varepsilon}\leq 1$.
Then, for $t\ge 0$, $x\in \R^3$, we have 
\begin{multline}\label{eq:decaywaveeq2}
(1+t+|x|)\sum_{|I|\leq 1}\,|\pa Z^I \phi(t,x)|
\leq  C\!\sup_{0\leq \tau\leq t}\Big(\frac{1+t}{1+\tau}\Big)^{C\tilde{\varepsilon}}
\sum_{|I|\leq 2}\!\|Z^I \phi(\tau,\cdot)\|_{L^\infty}\Big)\\
+ C\int_0^t\Big(\frac{1+t}{1+\tau}\Big)^{C\tilde{\varepsilon}}
\Big(\sum_{|I|\leq 1}(1+\tau)\| Z^I F(\tau,\cdot)\|_{L^\infty(D_\tau)}
+\sum_{|I|\leq 3}(1+\tau)^{-1} \| Z^I \phi(\tau,\cdot)\|_{L^\infty(D_\tau)}\Big)\, d\tau,
\end{multline}
where $D_t=\{x\in \bold{R}^3;\, t/2\leq |x|\leq 2t\}$.
\end{lemma}
\begin{proof} First when $r<t/2$ or $r>t/2$ the lemma
trivially follows from \eqref{eq:awaycone} with $\phi$ replaced by
$Z\phi$ so it only remains to prove the lemma when $t/2\leq r\leq
2t$. We have \beq\label{eq:sqZ} \Boxr_g Z \phi=F_Z=\hat{Z} F
+\big(\Boxr_g Z\phi-\hat{Z}\Boxr_g\phi\big), \eq where by
\eqref{eq:curvwaveeqcommutest1} the additional commutator term can
be estimated by \beq \big|\Boxr_g Z \phi- \hat{Z}\Boxr_g \phi|\les
\bigg (\frac {|Z H|+|H|}{1+t+|q|} +\frac {|Z H|_{LL}+|H|_{L\cal
T}}{1+|q|}\bigg ) \sum_{|I|\le 1}|\pa Z^I \phi|\les
\frac{\tilde{\varepsilon}}{1+t+q} \sum_{|I|\le 1}|\pa Z^I \phi|,
\eq where we used the decay assumption \eqref{eq:decaymetric3}.
Furthermore with the help of \eqref{eq:qrder}, applied to $Z^I\phi$ in place
of $\phi$, we obtain \beq \big|\Boxr_g Z \phi- \hat{Z}\Boxr_g
\phi|\les \frac{\tilde{\varepsilon}}{(1+t+|q|)^2}
\Big(\sum_{|I|\le 1}\big|\pa _q(rZ^I \phi)\big| +\sum_{|I|\leq
2}|Z^I\phi|\Big), \eq

Hence by \eqref{eq:waveoper6} applied to \eqref{eq:sqZ} in place of
\eqref{eq:waveeqscal} we get \beq\label{eq:quie} |(4\pa_s
-\frac{H_{LL}}{2g^{L\Lb}}\pa_{q})\pa_{q}( r Z\phi)| \les
\sum_{|I|\le 3} \frac{|Z^I \phi|}{1+t} +
\frac{\tilde{\varepsilon}}{1+t} \sum_{|I|\le 1}|\pa_q(r Z^I \phi)| +
t (|Z F|+|F|) \eq when $t/2+1/2\leq r\leq 2t-1$. Therefore \beq
\Big|(4\pa_s - \frac{H_{LL}}{2g^{L\Lb}}\pa_{q}) \sum_{|I|\leq
1}|\pa_{q}( r Z^I\phi)|\Big| \leq \frac{C\tilde{\varepsilon}}{1+t}
\sum_{|I|\le 1}|\pa_q(r Z^I \phi)| +C\sum_{|I|\le 3} \frac{|Z^I
\phi|}{1+t}
 + Ct (|Z F|+|F|)
\eq The desired result follows
multiplying  \eqref{eq:quie} by the  factor $(1+t)^{-C\tilde{\varepsilon}}$
and integrating as in the proof of the previous lemma. Along an
integral curve we have the equation \beq \Big|\frac{d}{dt}
\Big(\psi (1+t)^{-C\tilde{\varepsilon}}\Big)\Big| \leq
(1+t)^{-C\tilde{\varepsilon}} f, \eq where \beq \psi=\sum_{|I|\leq
1}|\pa_q (Z^I\phi)|,\quad f=C(1+t)(|Z F|+|F|)+C \sum_{|I|\leq
3}\frac{|Z^I \phi|}{1+t} \eq The lemma now follows as in the proof
of Lemma \ref{decaywaveeq1}.
\end{proof}

We observe that similar estimates hold for a system
\beq\label{eq:waveeqsyst} {\Boxr}_g \phi_{\mu\nu}=F_{\mu\nu} \eq
In particular, in our case, certain components of $F_{\mu\nu}$
expressed in the null-frame will decay better than others and for
these components we will also get better estimates for
$\phi_{\mu\nu}$. Since the vector fields $L$ and $\Lb$ commute
with contractions of any of the vector fields $\{L,\Lb,S_1,S_2\}$
proofs of the preceding lemmas imply the following
result:

\begin{cor} \label{decaywaveeq3} Let $\phi_{\mu\nu}$ be a solution
of reduced wave equation system \eqref{eq:waveeqsyst} on a curved background
with a metric $g$. Assume that
$H^{\alpha\beta}=g^{\alpha\beta}-m^{\alpha\beta}$ satisfies
\beq\label{eq:decaymetric4}
\sum_{|I|\leq 1}|Z^I H|\leq \frac{\tilde{\varepsilon}}{4},\qquad \text{and}\qquad
\sum_{|I|\leq 1}|Z^I H|_{LL}+|H|_{L{\cal T}}
\leq\frac{\tilde{\varepsilon}}{4}\, \frac{|q|+1}{1+t+|x|}.\eq
when $t/2\leq |x|\leq 2t$, for some $\tilde{\varepsilon}\leq 1$ and
\beq\label{eq:decaymetric5}
\int_0^{\infty}\!\! \|\,H(t,\cdot)\|_{L^\infty(D_t)}\frac{dt}{1+t}
\leq \frac{\tilde{\varepsilon}}{4}.
\eq
where $D_t=\{x\in\bold{R}^3;\, t/2\leq |x|\leq 2t\}$.
 Then for any $U,V\in\{L,\Lb,S_1,S_2\}$ and any $t\ge 0$, $x\in \R^3$:
\begin{align}
(1+t+|x|)\, &\,|\pa \phi(t,x)|_{UV}
\leq  C\!\sup_{0\leq \tau\leq t} \sum_{|I|\leq 1}\!
\|Z^I\! \phi(\tau,\cdot)\|_{L^\infty}\label{eq:decaywaveeq3}\\
&+ C\int_0^t
\Big((1+\tau)\||F|_{UV}(\tau,\cdot)\|_{L^\infty(D_\tau)}
+\sum_{|I|\leq 2}
(1+\tau)^{-1}\, \| Z^I \phi(\tau,\cdot)\|_{L^\infty(D_\tau)}\Big)\, d\tau,\nn\\
(1+t+|x|)&\sum_{|I|\leq 1}\,|\pa Z^I \phi|(t,x)|
\leq  C\!\sup_{0\leq \tau\leq t}\Big(\frac{1+t}{1+\tau}\Big)^{C\tilde{\varepsilon}}
\sum_{|I|\leq 2}\!\!\!\|Z^I\! \phi(\tau,\cdot)\|_{L^\infty}\label{eq:decaywaveeq4}\\
&+ C\int_0^t\Big(\frac{1+t}{1+\tau}\Big)^{C\tilde{\varepsilon}}
\Big(\sum_{|I|\leq 1}\| r(\cdot) |Z^I
F|(\tau,\cdot)\|_{L^\infty(D_\tau)} +\sum_{|I|\leq 3}
(1+\tau)^{-1}\, \| Z^I \phi(\tau,\cdot)\|_{L^\infty(D_\tau)}\Big)\,
d\tau.\nn
\end{align}
\end{cor}
\begin{proof} By Lemma \ref{waveeqframe} for each component
we have the estimate \beq \Big|\Big(4\pa_s -
\frac{H_{LL}}{2g^{L\Lb}}\pa_q -\frac{\overline{\operatorname{tr}}\,
H \,+H_{L\underline{L}}}{2g^{L\underline{L}}\,\, r} \Big) \pa_q
(r\phi_{\mu\nu})+\frac{rF_{\mu\nu}}{2g^{L\underline{L}}}\Big| \les
\Big(1+ \frac{r\,|H|_{\cal LT}}{1+|q|}+ |H|\Big) r^{-1} \sum_{|I|\le
2}|Z^I \phi_{\mu\nu}| \eq and since $\pa_s$ and $\pa_q$ commute with
contractions with the frame vectors $L, \Lb$ we get \beq
\Big|\big(4\pa_s - \frac{H_{LL}}{2g^{L\Lb}}\pa_q
-\frac{\overline{\operatorname{tr}}\, H
\,+H_{L\underline{L}}}{2g^{L\underline{L}}\,\, r} \Big) \pa_q
(r\phi_{UV})+\frac{rF_{UV}}{2g^{L\underline{L}}}\Big| \les \Big(1+
\frac{r\,|H|_{\cal LT}}{1+|q|}+ |H|\Big) r^{-1} \sum_{|I|\le 2}|Z^I
\phi| \eq As before it also follows that \beq
(1+t+|r)|\pa\phi|_{UV}\les \sum_{|I|\leq 1}|Z^I\phi|
+|\pa_q(r\phi)|_{UV} \eq The lemma now follows as before.
\end{proof}

\section {Energy estimates for the wave equation on a curved space time  }
 \label{section:energywaveeq}
In this section we derive the energy estimate for a solution
$\phi$ of the inhomogeneous wave equation \beq\label{eq:waveeq}
\Boxr_g \phi = F \eq under the following assumptions on the metric
$g^{\a\b}=m^{\a\b} + H^{\a\b}$:
\begin{align}
&(1+|q|)^{-1} |H|_{LL} +|\pa H|_{LL}+|\overline{\pa} H|\leq
C\varepsilon (1+t)^{-1},\nn\\
& (1+|q|)^{-1}\,|H|+ |\pa H|\leq C\varepsilon (1+t)^{-\frac 12}
(1+|q|)^{-\frac 12-\gamma}\label{eq:metricdecay}
\end{align}
\begin{proposition}
\label{prop:Decayenergy}
Let $\phi$ be a solution of the wave equation \eqref{eq:waveeq}
with the metric $g$ verifying the assumptions \eqref{eq:metricdecay}.
Then for any $0<\ga\le 1/2$, there is an $\varepsilon_0$
such that for $\varepsilon<\varepsilon_0$
\beq\label{eq:firstenergy}
 \int_{\Si_{t}} |\pa\phi|^{2} + \int_{0}^{t} \int_{\Si_{\tau}}
\frac {\gamma \,\, |\pab\phi|^{2}}{(1+|q|)^{1+2\gamma}}\leq
 8\int_{\Si_{0}} |\pa \phi|^{2} +  C\varepsilon
\int_0^t\int_{\Si_{t}} \frac {|\pa\phi|^{2}}{1+t}
 +16\int_0^t  \int_{\Sigma_t} |F||\pa_t\phi|
\eq
\end{proposition}
\begin{remark}
Observe that by the Gronwall inequality the energy estimate of the
above proposition implies $t^\varepsilon$ growth of the energy.
For similar estimates, proved under different assumptions, see also
\cite{S1},\cite{A2},\cite{A3}.
\end{remark}
\begin{proof}
The proof of the proposition relies on the energy estimate obtained
in Corollary \ref{cor:Energy}
\begin{align}
 \int_{\Si_{t}} \big (|\pa_{t}\phi|^{2}
+|\nab \phi|^{2}\big ) & + \int_{0}^{t} \int_{\Si_{\tau}}
\frac {\gamma\,\, |\pab\phi|^{2}}{(1+|q|)^{1+2\gamma}}\le
 4 \int_{\Si_{0}} \big (|\pa_{t}\phi|^{2}
+|\nab \phi|^{2}\big ) \nn \\ &+ 8
\int_{0}^{t} \int_{\Si_{\tau}} \big |\pa_{\a} g^{\a\b}\pa_{\b}\phi
\pa_{t}\phi -\frac 12 \pa_{t} g^{\a\b}\pa_{\a}\phi\pa_{\b}\phi
+ F\pa_{t}\phi \big |
\nn\\ &+ 2\int_{0}^{t} \int_{\Si_{\tau}}
\frac{\gamma}{(1+|q|)^{1+2\ga}}\big |
(g^{\a\b}-m^{\a\b})\pa_{\a}\phi\pa_{\b}\phi +
2(g^{\Lb\b} - m^{\Lb\b})\pa_{\b}\phi\pa_{t}\phi \big |\nn
\end{align}
We start by decomposing the terms on the right hand side with
respect to the null frame.
$$
|\pa_{\a} g^{\a\b}\pa_{\b}\phi \pa_{t} \phi|\le
\big (|H|\, |\pa H|+ |(\pa H)_{LL}| + |\pab H|\big ) \, |\pa \phi|^2 +
|\pa H|\,|\pab \phi|\,|\pa\phi|
$$
Similarly,
$$
|\pa_{t} g^{\a\b}\pa_{\a}\phi \pa_{\b}\phi |\le
\big (|g-m|\, |\pa g|+ |(\pa g)_{LL}| + |\pab g|\big ) \, |\pa \phi|^2 +
|\pa g|\,|\pab \phi|\,|\pa\phi|
$$
Therefore, using the assumptions \eqref{eq:metricdecay} on the metric
$g$, we obtain that
\beq\label{eq:spaceq}
|\pa_{\a} g^{\a\b}\pa_{\b}\phi \pa_{t} \phi -\frac 12
\pa_{t} g^{\a\b}\pa_{\a} \phi \pa_{\b} \phi | \les
\frac {\varepsilon}{1+t} |\pa\phi|^2 + \frac {\varepsilon}{(1+|q|)^{1+2\ga}}
|\pab\phi|^2
\eq
Decomposing the remaining terms we infer that
$$
|(g^{\a\b}-m^{\a\b})\pa_{\a}\phi \pa_{\b} \phi |\le
 |H_{LL}| \, |\pa\phi|^2 +
|H| |\pab\phi|\,|\pa\phi|
$$
Similarly,
$$
|(g^{\a\b}-m^{\a\b}) L_\a \pa_{\b} \phi \pa_{t} \phi |\le
|H_{LL}| \, |\pa\phi|^2 + |H|\, |\pab\phi|\,|\pa\phi|
$$
Once again, using the assumptions \eqref{eq:metricdecay}, we have
\beq\label{eq:timeq}
|2 (g^{\a\b}-m^{\a\b}) L_\a \pa_{\b} \phi \pa_{t} \phi +
(g^{\a\b}-m^{\a\b})\pa_{\a}\phi \pa_{\b} \phi| \les
\varepsilon \frac {1+|q|}{1+t} |\pa\phi|^2 + \frac \varepsilon
{(1+|q|)^{2\ga}} |\pab\phi|^2
\eq
Thus
$$
\int_{\Si_{t}}|\pa\phi|^2  + \int_{0}^{t} \int_{\Si_{\tau}}
\frac {\gamma\,\, |\pab\phi|^{2}}{(1+|q|)^{1+2\gamma}}
\leq 4\int_{\Si_0}|\pa\phi|^2 +
C\varepsilon
\int_0^t \int_{\Si_\tau}\Big ( \frac {|\pa\phi|^2} {1+t}  +
\frac {|\pab\phi|^{2}}{(1+|q|)^{1+2\gamma}}\Big ) + 8\int_0^t \int_{\Si_\tau}
|F|\, |\pa_t \phi|
$$
and the desired estimate follows if we take $\varepsilon$ so small that
$C\varepsilon<\gamma/2$.
\end{proof}

\section{Estimates from the wave coordinate condition} \label{section:decaywavec}
In previous sections we have shown that one only needs to control
certain components of the metric in order to establish decay
estimates for solutions of the reduced wave equation. In this section we
will see that the wave coordinate condition allows one to estimate
precisely those components in terms of tangential derivatives or
higher order terms with better decay better. Recall that the  wave coordinate condition
can be written in the form
\beq\label{eq:wavecdef} \pa_\mu\Big( g^{\mu\nu}
\sqrt{|\det{g}|}\Big)=0
\eq
We have the following decomposition
$$
g^{\mu\nu} \sqrt{|\det{g}|}
=\big(m^{\mu\nu}+H^{\mu\nu}\big)\big(1-\frac{1}{2}\tr H+O(H^2)\big)
$$
where
$H^{\alpha\beta}=g^{\alpha\beta}-m^{\alpha\beta}$,
$h_{\alpha\beta}=g_{\alpha\beta}-m_{\alpha\beta}$. Recall also
that $g^{\alpha\beta}$ is the inverse of $g_{\alpha\beta}$ and
$H^{\a\b}=-m^{\mu\a} m^{\nu\b}h_{\mu\nu}+ O(h^2)$.
Therefore we obtain the following expression for the wave coordinate
condition:
\beq\label{eq:wavecappr}
\pa_\mu\Big( H^{\mu\nu}-\frac{1}{2}
m^{\mu\nu}\operatorname{tr}\, H+O^{\mu\nu}(H^2)\Big)=0
\eq
Using that we can express the divergence in terms of the null frame
\beq\label{eq:divnull}
\pa_\mu F^\mu=L_\mu\pa_q F^\mu -\underline{L}_\mu \pa_s F^\mu
+A_\mu \pa_A F^\mu
\eq
we obtain:
\begin{lemma} \label{decaywavec0} Assume that
$|H|\leq 1/4$. Then \beq \label{eq:decaywavec1} |\pa H|_{L \cal T}
\les |\overline{\pa } H| + |H|\, |\pa H| ,\qquad\quad |\pa\,
\overline{\operatorname{tr}}\, H|\les |\bar\pa H|+|H|\,|\pa H| \eq
\end{lemma}
\begin{proof} It follows from \eqref{eq:wavecappr} and
\eqref{eq:divnull} that
\beq
\big|L_\mu\pa\big(H^{\mu\nu}-\frac{1}{2}m^{\mu\nu}
\operatorname{tr}\, H\big)\big|\leq |\bar\pa H|+|H||\pa H|
\eq
Contracting with $T_\nu$ and using that $m_{TL}=0$ gives the
first inequality and contracting with $\underline{L}_\mu$
and using that $m_{\underline{L} L}=-2$ gives the second since
\beq
H_{L\underline{L}}+\operatorname{tr} \, H
=\overline{\operatorname{tr}}\, H
\eq
\end{proof}
We now compute the commutators of the wave coordinate condition
with the vector fields $Z$.
\begin{lemma}\label{decaywavectwo}
Let $Z$ be one of the Minkowski Killing or conformally Killing
vector fields and let tensor $H$ satisfy the wave coordinate
condition. Then the estimate
\begin{align}
\big| \pa H^J\big|_{L\cal T}\label{eq:ZIHLL}
\les\sum_{|J|\leq |I|} |\pab Z^J H|+ \sum_{I_{1}+...+I_k=I,\,
k\geq 2} |Z^{I_k} H|\cdot\cdot\cdot|Z^{I_2}  H| \, |\pa Z^{I_1} H|
\nn
\end{align}
holds true for the expression
\beq \label{defHJ}
H_{\mu\nu}^J=Z^J \widetilde H_{\mu\nu}+\!\!\!\!\sum_{|J|<|I|}
c_{J\, \mu}^{I\,\,\,\,\ga }  Z^{J} \widetilde H_{\gamma \nu}, \qquad
{\text {where}}\qquad \widetilde H_{\mu\nu}=H_{\mu\nu} - \frac 12 m_{\mu\nu}
\tr H
\eq
with some constant tensors  $c_{J}^{I\,\ga\mu}$ such
that $c_{JL}^{I\,\,\,\underline{L}}=0$ if $|J|=|I|-1$.
\end{lemma}
\begin{proof}
The wave coordinate condition \eqref{eq:wavecdef} can be written
in the form
$$
\pa_\mu \big( \widetilde{G}^{\mu\nu} \big)=0,\qquad \text{where}\quad
\widetilde{G}^{\mu\nu}=(m^{\mu\nu}+H^{\mu\nu})\sqrt{|\det g|}.
$$
Let $Z$ be one of the Minkowski Killing or conformally Killing
vector fields. Then for any vector field $F$ we have that
$$
Z^I \pa_\alpha F^\alpha =\pa_\alpha\Big(Z^I F^\alpha
+\sum_{|J|<|I|} c_{J\,\,\gamma}^{\,I\,\,\,\alpha} Z^{J}
F^{\gamma}\Big) =\pa_\alpha\Big( \sum_{|J\le |I|}
c_{J\,\gamma}^{\,I\,\,\,\alpha} Z^{J} F^{\gamma}\Big),
$$
where $c_{J\,\gamma}^{\,\,\,\alpha}$ are constants such that
$$
 c_{J\, \gamma}^{\,I\,\,\,\alpha}=\delta^{\alpha}_{\,\, \gamma},
\quad\text{if}\quad |J|=|I|\qquad \text{and}\qquad c_{J\,
L}^{\,I\,\,\,\Lb}=0,\quad\text{if}\quad |J|=|I|-1
$$
The last identity is a consequence of the relation between $c_{J\,
\a}^{\,I\,\,\,\ga}$ and the commutator constants
$c_{\a\b}=[\pa_{\a}, Z]_{\b}$ for which we have established that
$c_{LL}=0$. It therefore follows that
$$
\pa_\mu \big(\sum_{|J|\le |I|} c_{J}^{\,I\,\,\,\mu\gamma}
Z^{J}\widetilde{G}_{\gamma\nu}\big)=0.
$$
Decomposing relative to the null frame $(L,\Lb, S_1, S_2)$ we obtain
$$
\pa_q \big( \sum_{|J|\le |I|} c_{J}^{\,I\,\,\,\Lb\,\gamma }
Z^{J}\widetilde{G}_{\gamma \nu} \big) =\pa_s  \Big(\sum_{|J|\le |I|}
c_{J}^{\,I\,\,\, L \gamma}
 Z^{J}\widetilde{G}_{\gamma\nu} \Big)
- A_{\mu} \pab_A \big( \sum_{|J|\le |I|}
c_{J}^{\,I\,\,\,\mu\gamma } Z^{J}\widetilde{G}_{\,\ga \nu} \big).
$$
We now contract the above identity with one of the tangential
vector fields $T^{\nu}$, $T\in \{ L,S_1,S_2\}$ to obtain
$$
\Big|L^\gamma \, T^\nu \pa_q  Z^I \widetilde{G}_{\gamma\nu}
 +\!\!\!\!\sum_{|J|<|I|}
c_{J}^{\,I\,\,\,\Lb\ga}  T^\nu\, \pa_q Z^{J}\widetilde{G}_{\gamma \nu}
\Big| \les \sum_{|J|\leq |I|} \big| \pab Z^{I}\widetilde{G}\big|
$$
We examine the expression
$$
L^\gamma \, T^\nu Z^J  \pa_q \widetilde{G}_{\gamma\nu}= L^\gamma \, T^\nu
\pa_q Z^J \bigg (
(m_{\gamma\nu}+H_{\gamma\nu})\sqrt{|\det g|} \bigg )=\sum_{J_1+J_2=J}
L^\gamma \, T^\nu \pa_q \bigg ( (Z^{J_1} H_{\gamma\nu}) Z^{J_2}
\sqrt{|\det g|}\bigg )
$$
since $m_{LT}= L^\gamma \, T^\nu m_{\gamma\nu}=0$. The desired
estimate now follows from the identity $\sqrt{|\det g|}=1+f(H)$,
which holds with a smooth function $f(H)$  such that
$f(H)=-\tr H/2+ O(H^2)$
\end{proof}
We now summarize the results of this section in the following
\begin{lemma}\label{decaywavecZ}
For a tensor $H$ obeying the wave coordinate condition \beq
\label{eq:decaywavec3} |\pa H|_{ LT}\les |\bar\pa H|+|H|\,|\pa H|,
\eq and
 \beq \label{eq:decaywavec3} |\pa Z H|_{ LL}\les |\pa
H|_{L\cal T} + \sum_{|I|\le 1} |\overline{\pa} Z^I H| +
\sum_{|I|+|J|\le 1} |Z^I H|\,|\pa Z^{J} H|. \eq
 In general,
 \begin{multline} \label{eq:decaywavec6}
|\pa Z^I H|_{ L\cal T}\les \sum_{|J|\leq |I|}|\overline{\pa} Z^J
H|+\!\!\!\!\sum_{|J|\leq |I|-1}\!\!\!|\pa Z^J H|\,\, +
\!\!\!\!\!\!\!\!\!\!\!\! \sum_{\,\,\,\,|I_1|+...+|I_m| \leq |I|,\,
m\geq 2}\!\!\!\!\!\!\!\!\!\!\!\!\!\! |Z^{I_m}
H|\!\cdot\cdot\cdot|Z^{I_{2}} H| |\pa Z^{I_1} H|,
\end{multline}
and
\begin{multline} \label{eq:decaywavec5}
|\pa Z^I H|_{ LL}\les  \sum_{|J|\leq |I|}|\overline{\pa} Z^J H|+
\sum_{|J|\leq |I|-1}|\pa Z^J H|_{L \cal T}
+\sum_{|K|\leq |I|-2} |\pa Z^{J} H| \\
+ \!\!\!\!\!\!\!\!\!\!\!\! \sum_{|I_1|+...+|I_m| \leq |I|,\, m\geq
2}\!\!\!\!\!\!\!\!\!\!\!\!\!\! |Z^{I_m}
H|\!\cdot\cdot\cdot|Z^{I_{2}} H| |\pa Z^{I_1} H|.
\end{multline}
 The same estimates also hold for $H$ replaced by $h$.
\end{lemma}
\begin{proof} This follows directly by the previous lemma with the help of
the identities $m_{LT}=0$ and $c_{J\,L}^{I\,\,\,\Lb}=0$.
\end{proof}

\section{Estimates for the inhomogeneous terms}\label{section:decayinhom}
In this section we will show that the inhomogeneous terms of the
reduced Einstein equations can be estimated in terms of tangential
derivatives, for which we have better decay estimates, or tangential
components which in turn can be expressed, using the wave coordinate
condition, in terms of tangential derivatives and lower order terms.
Recall that according to Lemma \ref{Einstwavecquad} the symmetric
two tensor $h_{\mu\nu}= g_{\mu\nu}- m_{\mu\nu}$ verifies the reduced
Einstein equations of the form:
\begin{align}
&\Boxr_g h_{\mu\nu} = F_{\mu\nu} (h) (\pa h, \pa h),\nn\\
& F_{\mu\nu}(h)(\pa h,\pa h)=P(\pa_\mu h,\pa_\nu h)+Q_{\mu\nu}(\pa
h,\pa h)
+G_{\mu\nu}(h)(\pa h,\pa h)\label{eq:F},\\
& P(\pa_\mu k,\pa_\nu p)=\frac 14 \pa_{\mu} \tr k\, \pa_{\nu}\tr p
-\frac 12 \pa_{\mu} k^{\a\b}\pa_{\nu} p_{\a\b},\label{eq:P}
\end{align}
Here $Q{\mu\nu}$ are linear combinations of the null-forms and
$G_{\mu\nu}(h)(\pa h,\pa h)$ is a quadratic form in $\pa h$ with
coefficients that are smooth functions of $h$ and vanishing at
$h=0$.

\begin{lemma}
\label{le:Pest}
The quadratic form $P$ satisfies the following pointwise estimate:
\begin{align}
&|P(\pa p,\pa k)|_{\cal TU}\les |\pab p\,|\, |\pa k| + |\pa p\,|\,
|\pab k|,
\label{eq:PTU}\\
&|P(\pa p,\pa k)|\les |\pa p\,|_{\cal TU} |\pa k|_{\cal TU} + |\pa
p\,|_{LL} |\pa k| + |\pa p\,| \, |\pa k|_{LL}\label{eq:Pall}
\end{align}
\end{lemma}
\begin{proof}
The first part of the statement follows trivially from \eqref{eq:P}.
To prove \eqref{eq:Pall} we use \eqref{eq:tanP} applied to
$R^\mu\pa_\mu p$ in place of $p$ and $S^\nu\pa_\nu k$ in place of
$k$, for any vector fields $T$ and $S$, to obtain \beq |T^\mu
S^\nu\, P(\pa_\mu p,\pa_\nu k)|\les |T^\mu\pa_\mu p\,|_{\cal TU}
|S^\nu\pa_\nu k|_{\cal TU}+|T^\mu\pa_\mu p\,|_{LL}|S^\nu\pa_\nu k|
+|T^\mu\pa_\mu p\,|\, |S^\nu\pa_\nu k|_{LL} \eq which proves the
lemma.
\end{proof}
Using the additional estimates
on the $h_{LL}$ component, derived in Lemma \ref{decaywavecZ} under the
assumption that the wave coordinate condition holds, we obtain the
following:
\begin{cor}
Under the additional assumption that $h$ satisfies the wave
coordinate condition \eqref{eq:wavecord}, the quadratic form $P$
obeys the estimate
\begin{align}
&|P(\pa h,\pa h)|_{\cal TU}\les |\pab h|\, |\pa h|,\label{eq:Ph}\\
&|P(\pa h,\pa h)|\les |\pa h|^2_{\cal TU} + |\pab h|\,
|\pa h| + |h|\,|\pa h|^2\label{eq:Pallcoord}
\end{align}
Moreover,
\begin{align*}
|Z^I P(\pa h, \pa h)|&\les \sum_{|J|+|K|\le |I|}
\big (|\pa Z^J h|_{\cal TU}\, |\pa Z^K h |_{\cal TU} +
|\pab Z^J h|\, |\pa Z^K h| \big )+ \sum_{|J|+|K|\le |I|-1}
|\pa Z^J h|_{L \cal T} |\pa Z^K h|\\ & +
\sum_{|J|+|K|\le |I|-2} |\pa Z^J h|\, |\pa Z^K h| +
\sum_{|J_1|+...+|J_m| \leq |I|,\, m\geq 3}
|Z^{J_m} h|\!\cdot\cdot\cdot |Z^{J_3} h|\,|\pa Z^{J_{2}} h| |\pa Z^{J_1} h| \label{eq:PZI}
\end{align*}
\end{cor}
\begin{proof} The inequality \eqref{eq:Ph} follows directly
from \eqref{eq:PTU}. To prove \eqref{eq:Pallcoord}
we use \eqref{eq:Pall} and
that by the wave coordinate condition
$|\pa h|_{LL}\les |\bar \pa h|+|h|\, |\pa h|$.

We now note that $Z^I P(\pa_\mu h,\pa_\nu h)$ is a sum of terms of the
form $P(\pa_\alpha Z^J h,\pa_\beta Z^K h)$ for some $\alpha$, $\beta$ and
$|J|+|K|\leq I$:
$$
|Z^I P(\pa h,\pa h)|\leq C\sum_{|J|+|K|\leq |I|}
|P(\pa  Z^J h,\pa Z^K h)|
$$
It follows from Lemma \ref{tander} and Lemma \ref{decaywavecZ} that
\begin{multline}
\!\!\!\!\!\sum_{|J|+|K|\leq |I|}\!\!\!\!\!
|P(\pa Z^J h,\pa Z^K h)|
\les \!\!\!\!\!\sum_{|J|+|K|\leq |I|}\!\!\!\!\!
|\pa Z^J h|_{\cal TU}|\pa Z^K h|_{\cal TU}
+|\pa Z^J h|_{LL} \,|\pa Z^K h |\\
\les \sum_{|J|+|K|\leq |I|} |\bar\pa Z^J h|\, |Z^K h|
+|\pa Z^J h|_{\cal TU}\, |\pa Z^K h|_{\cal TU} \\
+\!\!\!\!\! \sum_{|J|+|K|\leq |I|}\!\!\!\!\Big(
\sum_{|J^\prime|\leq |J|-1} \!\!\!\!\!\!
|\pa Z^{J^\prime} h|_{\cal LT}\,\,
+\!\!\!\!\!\!\sum_{|J^{\prime\prime}|\leq |J|-2} \!\!\!\!\!|\pa Z^{J^{\prime\prime}} h|\,\, +\!\!\!\!\!\!\!\!\!\!\!\!\!\!
\sum_{|J_1|+...+|J_m|\leq |J|,\, m\geq 2}
\!\!\!\!\!\!\!\!\!\!\!\!|Z^{J_m} h|\,\cdot\cdot\cdot
|Z^{J_2} h|\, |\pa Z^{J_1} h|\Big) |\pa Z^K h|
\end{multline}
which proves the lemma.
\end{proof}

\begin{prop} \label{decayinhom}
Let $F_{\mu\nu}=F_{\mu\nu}(h)(\pa h,\pa h)$ be  as in Lemma  \ref{Einstwavecquad}
and assume that the wave coordinate condition holds.
Then
\beq \label{eq:decayinhom1}
| F|_{\cal T U}\les |\overline{\pa} h|\, |\pa h|+|h|\,|\pa h|^2
\eq
and
\beq \label{eq:decayinhom2}
|F|\les |\pa h|^2_{\cal TU} + |\pab h|\,
|\pa h| + |h|\,|\pa h|^2
\eq
\beq  \label{eq:decayinhom3}
|Z F|\les \big(|\pa h|_{\cal TU}+|\overline{\pa} h|+|h|\,|\pa h|\big)
(|\pa Z h|+|\pa h|)
+\big(|\overline{\pa} Z h|+|Z h|\, |\pa h|\big)|\pa h|
\eq
\begin{multline} \label{eq:decayinhom4}
|Z^I F| \les
 \sum_{|J|+|K|\le |I|}
\big (|\pa Z^J h|_{\cal TU}\, |\pa Z^K h |_{\cal TU} +
|\pab Z^J h|\, |\pa Z^K h| \big )+ \sum_{|J|+|K|\le |I|-1}
|\pa Z^J h|_{L \cal T} |\pa Z^K h|\\  +
\sum_{|J|+|K|\le |I|-2} |\pa Z^J h|\, |\pa Z^K h| +
\sum_{|J_1|+...+|J_m| \leq |I|,\, m\geq 3}
|Z^{J_m} h|\!\cdot\cdot\cdot |Z^{J_3} h|\,|\pa Z^{J_{2}} h| |\pa Z^{J_1} h|
\end{multline}
\end{prop}
\begin{proof} First
$$
|Z^I G_{\mu\nu}(h)(\pa h,\pa h)|
\leq C\sum_{ |I_1|+...+|I_k|\leq |I|,\, k\geq 3}
 |Z^{I_k} h|\cdot\cdot\cdot |Z^{I_3} h|\,|\pa Z^{I_2} h|
\,|\pa Z^{I_1} h|.
$$
Since $Z Q(\pa u,\pa v)=Q(\pa u,\pa Zv)+Q(\pa Zu,\pa v)
+a^{ij}Q_{ij}(\pa u,\pa v)$, and
$|Q_{\mu\nu}( \pa h,\pa k)|\leq |\pa h|\, |\overline{\pa} k|+
 |\pa k|\, |\overline{\pa} h|$ it follows that
$$
|Z^I Q_{\mu\nu}(\pa h,\pa h)|\leq C\sum_{|J|+|k|\leq |I|}
|Q_{\mu\nu}(\pa  Z^J h,\pa Z^K h)|\leq  C\sum_{|J|+|k|\leq |I|}
|\pa Z^J h|\, |\overline{\pa} Z^K h|
$$
\end{proof}

\section
{The decay estimates for Einstein's
equations}\label{section:decayeinst} In this section we will
establish the improved decay estimates for Einstein's equations. Our strategy is to
use the weak decay estimates, obtained from the assumed energy
bounds, to prove sharper decay estimates and then to recover
the energy bounds in the next section.
\begin{theorem}\label{decayeinst}
 Suppose that for some $0< \gamma\leq 1/2$
\begin{align}
&| \pa Z^I h|\leq C\varepsilon (1+t+|q|)^{-1/2-\gamma}
(1+|q|)^{-1/2-\gamma},
\qquad  |I|\leq N/2+4, \label{eq:weakdecay1}\\
&|Z^I h|\leq C\varepsilon (1+t)^{-1},\qquad q= 1,\qquad |I|\leq
N/2+4 \label{eq:boundcond}
\end{align}
hold for $0\leq t\leq T$.
Then for $0\leq t\leq T$ we have
\begin{align}
&|Z^I h|\leq  C\varepsilon (1+t+|q|)^{-1/2-\gamma}
(1+|q|)^{1/2-\gamma},\qquad  |I|\leq N/2+4 ,\label{eq:weakdecay2}\\
&|\overline{\pa} Z^I h|\leq C\varepsilon (1+t+|q|)^{-3/2-\gamma}
(1+|q|)^{1/2-\gamma}, \qquad  |I|\leq N/2+3 .\label{eq:tandecay}
\end{align}
Assume also that
 $h$ satisfies the wave coordinate condition. Then
for $0\leq t\leq T$ we have \beq \label{eq:sharpdecay1} |\pa h|_{L
\cal T} +|\pa Z h|_{LL}\leq C\varepsilon (1+t)^{-1-2\gamma},
\qquad\text{and}\qquad |h|_{L \cal T} +|Z h|_{LL}\leq C\varepsilon
(1+t) ^{-1}(1+|q|) . \eq
Furthermore if in addition $h$ satisfies Einstein's equations then
for $\varepsilon$ sufficiently small and $0\leq t\leq T$
we also have
\begin{align}
&|\pa h|_{\cal T U}\leq  C\varepsilon (1+t)^{-1},
\qquad |h|_{\cal T U}\leq  C\varepsilon
(1+t)^{-1}(1+|q|), \label{eq:sharpdecay2}\\
&|\pa h|\leq  C\varepsilon (1+t)^{-1} \ln (2+t). \label{eq:sharpdecay3}
\end{align}
In general, there are constants $M_k$, $C_k$ and
$\varepsilon_k>0$ such that if $\varepsilon\leq \varepsilon_k$
then for $|I|=k\leq N/2+2$
\begin{equation}\label{eq:sharpdecay4}
| \pa Z^I h|\leq C_k\varepsilon (1+t)^{-1+M_k\varepsilon}  ,
\qquad \text{and}\qquad |Z^I h|\leq C_k\varepsilon
(1+t)^{-1+M_k\varepsilon}(1+|q|)
\end{equation}
\end{theorem}
\begin{remark}
We remind the reader that, as stated in the Remark \ref{rem:distinct},
our estimates make no distinction between the tensors $h$ and
$H=-h + O(h^2)$. In particular, one can directly verify
that the conclusions of the theorem also hold for the tensor $H$.
\end{remark}

First we note that all the estimates
\eqref{eq:weakdecay2}-\eqref{eq:sharpdecay4} trivially follow
from the assumptions \eqref{eq:weakdecay1}-\eqref{eq:boundcond}
away from the light cone, thus the theorem is only useful in the
region $t/2\leq |x|\leq 2t$. The estimate
\eqref{eq:weakdecay2} follows from integrating
\eqref{eq:weakdecay1} from $q=1$, where \eqref{eq:boundcond} hold.
Similarly the second parts of \eqref{eq:sharpdecay1},
\eqref{eq:sharpdecay2} and \eqref{eq:sharpdecay4} follow from
integrating the first and using \eqref{eq:boundcond}. It follows
from \eqref{eq:weakdecay2} and Lemma \ref{tanderZ} that we have
the better estimate \eqref{eq:tandecay} for the derivatives tangential
to the outgoing Minkowski cones.
The inequalities \eqref{eq:sharpdecay1}-\eqref{eq:sharpdecay4} for tangential
derivatives certainly follow from \eqref{eq:tandecay}, so it only
remains to prove these estimates for a derivative transversal to
the light cone.

\noindent
The missing improved estimates for a $(\pa_t-\pa_r)$ derivative transversal to the
light cones will be obtained, in the case of \eqref{eq:sharpdecay1},
from the wave coordinate condition, see section \ref{section:decaywavec}, and for
\eqref{eq:sharpdecay2}-\eqref{eq:sharpdecay4},
from integrating the reduced Einstein wave equations,  see section
\ref{section:decaywaveeq}.
The estimates from the
wave coordinate condition are easily obtained. In fact the first
estimate in \eqref{eq:sharpdecay1} follows directly from Lemma
\ref{decaywavec0} using the estimates \eqref{eq:weakdecay1},
\eqref{eq:weakdecay2} and \eqref{eq:tandecay} and the second
estimate in \eqref{eq:sharpdecay1} follows integrating the first
from $q=1$ where \eqref{eq:boundcond} holds. However, the wave
coordinate condition does not give estimates for a transversal
derivative of all components of the metric and the remaining
components have to be controlled by integrating the wave equation
expressed in polar coordinates. The estimates for the transversal
derivative obtained from the wave coordinate condition rely on a decomposition
of the metric with respect to the  null
frame.
On the other hand, the estimates obtained from
integrating the wave equation are based on a decomposition of the wave
operator in terms of tangential derivatives and a transversal
derivative.

\subsection{Proof of \eqref{eq:weakdecay2} and \eqref{eq:tandecay}}
For a fixed angular variable $\omega$
we integrate in the radial direction and use \eqref{eq:weakdecay1}
and \eqref{eq:boundcond}
\begin{multline}
 |Z^I
h(t,r,\omega)|\leq |Z^I h(t,t+1,\omega)| +\int_r^{t+1}|\pa_r Z^I
h(t,\rho,\omega)|\, d\rho\\
\les \frac{C\ve}{1+t}+\int_r^{t+1} \frac{C\ve\, d\rho
}{(1+t+|t-\rho|)^{1/2+\gamma}(1+|t-\rho|)^{1/2+\gamma}}\les
\frac{C\ve}{1+t}+\frac{C\ve(1+|t-r|)^{1/2+\gamma}}{(1+t+r)^{1/2+\gamma}}
\end{multline}
The estimate \eqref{eq:weakdecay2} now follows. By Lemma
\ref{tanderZ} and \eqref{eq:weakdecay2}
$$
|\overline{\pa} Z^I h|\les \frac{1}{1+t+|q|} \sum_{|J|\leq |I|+1}
|Z^J h|\les\frac{\ve (1+|q|)^{1/2-\gamma}}{(1+t+|q|)^{3/2+\gamma}}
$$
which proves \eqref{eq:tandecay}.

\subsection{Proof of \eqref{eq:sharpdecay1}.}
We now show that the wave coordinate condition allows one
to control certain components by  lower order terms and
terms with fast decay.
\begin{lemma} \label{wavecdecay4} Suppose that the estimates \eqref{eq:weakdecay1}-
\eqref{eq:tandecay} hold and that $h$ satisfies the wave
coordinate condition. Then
\beq \label{eq:absent}
\sum_{|I|\leq \, k}|\pa Z^I
h|_{LL} + \sum_{|J|\leq \, k-1} |\pa Z^{J} h|_{\cal LT}\les
\sum_{|K|\leq \, k-2} |\pa Z^K
h|+\varepsilon(1+t+|q|)^{-1-2\gamma}
\eq
Here the sum over $k-2$
is absent if $k\leq 1$ and the sum over $k-1$ is absent if $k=0$.
Furthermore
\beq \label{eq:absent1}
\frac{1}{1+|q|}\Big(\sum_{|I|\leq\,  k}|Z^I
h|_{LL} + \sum_{|J|\leq \, k-1} |Z^{J} h|_{\cal LT}+ \sum_{|K|\leq
\, k-2} |Z^K h|\Big)(t,x) \les \sup_{ t/2
\leq |y|\leq 3t/2}
\sum_{|K|\leq \, k-2}  |\pa Z^K h(t,y)| +\frac{\varepsilon}{1+t}
\eq
\end{lemma}
\begin{proof}
We first prove \eqref{eq:absent}.
Using the estimates of Lemma \ref{decaywavectwo}
derived from the wave coordinate condition followed by  \eqref{eq:weakdecay1}-
\eqref{eq:tandecay} we obtain
 \begin{multline}
\sum_{|I|\leq k}|\pa Z^I H|_{LL}+\!\!\sum_{|J|\leq k-1}\!\!|\pa
Z^{J} H|_{L\cal T}\les \sum_{|I|\leq k} |\bar\pa Z^I
h|\,+\!\!\!\sum_{|K|\leq \, k-2}\!\!\!|\pa Z^K
h|+\!\!\!\!\!\!\sum_{|I_1|+...+|I_m|\leq k,\, m\geq 2}
\!\!\!\!\!\! |Z^{I_m} h|\cdot\cdot\cdot |Z^{I_2} h|\, |\pa Z^{I_1}
h|\\
\leq \sum_{|K|\leq \, k-2}\!\!\!|\pa Z^K h| +\varepsilon
(1+t+|q|)^{-1-2\gamma}+\varepsilon (1+t+|q|)^{-1-2\gamma}
 \end{multline}
The proof of estimate \eqref{eq:absent1} for  $|q|\geq t/2$
follows directly from \eqref{eq:weakdecay2}.  Thus we may assume
that $|q|\leq t/2$. We now use the inequality
\beq\label{eq:derest} |H(t,r\omega)|\leq
|H\big(t,(t+1)\omega\big)|+(1+|q|)\sup_{|\rho|\leq |q|+1}
|\pa_\rho H(t,(t+\rho)\omega)|, \eq and the boundary condition
\eqref{eq:boundcond} to conclude that
 \beq \frac{|Z^I
H|_{LL}+|Z^{J} H|_{L\cal T} +|Z^{K} H|}{1+|q|}\les \sup_{t/2\leq
|y|\leq 2t} \Big(|\pa_r Z^I H|_{LL}+|\pa_r Z^{J} H|_{L\cal T}
+|\pa_r Z^{K} H| \Big)(t,y)+\frac{\varepsilon}{1+t} \eq
The desired result now follows from \eqref{eq:absent}.
\end{proof}

The first part of \eqref{eq:sharpdecay1} now follows directly from
the Lemma with $k=0,1$ and the second part follows from
integrating the first and using the boundary assumption
\eqref{eq:boundcond} as in the proof of \eqref{eq:weakdecay2}.

\subsection{Proof of \eqref{eq:sharpdecay2}-\eqref{eq:sharpdecay3}.}
We will appeal to  the $L^\infty$
estimates of section \ref{section:decaywaveeq} for the reduced
wave equation
$$
{\Boxr}_g h_{\mu\nu}=F_{\mu\nu},
$$
where $F_{\mu\nu}$ is as in Lemma \ref{Einstwavecquad}.
We will now prove \eqref{eq:sharpdecay2} and \eqref{eq:sharpdecay3} assuming
\eqref{eq:weakdecay1}-\eqref{eq:sharpdecay1}.
\begin{lemma} \label{inhomdecay}
Suppose that the assumptions of Proposition \ref{decayeinst} hold
and let $F_{\mu\nu}=F_{\mu\nu}(h)(\pa h,\pa h)$ be  as in Lemma
\ref{Einstwavecquad}.
Then
\beq \label{eq:decayinhom1}
| F|_{\cal T U}\leq
 C\varepsilon t^{-1-2\gamma}|\pa h|
\eq
and
\beq \label{eq:decayinh}
|F|\leq C\varepsilon t^{-1-2\gamma}|\pa h| +C|\pa h|_{\cal TU}^2
\eq
\end{lemma}
\begin{proof} This follows from Lemma \ref{decayinhom} using
\eqref{eq:weakdecay1}-\eqref{eq:sharpdecay1}.
\end{proof}

Using the first part of Corollary \ref{decaywaveeq3}; \eqref{eq:decaywaveeq3},
and \eqref{eq:weakdecay1}-\eqref{eq:sharpdecay1} and
the previous lemma we get
\begin{lemma} With a constant depending on $\gamma>0$ we have
\begin{equation}
(1+t)\|\, |\pa h|_{\cal T U} (t,\cdot)\|_{L^\infty}
\leq  C\varepsilon
+ C\varepsilon
\int_0^t(1+\tau)^{-2\gamma}\|\pa h(\tau,\cdot)\|_{L^\infty}\, d\tau,
\end{equation}
and
\begin{equation}
(1+t)\|\pa h(t,\cdot)\|_{L^\infty}
\leq  C\varepsilon
+ C\int_0^t
\Big(\varepsilon (1+\tau)^{-2\gamma}\|\pa h(\tau,\cdot)\|_{L^\infty}+
(1+\tau)\||\pa h|_{\cal T U}(\tau,\cdot)\|_{L^\infty}^2\Big)\, d\tau.
\end{equation}
\end{lemma}

The estimates \eqref{eq:sharpdecay2} and \eqref{eq:sharpdecay3} now follow from
the above lemma and the following technical result applied to
$n_{00}(t)=(1+t)\|\, |\pa h|_{\cal T U} (t,\cdot)\|_{L^\infty}$ and
$n_{01}(t)=(1+t)\|\pa h(t,\cdot)\|_{L^\infty}$:
\begin{lemma} \label{ode1} Suppose that $n_{00}\geq 0$ and $n_{01}\geq 0$ satisfy
\begin{align}\label{eq:ode1}
n_{00}(t)&\leq C\varepsilon \Big(\int_0^t (1+s)^{-1-\gamma} n_{01}(s)\, ds+1\Big) \\
n_{01}(t)&\leq  C\varepsilon \Big(\int_0^t (1+s)^{-1-\gamma} n_{01}(s)\, ds+1\Big)
+C\int_0^t  (1+s)^{-1} n_{00}(s)^2\, ds
\end{align}
for some positive constants such that $0<16(C^2+C)\varepsilon<\gamma\leq 1$.
Then
\beq\label{eq:odebound}
n_{00}(t)\leq 2C\varepsilon,  \qquad\text{and}\qquad
n_{01}(t)\leq  2C\varepsilon \big(1+\gamma \ln{(1+t)}\big)
\eq
\end{lemma}
\begin{proof} Let $T$ be the largest time such that
\beq\label{eq:N1bound}
N_{01}(t)=\int_0^t (1+s)^{-1-\gamma} n_{01}(s)\, ds+1 \leq 2
,\qquad\text{for}\qquad t\leq T
\eq
Then for $t\leq T$ \eqref{eq:odebound} holds and since,
$$
\int_0^\infty (1+s)^{-1-\gamma}   \big(1+\gamma\ln{(1+s)}\big)\, ds=
\gamma^{-1}\int_0^\infty (1+\tau) \, e^{-\tau}\, d\tau =2\gamma^{-1}+1
$$
it follows that
$$
N_{01}(t)\leq 2C\varepsilon\big(2\gamma^{-1}+1\big)+1\leq 3/2
,\qquad\text{for}\qquad t\leq T .
$$
Since $N_{01}(t)$ is continuous this contradicts that $T$ is the
maximal number such that \eqref{eq:N1bound} holds. Thus
$T=\infty$ and \eqref{eq:odebound} holds for all $t<\infty$.
\end{proof}

This proves the first part of \eqref{eq:sharpdecay2} and
\eqref{eq:sharpdecay3}.  The second part of
\eqref{eq:sharpdecay2} follows from integrating the first using
the boundary assumption \eqref{eq:boundcond} as in the
proof of \eqref{eq:weakdecay2}.

\subsection{Proof of \eqref{eq:sharpdecay4} in case $k=1$.}
We will now prove the first part of \eqref{eq:sharpdecay4} for
$|I|=1$ assuming \eqref{eq:weakdecay1}-\eqref{eq:sharpdecay3}.
\begin{lemma} \label{inhomdecay}
Suppose that the assumptions of Proposition \ref{decayeinst} hold
and let $F_{\mu\nu}=F_{\mu\nu}(h)(\pa h,\pa h)$ be  as in Lemma
\ref{Einstwavecquad}. Then \beq  \label{eq:decayinhom3b} |Z F|\leq
C\varepsilon t^{-1}\big(|\pa Z h| +|\pa h| \big) \eq
\end{lemma}
\begin{proof} This follows from Lemma \ref{decayinhom}.
\end{proof}

Using the second part of Corollary \ref{decaywaveeq3};
\eqref{eq:decaywaveeq4},
and \eqref{eq:weakdecay1}-\eqref{eq:sharpdecay1} and
the previous lemma we get
\begin{lemma} If $\varepsilon>0$ is sufficiently small then
\begin{multline}
(1+t)\sum_{|I|\leq 1}\|\pa Z^I h(t,\cdot)\|_{L^\infty}
\leq  C\varepsilon (1+t)^{C\varepsilon}\Big(1+\int_0^t (1+\tau)^{-C\varepsilon}
\sum_{|I|\leq 1}\|\pa Z^I h(\tau,\cdot)\|_{L^\infty}\, d\tau\Big).
\end{multline}
\end{lemma}
The estimate
\eqref{eq:sharpdecay4} for $|I|=1$ is now a consequence of the above lemma
and the following technical result applied to
$n_1(t)=(1+t)^{1-C\varepsilon}\sum_{|I|\leq 1}\|\pa Z^I h(t,\cdot)\|_{L^\infty}$:

\begin{lemma} Suppose that $n_1(t)\geq 0$ satisfies
\beq
n_1(t)\leq C\varepsilon \Big(1+\int_0^t (1+\tau)^{-1} n_1(\tau)\, d\tau\Big)
\eq
Then
\beq
n_1(t)\leq C \varepsilon (1+t)^{C\varepsilon}
\eq
\end{lemma}
\begin{proof}
\beq
N_1(t)=1+\int_0^t (1+\tau)^{-1} n_1(\tau)\, d\tau
\eq
satisfies $\dot{N}_1(t)\leq C\varepsilon (1+\tau)^{-1}N_1(t)$.
Multiplying by the integrating factor $(1+t)^{-C\varepsilon} $ and integrating
we get $N_1(t)\leq N_1(0)(1+t)^{C\varepsilon}=(1+t)^{C\varepsilon}$ and the
lemma follows.
\end{proof}

\subsection{Proof of \eqref{eq:sharpdecay4} in case $k\geq 1$.}
We will now use induction to prove the first part of
\eqref{eq:sharpdecay4} for $|I|= k+1$ assuming that
\eqref{eq:weakdecay1}-\eqref{eq:sharpdecay1}, the first part of
\eqref{eq:sharpdecay2}, \eqref{eq:sharpdecay3} and the first part
of \eqref{eq:sharpdecay4} for $|I|\leq k$ hold.

\begin{lemma} \label{inhomdecay}
Suppose that the assumptions of Proposition \ref{decayeinst} hold
and let $F_{\mu\nu}=F_{\mu\nu}(h)(\pa h,\pa h)$ be  as in Lemma
\ref{Einstwavecquad}. Then \beq
 |Z^I F|\leq C\varepsilon t^{-1} \sum_{|K|\leq |I|}|\pa Z^K h|
+ C\sum_{|J|+|K|\leq |I|,\, |J|\leq |K|<|I|}|\pa Z^{J} h||\pa Z^{K} h|
\eq
\end{lemma}
\begin{proof} This follows from Lemma \ref{decayinhom} using
\eqref{eq:weakdecay1}-\eqref{eq:sharpdecay3}.
\end{proof}
By Corollary \ref{commut3}
\begin{multline}
|{\Boxr}_g Z^I h|\les |\hat{Z}^I F| +(1+t)^{-1}\sum_{|K|\leq
|I|,}\,\,\sum_{|J|+(|K|-1)_+\leq |I|} \!\!\!\!\!\!\!\!\!\!
|Z^{J} H| |\pa Z^{K} h|\\
+ C(1+q)^{-1}
 \sum_{|K|\leq |I|}\Big(\sum_{|J|+(|K|-1)_+\leq |I|} \!\!\!\!\!|Z^{J} H|_{LL}
+\!\!\!\!\!\sum_{|J^{\prime}|+(|K|-1)_+\leq |I|-1}\!\!\!\!\!|Z^{J^{\prime}} H|_{L\cal T}
+\!\!\!\!\!\sum_{|J^{\prime\prime}|+(|K|-1)_+\leq |I|-2}\!\!\!\!\!
|Z^{J^{\prime\prime}} H|\Big)|\pa Z^{K} h|
\end{multline}
where $(|K|-1)_+=|K|-1$, if $|K|\geq 1$, and $0$, if $|K|=0$.
Using Lemma \ref{wavecdecay4} we get \beq
(1+q)^{-1}\!\!\!\!\!\!\!\!\!\! \sum_{|J|\leq k,\,|J^\prime|\leq
k-1,\,|J^{\prime\prime}|\leq k-2}
 |Z^J H|_{LL}+|Z^{J^\prime} H|_{L\cal T}
+|Z^{J^{\prime\prime}} H|\leq \frac{C\varepsilon}{1+t}+
\sum_{|J^{\prime\prime}|\leq k-2} \, \sup_{t/2\leq |y|\leq 2t} \,
|\pa Z^{J^{\prime\prime}} H(t,y)| \eq
 We hence obtain \beq |{\Boxr}_g Z^I
h|\leq C\varepsilon (1+t)^{-1} \sum_{|K|\leq |I|}|\pa Z^K h|
+\,\,\,\sum_{|J|+|K|\leq |I|-1} \, \sup_{t/2\leq |y|\leq 2t} \, |
\pa Z^{J} H(t,y)|\, |\pa Z^K h| \eq Then we have proven that
\begin{lemma} Let
\beq n_k(t)=(1+t)\sum_{|I|\leq k}\| \pa Z^I
h(t,\cdot)\|_{L^\infty}. \eq Then for $|I|=k$: \beq |\Boxr_g Z^I
h|\leq C (1+t)^{-2}\big( \varepsilon n_k(t)+  n_{k-1}(t)^2\big)
\eq
\end{lemma}

By the first part of Corollary \ref{decaywaveeq3};
\eqref{eq:decaywaveeq3}, it therefore follows that:
\begin{lemma}
\beq
n_k(t)\leq C\varepsilon+ C\int_0^t (1+\tau)^{-1}
\big( \varepsilon n_k(\tau)+  n_{k-1}(\tau)^2\big)
\, d\tau
\eq
\end{lemma}
Our inductive hypothesis is
$n_{k-1}(t)^2\!\leq\! C\varepsilon^2 (1+t)^{C\varepsilon}$ so
the bound $n_k(t)\!\leq \!C\varepsilon (1+t)^{2C\varepsilon}$ follows from:
\begin{lemma} Suppose that
\beq
n_k(t)\leq C\varepsilon(1+ t)^{C\varepsilon}
+ C\varepsilon\int_0^t (1+\tau)^{-1}  n_k(\tau)\, d\tau
\eq
then
\beq
n_k(t)\leq C\varepsilon (1+t)^{2C\varepsilon}.
\eq
\end{lemma}
\begin{proof}
Let $N_k(t)=\int_0^t (1+\tau)^{-1} n_k(\tau)\, d\tau$.
Then $|\dot{N}_k(t)|
\leq C\varepsilon (1+t)^{-1} \big((1+ t)^{C\varepsilon}+N_k(t)\big)$.
Multiplying by an integrating factor gives
$\big(N_k(t) (1+t)^{-2C\varepsilon} \big)^{\prime}
\leq C\varepsilon (1+t)^{-1-C\varepsilon}$
so $N_k(t) (1+t)^{-2C\varepsilon} \leq C $
and hence $N_k(t)\leq C (1+t)^{2C\varepsilon}$
and   $n_k(t)\leq 2C\varepsilon (1+t)^{2C\varepsilon}$.
\end{proof}

\noindent
This proves the first part of \eqref{eq:sharpdecay4}. The
second part of \eqref{eq:sharpdecay4} follows from integrating the
first and using the boundary assumption \eqref{eq:boundcond}
as in the proof of \eqref{eq:weakdecay2}.

\section {Energy estimates for Einstein's equations}\label{section:energyeisnt}
Recall the definitions
\begin{align}
E_{N}(t)=\sup_{0\le\tau\le t}\sum_{|I|\le N}\int_{\Si_{t}} |\pa Z^{I} h|^{2},
\label{eq:energyN} \\
S_{N}(t) =
 \sum_{|I|\le N}\int_{0}^{t} \int_{\Si_{\tau}}
\frac {\gamma\,\, |\pab Z^{I} h|^{2}}{(1+|q|)^{1+2\gamma}}\label{eq:senergyN}
\end{align}
In this section we prove the following theorem.
\begin{theorem}\label{energyest}
 Assume that $g=h+m$ satisfies both Einstein's equations and the wave coordinate condition
 for $0\leq t\leq T$. Suppose also that for some $0< \gamma\leq
 1/2$ we have the following estimates for $0\leq t\leq T$:
\begin{enumerate}
\item For all multi-indices $I,\,\, |I|\le N/2+4$
\beq
| \pa Z^I h| + (1+|q|)^{-1} |Z^I h| + (1+t)(1+|q|)^{-1} |\pab Z^{I} h|
\leq C\varepsilon (1+t)^{-1/2-\ga} (1+|q|)^{-1/2-\gamma},\label{eq:weakpoint}
\eq
\item For all multi-indices $I,\,\, |I|\le N$
\beq
|Z^{I} H(s,q,\omega)|\le C \ve (1+t)^{-1}, \quad {\text{for}}\,\,\,
q=1, \label{eq:ZIboundary}
\end{equation}
\item
\begin{equation}
|\pa H|_{\cal TU} + (1+|q|)^{-1} | H|_{\cal TU}+ (1+|q|)^{-1}
|Z H|_{\cal LL}\le C\varepsilon
(1+t)^{-1},\label{eq:specialcomp}
\end{equation}
\item For all multi-indices $I,\,\, |I|\le N/2+2$ \beq | \pa Z^I h|
+ (1+|q|)^{-1} |Z^{I} h|\leq C\varepsilon
(1+t)^{-1+C\varepsilon},\label{eq:pointimprove} \eq
 \item
 \beq
E_N(0)\leq \ve^2.\label{eq:initialenergy}
 \eq
\end{enumerate}
Then there are  positive constants $C_{k}$ independent of $T$ such
that if $\ve \leq C_k^{-2}$ we have the energy estimate \beq\label{eq:energyestk} E_k(t)
+ S_{k}(t) \leq 16 \varepsilon^{2} (1+t)^{C_{k}\varepsilon},
\end{equation}
for $0\leq t\leq T$ and for all $k\leq N$.
\end{theorem}
\begin{remark}
Once again we recall that our estimates hold simultaneously for the tensors
$h$ and $H=-h+O(h^2)$. We shall freely interchange $h$ and $H$ in the
proof below.
\end{remark}
\begin{proof}
Recall that the components of the tensor $h_{\mu\nu} = g_{\mu\nu}
- m_{\mu\nu}$ satisfy the following wave equations:
\begin{align}
&g^{\alpha\beta}\pa_{\alpha}\pa_{\beta} h_{\mu\nu}=
F_{\mu\nu},\label{eq:equat},\\
&F_{\mu\nu}= P(\pa_{\mu} h,\pa_{\nu} h) +
Q_{\mu\nu}(\pa h,\pa h)+G_{\mu\nu}(h)(\pa h,\pa h).\nn
\end{align}
where \beq \label{eq:defiP} P(\pa_{\mu}h,\pa_{\nu} h) =
\frac{1}{4}  m^{\alpha\alpha^\prime}\pa_\mu
h_{\alpha\alpha^\prime} \,  m^{\beta\beta^\prime}\pa_\nu
h_{\beta\beta^\prime} -\frac{1}{2}
m^{\alpha\alpha^\prime}m^{\beta\beta^\prime} \pa_\mu
h_{\alpha\beta}\, \pa_\nu h_{\alpha^\prime\beta^\prime} \eq We
prove the desired estimate by induction on $k$. We first establish
the estimate \beq \label{eq:Ener0} E_{0}(t) + S_{0}(t) \le 8
\ve^{2} (1+t)^{C_{0}\ve} \ee for some constant $C_{0}$. After that
we shall assume that the statement \eqref{eq:energyestk} for
$k\leq N^\prime-1$ and prove the corresponding statement for $k\le
N^\prime$ with some constant $C_{N'}$. We shall base our argument
on the energy estimate \eqref{eq:firstenergy} for the solution of
the wave equation $\Boxr_{g}\phi =F$ proved in Proposition
\ref{prop:Decayenergy}. Observe that the conditions of our
Proposition on the  tensor $h=g-m$ imply the assumptions of
Proposition \ref{prop:Decayenergy} for the metric $g$.
\beq\label{eq:energyuse}
 \int_{\Si_{t}} |\pa\phi|^{2} + \int_{0}^{t} \int_{\Si_{\tau}}
\frac {\gamma\,\, |\pab\phi|^{2}}{(1+|q|)^{1+2\gamma}}\leq
 8\int_{\Si_{0}} |\pa \phi|^{2} +  C\varepsilon
\int_0^t\int_{\Si_{t}} \frac {|\pa\phi|^{2}}{1+t}
 + 16\int_0^t \int_{\Si_{t}}   |F|\, |\pa \phi|
\eq
\subsection{The case of $N'=0$.}
In this section we prove the basic energy estimate for a solution
of the equation \eqref{eq:equat}.
$$
\Boxr_{g} h_{\mu\nu}= F_{\mu\nu}:= P(\pa_{\mu} h,\pa_{\nu} h) +
Q_{\mu\nu}(\pa h,\pa h)+G_{\mu\nu}(h)(\pa h, \pa h).
$$
Recall that according to \eqref{eq:decayinhom2} of Lemma \ref{decayinhom}
we have a pointwise bound
$$
|F|\les |\pa h|_{\cal TU}^{2} + |\pab h||\pa h| + h |\pa h|^{2}
$$
Using the assumptions of the proposition we infer that
\beq
\label{eq:Fzero}
|F|\les \varepsilon \frac {|\pa h|}{1+t}
\end{equation}
Therefore, the energy estimate \eqref{eq:energyuse} with
$\phi=h_{\mu\nu}$ implies that
\begin{equation}
\label{eq:inductener0}
 \int_{\Si_{t}} |\pa h|^{2} + \int_{0}^{t} \int_{\Si_{\tau}}
\frac {\gamma \,\,|\pab h|^{2}}{(1+|q|)^{1+2\gamma}}\leq
 8\int_{\Si_{0}} |\pa h|^{2} +  C_0\varepsilon
\int_0^t\int_{\Si_{t}} \frac {|\pa h|^{2}}{1+t}.
\end{equation}
Using the smallness assumption on the initial data and the Gronwall
inequality this, in turn, leads to the desired estimate
\eqref{eq:Ener0}.
$$
E_{0}(t) + S_{0}(t) \le 8\ve^{2} (1+t)^{C_{0}\ve}
$$
\subsection{The case of $N'=1$.}
To facilitate the exposition we  first consider the
case $N'=1$.
We start by noting that according to \eqref{eq:curvwaveeqcommutest1}
of Corollary \ref{commut3} we have that
$$
\Boxr_{g} Z h_{\mu\nu} = \hat Z F_{\mu\nu} + D_{\mu\nu},
$$
where the term $D_{\mu\nu} = \Boxr_{g} Z h_{\mu\nu} - \hat Z
\Boxr_{g} h_{\mu\nu}$ satisfies the estimate
$$
|D|\les \Big (\frac {|ZH| + |H|}{1+t} + \frac {|ZH|_{\cal LL} +
|H|_{\cal LT}}{1+|q|}\Big )\sum_{|I|\le 1} |\pa Z^{I} h|
$$
Recall that the tensor $H^{\a\b}= - h^{\a\b} + O(h^{2})$. Thus
using the assumptions on $h$ of the proposition we derive that
$$
|D|\les \varepsilon \sum_{|I|\le 1} \frac {|\pa Z^{I} h|}{1+t}
$$
On the other hand, inequality \eqref{eq:decayinhom3} gives the
estimate
$$
|Z F|\leq \big(|\pa h|_{\cal TU}+|\overline{\pa} h|+|h|\,|\pa h|\big)
(|\pa Z h|+|\pa h|)
+C|\pa h |\,|\overline{\pa} Z h|
+C|\pa h |^2\,|Z h|
$$
Using the assumptions of the proposition we conclude that
$$
|\hat Z F|= |(Z+c_{Z}) F|\les  \varepsilon \sum_{|I|\le 1} \frac
{|\pa Z^{I} h|}{1+t} + \,\varepsilon\, \frac {|\pab Zh|}{
(1+t)^{\frac 12}(1+|q|)^{\frac 12 +\ga}}
$$
Now using the energy estimate \eqref{eq:energyuse} with
$\phi=Z h_{\mu\nu}$ and $F=\hat Z F_{\mu\nu} + D_{\mu\nu}$ we obtain
\begin{align*}
 \int_{\Si_{t}} |\pa Z h|^{2} + \int_{0}^{t} \int_{\Si_{\tau}}
\frac {\gamma\,\, |\pab Z h|^{2}}{(1+|q|)^{1+2\gamma}}& \leq
 8\int_{\Si_{0}} |\pa Z h|^{2} +  C\varepsilon  \sum_{|I|\le 1}
\int_0^t\int_{\Si_{t}}\frac {|\pa Z^{I}h|^{2}}{1+t} +
 C\varepsilon \int_0^t\int_{\Si_{t}} \frac {|\pab Zh|\, |\pa
 Zh|}{t^{\frac 12} (1+|q|)^{\frac 12 +\gamma}}\\ &\leq
 8\int_{\Si_{0}} |\pa Z h|^{2} +  C\varepsilon  \sum_{|I|\le 1}
\int_0^t\int_{\Si_{t}}\frac {|\pa Z^{I}h|^{2}}{1+t} +
C\varepsilon \int_0^t\int_{\Si_{t}}\frac {|\pab
Z^{I}h|^{2}}{(1+|q|)^{1+2\gamma}},
\end{align*}
where we used the Cauchy-Schwarz inequality to pass to the last line.
Combining this with  the energy inequality \eqref{eq:inductener0}
we infer that if $C\varepsilon\leq \gamma/2$ then
\begin{equation}
\label{eq:inductener1}
\sum_{|I|\le 1}\int_{\Si_{t}} |\pa Z^{I} h|^{2} +
\sum_{|I|\le 1} \int_{0}^{t} \int_{\Si_{\tau}}
\frac {\gamma\,\, |\pab Z^{I} h|^{2}}{(1+|q|)^{1+2\gamma}} \leq
16\sum_{|I|\le 1}\int_{\Si_{0}} |\pa Z^{I} h|^{2} +
C_1\varepsilon  \sum_{|I|\le 1}
\int_0^t\int_{\Si_{t}}\frac {|\pa Z^{I}h|^{2}}{1+t}
\end{equation}
The desired estimate
$$
E_{1}(t) + S_{1}(t) \le 16 \ve^{2} (1+t)^{C_{1}\ve}
$$
now follows from the Gronwall inequality and the smallness assumption
on the initial data.

\subsection{The case of $N'> 1$.}
In what follows we assume that we have already shown that \beq
E_{N'-1}(t) + S_{N'-1}(t)\le 16\ve^{2}
(1+t)^{C_{N'-1}\ve},\label{eq:enerN'} \eq and prove that there
exists a constant $C_{N'}$ such that \beq E_{N'}(t) + S_{N'}(t)\le
16\ve^{2} (1+t)^{C_{N'}\ve}, \label{eq:enerind} \eq
 We start this
section by writing the wave equation for the quantity $Z^{I}
h_{\mu\nu}$ with $|I|=N'$
$$
\Boxr_{g} Z^{I} h_{\mu\nu} = \hat Z^{I} F_{\mu\nu} +
D^{I}_{\mu\nu},
$$
where
$$
D^{I}_{\mu\nu}=\Boxr_{g} Z^{I} h_{\mu\nu}-\hat Z^{I}\Boxr_{g}
h_{\mu\nu}
$$

We apply the energy estimate \eqref{eq:energyuse} with the
functions $\phi=Z^{I} h_{\mu\nu}$ and $F=\hat Z^{I} F_{\mu\nu} +
D^{I}_{\mu\nu}$
 \beq\label{eq:energyuse2} \int_{\Si_{t}} |\pa
Z^{I}h|^{2} + \int_{0}^{t} \int_{\Si_{\tau}} \frac {\gamma\,\,
|\pab Z^{I}h|^{2}}{(1+|q|)^{1+2\gamma}}\leq
 8\int_{\Si_{0}} |\pa Z^{I}h|^{2} +  C\varepsilon
\int_0^t\int_{\Si_{t}} \frac {|\pa Z^{I} h|^{2}}{1+t}
 + 16\int_0^t \int_{\Si_{t}}   \big (|\hat Z^{I} F| + |D^{I}|\big )\,
|\pa Z^{I}h|
 \eq
 Note that we can estimate
 \beq \label{eq:energypart}
 \int_0^t\int \big(|\hat{Z}^I F|+|D^I|\big)\, |\pa Z^I h|\, dx\, dt\les
\int_0^t\frac{\varepsilon}{1+t}|\pa Z^I h|^2\, dx \, dt+\int_0^t
\int \varepsilon^{-1}(1+t)\big(|\hat{Z}^I F|^2+|D^I|^2\big)\, dx\,
dt \eq
 Here the first term is of the type that appears already in the energy estimate
 \eqref{eq:energyuse2}.
  Thus it remains to handle the second term.

According to \eqref{eq:curvwaveeqcommutest2} of Corollary \ref{commut3}
we have that
\begin{align}
&D^I=\sum_{k=0}^{|I|} D^{I}_{k},\\
&D^{I}_k=D^{I1}_k + D^{I2}_k + D^{I3}_k+D^{I4}_k,
\label{eq:commutators}\\
&|D^{I1}_k|\les \sum_{|K|=k}\,\,\,\,\sum_{|J|+(|K|-1)_+\leq |I|}
\frac{|Z^{J} H|}{1+t+|q|}\,\, {|\pa Z^{K} h|},\label{eq:commute1}\\
&|D^{I2}_k|\les \sum_{|K|=k}\,\,\,\, \sum_{|J|+(|K|-1)_+\leq |I|}
\frac{|Z^{J} H|_{LL}}{1+|q|}\,{|\pa
Z^{K} h|},\label{eq:commute2}\\
&|D^{I3}_k|\les \sum_{|K|=k}\,\,\,\,\sum_{|J|+(|K|-1)_+\leq |I|-1}
\frac {|Z^{J} H|_{L\cal T}}{1+|q|}\,{|\pa Z^{K} h|},\label{eq:commute3}\\
&|D^{I4}_k|\les \sum_{|K|=k}\,\,\,\,\sum_{|J|+(|K|-1)_+\leq |I|-2}
\frac {|Z^{J} H|}{1+|q|}\,{|\pa Z^{K} h|},\label{eq:commute4}
\end{align}
 {\bf The estimates for $D^I_k$ with $k\le N/2$.}
 We must now estimate
 \beq \int_0^t \int \varepsilon^{-1}(1+t)|D^I_k|^2\, dx\, dt \eq
 Since $k=|K|\leq N/2$
 in \eqref{eq:commute1}-\eqref{eq:commute4} it follows from the assumptions
 in the theorem that we can estimate
 \beq
\varepsilon^{-1} (1+t)|\pa Z^K h|^2\les
\min{\Big(\frac{\varepsilon}{(1+t)^{1-C\varepsilon}},
\frac{\varepsilon}{(1+|q|)^{1+2\gamma}}\Big)}
 \eq
and it thus suffices to estimate
 \begin{align}\label{eq:DkI1}
&\int_0^t\int \varepsilon^{-1}(1+t)|D^{I1}_k|^2\, dx\, dt \les
 \sum_{|J|\leq |I|}
 \int_0^t\int\frac{\varepsilon}{(1+|q|)^{1+2\gamma}} \frac{|Z^J
 H|^2}{(1+t+|q|)^2}\, dx\, dt,\\
 \label{eq:DkI34}
 &\int_0^t\int \varepsilon^{-1}(1+t)\big(|D^{I3}_k|^2+|D^{I4}_k|^2)\,
dx\, dt \les
 \sum_{|J|\leq |I|-1}
 \int_0^t\int\frac{\varepsilon}{(1+t)^{1-C\varepsilon}} \frac{|Z^J
 H|^2}{(1+|q|)^2}\, dx\, dt,\\
 \label{eq:DkI2}
 &\int_0^t\int \varepsilon^{-1}(1+t)|D^{I2}_k|^2\, dx\, dt \les
 \sum_{|J|\leq |I|}
 \int_0^t\int \min{\Big(\frac{\varepsilon}{(1+t)^{1-C\varepsilon}},
\frac{\varepsilon}{(1+|q|)^{1+2\gamma}}\Big)}\frac{|Z^J
 H|_{\cal LL}^2}{(1+|q|)^2}\, dx\, dt,
 \end{align}
 \begin{lemma} Let $f$ be a smooth function satisfying the condition
 \beq\label{eq:boundcondrepeat}
|f|\les \varepsilon (1+t)^{-1},\qquad \text{for}\qquad q=1
 \eq
 Then
 \beq \label{eq:commuteenergy1} \int_0^t \int
\frac{\varepsilon}{(1+|q|)^{1+2\gamma} }\frac{|f
|^2}{(1+t+|q|)^2}\, dx\, dt \les \int_0^t
\frac{\varepsilon}{(1+t)^{1+2\gamma}}\int |\pa
f|^2\, dx \, dt +\varepsilon^3\eq
 and
\beq \label{eq:commuteenergy2} \int_0^t \int
\frac{\varepsilon}{(1+t)^{1-C\varepsilon}}\frac{|f|^2}{(1+|q|)^2}\,
dx\, dt \les \int_0^t
\frac{\varepsilon}{(1+t)^{1-C\varepsilon}}\Big(\varepsilon^2+\int
|\pa f|^2\, dx \Big)\, dt\eq Furthermore,
\begin{align} \int_0^t \int
\min{\Big(\frac{\varepsilon}{(1+t)^{1-C\varepsilon}},
\frac{\varepsilon}{(1+|q|)^{1+2\gamma}}\Big)}\frac{|f|^2}{(1+|q|)^2}\,
dx\, dt &\les \int_0^t \int \frac{\varepsilon|\pa_r f|^2
}{(1+|q|)^{1+2\gamma}}\, dx dt\nn\\ &+ \varepsilon^2
\int_0^t \frac{\varepsilon}
{(1+t)^{1-C\varepsilon}}\,dt\label{eq:commuteenergy3}
\end{align}
 \end{lemma}
 \begin{proof}
We shall repeatedly use the Poincar\'e inequality \eqref{eq:Poinc}
of Lemma \ref{le:Poinc} \beq \label{eq:Poinc2} \il_{\Si_{t}} \frac
{|f(x)|^{2}\, dx}{(1+|q|)^{2+2\si}} \les \il_{S_{(t+1)}} |f|^{2}\,
dS + \il_{\Si_{t}} \frac {|\pa_{r} f(x)|^{2}\, dx}{(1+|q|)^{2\si}}
\ee which holds for any value of $\si >-1/2,\,\si\ne 1/2$. In
particular, using \eqref{eq:boundcondrepeat}, we
obtain that
\beq \label{eq:PoincZI} \il_{\Si_{t}} \frac
{|f|^{2}\, dx}{(1+|q|)^{2+2\si}} \les \varepsilon^{2}+
\il_{\Si_{t}} \frac {|f|^{2}\, dx}{(1+|q|)^{2\si}} \eq
The estimates \eqref{eq:commuteenergy1} and
\eqref{eq:commuteenergy2} now
follow from \eqref{eq:PoincZI} with $\sigma=0$.
\end{proof}
We now note the following generalization of \eqref{eq:Poinc2} \beq
\label{eq:Poinc3} \il_{\Si_{t}}
\min{\Big(\frac{\varepsilon}{(1+t)^{1-C\varepsilon}},
\frac{\varepsilon}{(1+|q|)^{1+2\gamma}}\Big)}\frac {|f(x)|^{2}\,
dx}{(1+|q|)^{2}} \les
\frac{\varepsilon}{(1+t)^{1-C\varepsilon}}\il_{S_{(t+1)}}
|f|^{2}\, dS + \varepsilon \il_{\Si_{t}} \frac {|\pa_{r} f(x)|^{2}\,
dx}{(1+|q|)^{1+2\gamma}} \ee
 The proof of \eqref{eq:Poinc3} can be reduced to \eqref{eq:Poinc2} by
 subtracting a term which picks up the boundary value. We define
 \beq
\widetilde{f}=f-\overline{f},\qquad\text{where}\qquad
\overline{f}(r,\omega)=f\big((t+1),\omega\big)\chi(r/t)
 \eq
 and $\chi(s)=1$, when $3/4\leq s\leq 3/2$ and $\chi(s)=0$ when
 $s\leq 1/2$ or $s\geq 2$. Then
\begin{align}
\il_{\Si_{t}}
\min{\Big(\frac{\varepsilon}{(1+t)^{1-C\varepsilon}},
\frac{\varepsilon}{(1+|q|)^{1+2\gamma}}\Big)}\frac {|f(x)|^{2}\,
dx}{(1+|q|)^{2}} \les \varepsilon \il_{\Si_{t}}
\frac {|\tilde f(x)|^{2}\,dx}{(1+|q|)^{3+2\ga}} +
\frac{\varepsilon}{(1+t)^{1-C\varepsilon}}\il_{\Si_{t}}
\frac {|\bar f(x)|^{2}\,
dx}{(1+|q|)^{2}}
\end{align}
We now apply \eqref{eq:Poinc2} to the function $\tilde f$,
which vanishes at $r=t+1$, and observe that
\beq
\il_{\Si_{t}} \frac {|\pa_{r} \overline{f}(x)|^{2}\,
dx}{(1+|q|)^{1+2\gamma}}\les \int
f\big((t+1),\omega\big)^2 d\omega
\les\frac{1}{(1+t)^{2}}\int_{S_{t+1}}|f|^2\, dS
 \eq
On the other hand,
\beq \frac{\varepsilon}{(1+t)^{1-C\varepsilon}}\,
\int_{\Si_{t}}\frac {|\overline{f}(x)|^{2}\, dx}{(1+|q|)^{2}} \les
\frac{\varepsilon}{(1+t)^{1-C\varepsilon}} \int_{S_{t+1}}|f|^2\,
dS\eq which proves \eqref{eq:Poinc3}.

Using the lemma above with $f= Z^J H$, together with \eqref{eq:DkI1},
\eqref{eq:DkI34}
and the assumption that $E_{N'-1}\leq
16(1+t)^{C_{N'-1}\varepsilon}$ we see that we can estimate
 \begin{align*}
\int_0^t\int
\varepsilon^{-1}(1+t)\big(|D^{I1}_k|^2+|D^{I3}_k|^2+|D^{I4}_k|^2)\, dx\,
dt &\les \int_0^t\frac{\varepsilon }{(1+t)^{1+2\gamma}}E_{N'}(t)\,
dt+\varepsilon^2\int_0^t
\frac {\varepsilon}{(1+t)^{1-C\varepsilon}}\,dt\nn\\
&\les \varepsilon\, E_{N'}(t) + \varepsilon^2\int_0^t
\frac {\varepsilon}{(1+t)^{1-C\varepsilon}}\,dt
\end{align*}
for all $\le N/2$.

It thus remains the term \eqref{eq:DkI2} containing $D^{I2}_k$. We
shall use the version of the Poincar\'e inequality
\eqref{eq:commuteenergy3}  to create the term $\pa_{q} (Z^{J}
H)_{LL}$, which can be then converted to a tangential derivative
of $Z^{J} H$ via the wave coordinate condition. However, in order
to implement this strategy we modify the term $Z^{J} H_{LL}$
according to Lemma \ref{decaywavectwo}. We recall the notation
\beq\label{eq:HJ}
 H^{J}_{\mu\nu}=Z^{J} H_{\mu\nu} +
\sum_{|J'|<|J|} c_{J^\prime\, \, \mu}^{\, J \,\,\,\ga } \, Z^{J}
H_{\gamma \nu} \eq
 If $|J|\leq N^\prime $ then the lower order terms in the right
 hand side of \eqref{eq:HJ} may be estimated using
 \eqref{eq:DkI34} and \eqref{eq:commuteenergy2} as before.
 According to  Lemma  \ref{decaywavectwo} and the pointwise
estimates in \eqref{eq:pointimprove} and \eqref{eq:ZIboundary}
\begin{multline}\label{eq:HJr}
|\pa_r H^{J}_{LL}|\les \sum_{|J^\prime|\le |J|} |\pab Z^{J^\prime}
H|+ \sum_{|J_{1}|+..+|J_{m}|\leq |J|, \, \, m\ge 2} |Z^{J_{m}}
H|\cdots |Z^{J_{2}} H|\, |\pa Z^{J_{1}} H|\\
  \les \sum_{|J^\prime|\le
|J|} |\pab Z^{J^\prime} H|+ \sum_{|J_{1}|+|J_{2}|\leq |J|}
|Z^{J_{1}} H|\,|\pa Z^{J_{2}} H|\\
\les \sum_{|J^\prime|\le |J|} |\pab Z^{J^\prime} H|+
\frac{\varepsilon(1+|q|)^{1/2-\gamma}}{(1+t)^{1/2+\gamma}} |\pa
Z^{J^\prime} H|+ \frac{\varepsilon|Z^{J^\prime}
H|}{(1+t)^{1/2+\gamma}(1+|q|)^{1/2+\gamma}}
\end{multline}
Hence \beq \int_0^t \int \frac{\varepsilon|\pa_r H^J_{LL}|^2
}{(1+|q|)^{1+2\gamma}}\, dx dt \les \sum_{|J^\prime|\leq |J|}
\int_0^t \int \Big(\frac{\varepsilon|\overline{\pa} Z^{J^\prime}
H|^2 }{(1+|q|)^{1+2\gamma}}  + \frac{\varepsilon|{\pa}
Z^{J^\prime} H|^2 }{(1+t)^{1+2\gamma}} + \frac{\varepsilon|
Z^{J^\prime} H|^2 }{(1+t)^{1+2\gamma}(1+|q|)^{2}}\Big)\, dx dt \eq
If we use \eqref{eq:commuteenergy2} with $C\varepsilon$ in the
exponent replaced by $2\gamma$ we see that the last term can be
estimated by the second term from the right plus a term from the
boundary: \beq\label{eq:finis2}
 \int_0^t \int \frac{\varepsilon|\pa_r
H^J_{LL}|^2 }{(1+|q|)^{1+2\gamma}}\, dx dt \les
\sum_{|J^\prime|\leq |J|} \int_0^t \int
\Big(\frac{\varepsilon|\overline{\pa} Z^{J^\prime} H|^2
}{(1+|q|)^{1+2\gamma}}  + \frac{\varepsilon|{\pa} Z^{J^\prime}
H|^2 }{(1+t)^{1+2\gamma}} + \frac{\varepsilon}
{(1+t)^{1+2\gamma}}\varepsilon^2\Big)\, dx dt \eq
 As we argued,
when estimating \eqref{eq:DkI2} we can replace $|Z^J H|_{LL}$ by
the left hand side of \eqref{eq:HJ}. After that we use the version
of the Poincare inequality \eqref{eq:commuteenergy3} applied to
$H^J_{LL}$ and this together with \eqref{eq:finis2} gives
 \beq
\int_0^t \int \varepsilon^{-1} (1+t)|D^{I2}_k|^2\, dx dt \les
\varepsilon S_{N^\prime}(t)+\varepsilon E_{N^\prime}(t)+
\varepsilon^2\int_0^t\frac{\varepsilon}{(1+t)^{1-C\varepsilon}}\,dt
 \eq
 Summarizing, we have proven that
\beq \label{eq:commutatorone}\int_0^t \int \varepsilon^{-1}
(1+t)|D^{I}_k|^2\, dx dt \les \varepsilon
S_{N^\prime}(t)+\varepsilon E_{N^\prime}(t)+
\varepsilon^2\int_0^t\frac{\varepsilon}{(1+t)^{1-C\varepsilon}}\,dt,
\qquad k\leq N/2
 \eq
\noindent
This concludes the estimates in the case $k\le N/2$.

{\bf The commutator in case $k\ge N/2$.}
 \noindent We isolate the case when $|K|=N^\prime=|I|$.
We can estimate its contribution to the $D^{I}_{N^\prime}$ by the
following expression:
$$
|D^{I}_{N^\prime}|\les \sum_{|K|=|I|}\Big (\frac {|H| +
|ZH|}{1+t+|q|} + \frac {|ZH|_{\cal LL} +
 |H|_{\cal LT}}{1+|q|} \Big )|\pa Z^{K} h|\les \ve
\sum_{|K|=|I|}\frac {|\pa Z^{K} h|}
{1+t},
$$
where to pass to the last line we used pointwise estimates from
\eqref{eq:specialcomp}, \eqref{eq:weakpoint}, and
\eqref{eq:ZIboundary}. In the case when $N/2\leq k<|I|$ we
estimate the contribution of the corresponding term in $D^{I}_k$,
with the help of \eqref{eq:pointimprove} as follows:
$$
|D^{I}_k|\les \sum_{|K|<|I|} \sum_{|J|\le N/2} \frac {|ZH|}{1+|q|}
|\pa Z^{K} h|\les \ve \sum_{|K|<|I|}\frac {|\pa Z^{K}
h|}{(1+\tau)^{1-C\ve}}
$$
Therefore,
 \beq \int_{0}^{t}\int_{\Si_{\tau}}
\varepsilon^{-1}(1+t)|D^{I}_k|^2\, dx dt \les \ve \int_0^t \int
\sum_{|K|<|I|} \frac {|\pa Z^{K} h|^2}{(1+\tau)^{1-2C\ve}}
+\sum_{|K|=|I|} \frac {|\pa Z^{K} h|^2}{1+\tau}\, dx dt
 \eq
 Using the inductive assumption \eqref{eq:enerN'} we can therefore
 estimate
\beq \label{eq:commutatortwo}
\int_{0}^{t}\int_{\Si_{\tau}}
\varepsilon^{-1}(1+t)|D^{I}_k|^2\, dx dt \les \varepsilon \int_0^t
\frac {E_{N^\prime}(\tau)}{1+\tau}\, dt+ \ve^2 \int_0^t \frac
{\varepsilon\, dt}{(1+\tau)^{1-2C\ve}} ,\qquad N/2\leq k\leq
N^\prime
 \eq

{\bf The inhomogeneous term.}
By \eqref{eq:decayinhom4}
\begin{multline}
|\hat{Z}^I F| \les
 \sum_{|J|+|K|\le |I|}
\big (|\pa Z^J h|_{\cal TU}\, |\pa Z^K h |_{\cal TU} +
|\pab Z^J h|\, |\pa Z^K h| \big )+ \sum_{|J|+|K|\le |I|-1}
|\pa Z^J h|_{L \cal T} |\pa Z^K h|\\  +
\sum_{|J|+|K|\le |I|-2} |\pa Z^J h|\, |\pa Z^K h| +
\sum_{|J_1|+...+|J_m| \leq |I|,\, m\geq 3}
|Z^{J_m} h|\!\cdot\cdot\cdot |Z^{J_3} h|\,|\pa Z^{J_{2}} h| |\pa Z^{J_1} h|
\end{multline}
The highest order terms with one of $|J|$, $|K|$ or $|I_i|$ equal to $N=|I|$
are bounded by
\begin{multline}
\big(|\pa h|_{\cal T U}+|\pab h|+|h||\pa h|\big)
\sum_{|I|=N}|\pa Z^I h|
+|\pa h|^2\sum_{|I|=N} |Z^I h|+|\pa h|\, \sum_{|I|=N}|\pab Z^I h|\\
\leq \frac{\ve}{1+t}\sum_{|I|=N}|\pa Z^I h| +
\frac{\ve^2}{(1+t)^{1+2\gamma}(1+q)^{1+2\gamma} }\sum_{|I|=N} |Z^I
h| + \frac{\ve}{(1+t)^{1/2+\gamma}(1+q)^{1/2+\ga} }\sum_{|I|=N}
|\pab Z^I h|
\end{multline}
The remaining lower order terms are of the form
\begin{multline}
\sum_{|K|<N,\, |J|\leq N/2 }|\pa Z^J h|\, |\pa Z^K h|
+\sum_{|K|<N,\,\, |J|,\,|L|\leq N/2 }|\pa Z^J h|\,|\pa Z^L h|\, |Z^K h| \\
\leq \frac{\ve }{(1+t)^{1-C\ve} }\sum_{|K|<N}|\pa Z^K h| +
\frac{\ve^2}{(1+t)^{1+2\gamma}(1+q)^{1+2\gamma} }\sum_{|I|<N} |Z^I
h|.
\end{multline}

It therefore follows that
 \begin{multline}\int_{0}^{t}\int_{\Si_{\tau}}
\varepsilon^{-1}(1+t)|\hat{Z}^I F|^2\, dx dt \les \sum_{|K|\leq
|I|} \ve \int_0^t\int \frac {|\pa Z^{K} h|^2}{1+\tau}+ \frac
{|\overline{\pa} Z^{K} h|^2}{(1+|q|)^{1+2\gamma}}
+\frac {| Z^{K} h|^2}{(1+t)^{1+2\gamma}(1+|q|)^{2}}\, dx dt\\
 +\int_0^t \frac {\varepsilon\, dt}{(1+\tau)^{1-2C\ve}}
\sum_{|I|<N} |\pa Z^I h|^2\, dx dt\\
\les \sum_{|K|\leq |I|} \ve \int_0^t\int \frac {|\pa Z^{K}
h|^2}{1+\tau}+ \frac {|\overline{\pa} Z^{K}
h|^2}{(1+|q|)^{1+2\gamma}} +\frac
{\varepsilon^2}{(1+t)^{1+2\gamma}}\, dx dt+
 \int_0^t \frac {\varepsilon\, dt}{(1+\tau)^{1-2C\ve}}
\sum_{|I|<N} |\pa Z^I h|^2\, dx dt
 \end{multline}
 Here, to estimate the last term in the first row we used
 \eqref{eq:commuteenergy2} with $-C\ve$ in the exponent replaced by $2\gamma$,
 which produced a term similar to the first term of the first line plus a boundary term.
 Using the inductive assumption \eqref{eq:enerN'} we thus obtain
 \beq\label{eq:inhomenergy}
\int_{0}^{t}\int_{\Si_{\tau}} \varepsilon^{-1}(1+t)|\hat{Z}^I
F|^2\, dx dt \les \ve\int_0^t \frac{E_{N^\prime}(\tau)\,
d\tau}{1+\tau}+\ve S_{N^\prime}(t)+\ve^2\int_0^t\frac{\ve\,
d\tau}{(1+\tau)^{1-C\ve}}
 \eq

{\bf The conclusion of the proof in case $N^\prime>1$}
 The inequalities \eqref{eq:energyuse2}-\eqref{eq:energypart} and
 \eqref{eq:commutatorone}, \eqref{eq:commutatortwo} and
 \eqref{eq:inhomenergy} imply that
for some constant $C$: \beq\label{eq:energyuse3}
E_{N^\prime}(t)+S_{N^\prime}(t)\leq 8
E_{N^\prime}(0)+C\ve\big(E_{N^\prime}(t)+S_{N^\prime}(t)\big)
 +C\ve \int_0^t\frac{E_{N^\prime}(\tau)\, d\tau}{1+\tau}
 +C\ve^2\int_0^t\frac{\ve\,
 d\tau}{(1+\tau)^{1-C\ve}}
 \eq
 If we now choose $\ve$ so small that $C\ve\leq 1/9$ we can move
 the second term on the right to the left and multiply by $9/8$ to
 obtain for some new constants
\beq\label{eq:energyuse4} E_{N^\prime}(t)+S_{N^\prime}(t)\leq 9
E_{N^\prime}(0)
 +C\ve \int_0^t\frac{E_{N^\prime}(\tau)\, d\tau}{1+\tau}
 +C\ve^2\int_0^t\frac{\ve\,
 d\tau}{(1+\tau)^{1-C\ve}}
 \eq
 This can now be integrated using a Gr\"onwall type of argument.
 If $G(t)$ denotes the right hand side then we have
 $$
 G^{\,\prime}(t)\leq\frac{C\ve}{1+t}G(t)+\frac{C\ve^3}{(1+t)^{1-C\ve}}
 $$
 Multiplying with the integrating factor we get
 $$
 \frac{d}{dt}\Big( G(t) (1+t)^{-C\ve}\Big)
 \le \frac{C\ve^3}{1+t}
 $$
 and hence if integrate and use that $ C\ve \ln{(1+t)}\leq
 (1+t)^{C\ve}$, for $t\geq 0$ (as is seen by differentiating both
 sides), and use that by assumption \eqref{eq:initialenergy}
 $G(0)\leq 9\ve^2$,
 we obtain
 $$
 G(t)\leq G(0) (1+t)^{C\ve}+ C\ve^3 \ln{(1+t)} (1+t)^{C\ve}\leq
9\ve^2(1+t)^{C\ve}+ \ve^2 (1+t)^{2C\ve}\leq 10\ve^2(1+t)^{2C\ve}
$$
Hence we have proven that
$$
E_{N^\prime}(t)+S_{N^\prime}(t)\leq 10\ve^2 (1+t)^{2C\ve}
$$
This concludes the induction and the proof of the theorem.
\end{proof}

\section{Geodesic completeness}
Having constructed a solution metric $g=m+h$ of the Einstein
equations we need to verify that the resulting space-time $({\Bbb R}^4, g)$
is causally geodesically complete.
Let
$$
X(\tau)= (x^{0}(\tau),x(\tau))=
(t(\tau), x(\tau) = (t(\tau), r\omega(\tau))
$$
be a causal geodesic parameterized by the affine parameter $\tau$.
Such geodesics satisfy the equations
\begin{align}
&\ddot X^{\a}(\tau) + \Gamma^{\a}_{\b\ga}(X(\tau))
\dot X^{\b} \dot X^{\ga} = 0, \label{eq:geod}\\
& X(0)= Y,\quad \dot X(0)=\xi\nn
\end{align}
where $Y$ is the point of the origin of the geodesic $X(\tau)$ and
$\xi$ is the initial velocity satisfying the condition \beq
\label{eq:timelike} g_{\a\b}(Y) \xi^{\a}\xi^{\b} = -A^{2}\le 0 \ee
for some constant $A$. Condition \eqref{eq:timelike} is preserved
in time, i.e., \beq \label{eq:time} g_{\a\b} (X(\tau)) \dot X^{\a}
\dot X^{\b} =-A^{2} \ee In the following lemma we show that a
vector $\eta$ causal with respect to the metric $g$ is "almost"
causal with respect to the Minkowski metric $m$.
\begin{lemma} \label{le:Cause}
Let $\eta$ be a causal 4-vector, i.e,
\beq
\label{eq:causal}
g_{\a\b} \eta^{\a} \eta^{\b} \le - A^{2} \le 0
\ee
for some non-negative constant $A$.
Then
\beq
\label{eq:timecontrol}
A + |\eta^{i}|\le 2 |\eta^{0}|,\qquad \forall i=1,..,3
\ee
\end{lemma}
\begin{proof}
Expanding $g= m + h$ we obtain from \eqref{eq:causal}
that
$$
-|\eta^{0}|^{2} + \sum_{i=1}^{3} |\eta^{i}|^{2} \le
|h|\cdot (|\eta^{0}|^{2} +\sum_{i=1}^{3} |\eta^{i}|^{2})
$$
and the desired estimate follows provided that
$|h|\le 1/4$.
\end{proof}
\noindent
We choose a future oriented initial velocity $\xi$,
i.e., $\dot x^{0}(0) > 0$.
\begin{prop}\label{pr:Inext}
Assume that $h=g-m$ satisfies the estimates\footnote{These assumptions are
consistent with the decay estimates for $h$ proved in Theorem \ref{decayeinst}.}
\begin{align*}
&|h| |\pa h|+|\pa h|_{TU} +
|\pab h|_{\Lb\Lb}\les \ve t^{-1},\\
&|\pa h(t,x)|\les \ve t^{-1},\quad \text{for}\quad
|x|\le t/2
\end{align*}
Let $X(\tau)$ is a future inextendible causal geodesic. Then the
values of the affine parameter $\tau$ span the interval $[0,\infty)$.
\end{prop}
\begin{proof}
We start by considering a time-like geodesic $X(\tau)$.
Reparameterizing, if necessary, we can assume that the constant
$A=1$ in \eqref{eq:time}. Then equation \eqref{eq:time} and
inequality \eqref{eq:timecontrol} with $A=1$ imply that for all
$\tau \ge 0$. \beq \label{eq:timemon} \dot x^{0}(\tau) \ge \frac
12  + |\dot x(\tau)| \ee We removed the absolute value from $\dot
x^{0}(\tau)$, since $\dot x^{0}(0)>0$. This is the only part of
the argument, which uses the fact that $X(\tau)$ is a time-like
geodesic. The case of a null geodesic will require an additional
argument.

Assume that $X(\tau)$ is a time-like geodesic of finite
length $\tau_{*}$.
We first observe that
$$
\lim_{\tau\to \tau_{*}} |X(\tau)|=\infty
$$
which means that $X(\tau)$ escapes to infinity\footnote{viewed from the point
of view of the global system of wave coordinates on ${\Bbb R}^4$.}  in finite
proper time $\tau_{*}$. This easily follows from the standard
ODE theory.
The inequality \eqref{eq:timemon} implies that
$\dot x^{0} (\tau)$ controls $\dot X(\tau)$.
Thus to obtain contradiction it suffices to show that
$$
\lim_{\tau\to \tau_{*}} x^{0}(\tau) <\infty
$$
Throughout this section we will
use consistently use the notation $x^{0} = t$.
We recall that
$$
\Gamma^{\a}_{\b\ga} = g^{\a\si}(\pa_{\b} h_{\ga\si} +
\pa_{\ga} h_{\b\si} -\pa_{\si} h_{\b\ga} )
$$
Thus, expanding the metric $g=m+h$,
$$
\ddot x^{0} - (2 \pa_{\b} h_{0\ga} - \pa_{0} h_{\b\ga} )
\dot x^{\b} \dot x^{\ga} = h \cdot \pa h \cdot |\dot X|^{2}
$$
We further observe that
\beq
\label{eq:fulld}
\pa_{\b} h_{0\ga} \dot x^{\b} \dot x^{\ga} =
\frac d{d\tau} \big (h_{0\ga} \dot x^{\ga} ) -
h_{0\ga} \ddot x^{\ga}
\ee
We now additionally recall that $\pa_{q} h_{\Lb\Lb}$ is the only
derivative of
$h$ that does not have the decay rate of at least $(x^{0})^{-1}$.
Thus
$$
\pa_{0} h_{\b\ga} \dot X^{\b} \dot X^{\ga} =
\pa_{q} h_{\b\ga} \dot X^{\b} \dot X^{\ga} +
\ve O(t^{-1}) |\dot X|^{2} = \pa_{q} h_{\Lb\Lb} |\dot X^{\Lb}|^{2} +
\ve O((x^{0})^{-1}) |\dot X|^{2}
$$
The expression
$$
\dot X^{\Lb} = \dot X^{\a} L_{\a} = - \dot x^{0} + \frac
{x^{i}}{|x|} \dot x^{i} = -\frac d{d\tau} (t-r) = -\dot q
$$
Moreover,
$$
\pa_{q} h_{\Lb\Lb} = 4\pa_{q} h_{00} + \ve O((x^{0})^{-1})
$$
Furthermore, introduce $\zeta (x^{0}/r)$ a cut-off function of
the set $r\ge x^{0}/2$. Then
$$
\pa_{q} h_{00} = (1-\zeta)\pa_{q} h_{00} + \zeta \pa_{q} h_{00}
= \ve O(t^{-1})  + \pa_{q} (\zeta h_{00}) -
(\pa_{q} \zeta (x^{0}/r)) h_{00}
$$
We compute
$$
\pa_{q}\zeta(x^{0}/r) = (r^{-1} + x^{0} r^{-2}) \zeta'(x^{0}/r) \les
(x^{0})^{-1}
$$
since $r\ge x^{0}/2$ on the support of $\zeta'(x^{0}/r)$. Thus
$\pa_{q} h_{00}$ can be replaced by $\pa_{q} (\zeta h_{00})$
at the expense of a term of order $\ve O((x^{0})^{-1})$.
Therefore,
\begin{align*}
\pa_{q} h_{\Lb\Lb}  |\dot X^{\Lb} |^{2} &=
4\pa_{q} h_{00} |\dot q|^{2}  = 4\pa_{q} (\zeta h_{00}) |\dot q|^{2} + \ve
O((x^{0})^{-1}) |\dot q|^{2}\\ & =
4\frac d{d\tau} \big (\zeta h_{00} \,\dot q\big ) -
4\zeta h_{00} \ddot q - 4\pa_{L} (\zeta h_{00}) \dot X^{L}\dot X^{\Lb} -
4\pa_{\Omega} (\zeta h_{00} )\dot X^{\omega} \dot X^{\Lb} +\ve O((x^{0})^{-1})
\end{align*}
Here,
$$
h(X(\tau)) = h(q((\tau), v(\tau), \omega (\tau))
$$
where $q=x^{0}-r$, $v=x^{0}+r$, and $\omega = \frac {x^{i}}{r}$.
The advantage is that $\pa_{\Omega} h_{00}, \pa_{L} h_{00}$
already decay faster than $(x^{0})^{-1}$ and
$\pa_{\Omega} \zeta (x^{0}/r) =0$,
while $|\pa_{L} \zeta (x^{0}/r) |\les (x^{0})^{-1}$.
Thus
$$
\pa_{q} h_{\Lb\Lb}  |\dot X^{\Lb} |^{2} =
\frac d{d\tau} \big (\zeta h_{00} \dot q\big ) -
\zeta h_{00} \ddot q + \ve O((x^{0})^{-1})
$$
It remains to analyze the term
\beq
\label{eq:dotq}
\ddot q = \frac d{d\tau} (\dot x^{0} - \frac {x^{i}}{r} \dot x^{i})
=\ddot x^{0} - \ddot x^{i} \frac {x^{i}}r +
r^{-1} (|\dot x|^{2} - r^{-2}|x\cdot \dot x|^{2})
\ee
From the geodesic equation \eqref{eq:geod} we can estimate
$$
|\ddot x^{\a}|\le |\pa h| |\dot X|^{2}
$$
Additionally, since on the support of $\zeta (x^{0}/r)$,
$r\ge x^{0}/2$, we have that the last term in \eqref{eq:dotq}
multiplied by $\zeta h_{00}$ contributes at most\footnote{This
is the reason for introducing the cut-off function $\zeta$.}
$\ve (x^{0})^{-1}$.
Thus combining everything together we have
$$
\frac d{d\tau}\Big (\dot x^{0} -  2 h_{0\ga}\dot x^{\ga} + \zeta
h_{00} \,\dot q \big )
= O(\ve (x^{0})^{-1})|\dot X|^2
$$
We integrate this identity between proper times $0<\tau$.
Observe that
$|\dot X|\le |\dot x^{0}|$ and that
$$
(x^{0})^{-1} |\dot x^{0}| = \frac d{d\tau} \ln x^{0}
$$
Thus
$$
\dot x^{0} (\tau) \les  \dot x^{0} (0) +
(2 h_{0\ga}\dot x^{\ga} - \zeta h_{00} \dot q )|_{0}^{\tau} +
\ve \int_{0}^{\tau} \frac d{d\tau} \ln x^{0} \dot x^{0} \,
d\tau'
$$
It follows that
$$
\dot x^{0} (\tau) \les \Big (\frac {|x^{0}(\tau)|}{|x^{0}(0)|}\Big )^{\ve}
\dot x^{0} (0)
$$
Integrating one more time  and assuming that $x^{0}(0)=t(0)=1$ we obtain
that
$$
(x^{0}(\tau))^{1-\ve} \les 1 + {\dot x^{0}(0)  }\tau
$$
From this we conclude that the time $x^{0}=t$ remains finite with $\tau$.
This concludes the proof for time-like geodesics.
\end{proof}

We now address the issue of null geodesics $X(\tau)$,
$$
g^{\a\b} \dot X^{\a} \dot X^{\b} =0
$$
Examining the proof above leads to the conclusion that
is suffices to establish that the condition
$\dot x^{0}(\tau) > 0$ is preserved in time.
\begin{lemma}\label{le:nullgeo}
For a future oriented inextendible null geodesic $X(\tau)$ defined
on the interval $[0,\tau_{*})$ we have $\dot x^{0}(\tau) >0$
for all $\tau\in [0,\tau_{*})$.
\end{lemma}
\begin{proof}
Let $\tau_{0}<\tau_{*}$ be the first time when $\dot x^{0}(\tau_{0})=0$.
Fix a sufficiently small small constant $c$. Then there exists
a small interval of size $\delta$ such that
$$
0\le \dot x^{0}(\tau) \le c,\qquad \forall \tau\in [\tau_{0}-\de,\tau_{0}]
$$
and
\beq
\label{eq:nonzero}
\dot x^{0}(\tau_{0}-\de) =c
\ee
We observe that \eqref{eq:timecontrol} with $A=0$ implies that
$|\dot X(\tau)|\le 2 |\dot x^{0}(\tau)|$ and therefore,
$$
|\dot X(\tau)| \le 2 c,\qquad \forall \tau\in [\tau_{0}-\de,\tau_{0}]
$$
Integrating the geodesic equation \eqref{eq:geod} we obtain
$$
|\dot x^{0}(\tau_{0}) - \dot x^{0}(\tau_{0}-\de) |\le
\int_{\tau_{0}-\de}^{\tau_{0}} |\Gamma| |\dot X|^{2} \le
\ve c^{2} \de
$$
Thus, using \eqref{eq:nonzero},
$$
\dot x^{0} (\tau_{0}) \ge c - \ve c^{2} \de >0
$$
Contradiction.
\end{proof}

This completes the proof of Proposition \ref{pr:Inext}.

We have shown that all future inextendible causal geodesics $X(\tau)$ exist
for all values of the affine parameter $\tau\in [0,\infty)$. This
means that the constructed space-time is future causally geodesically
complete. Next we establish that all future oriented causal geodesics
escape to infinity.
\begin{prop}\label{pr:Escape}
Let $X(\tau)$ be a future oriented causal geodesic. Then
\beq
\label{eq:escape}
\lim_{\tau\to \infty} |X(\tau)|=\infty
\eq
\end{prop}
\begin{proof}
The inequality \eqref{eq:timemon} immediately gives the desired result
for time-like geodesics.
Recall that by Lemma \ref{le:nullgeo} we have that
$\dot x^{0}(\tau)>0$ and thus $x^{0}(\tau)$ is monotonically increasing
in $\tau$. We now argue by contradiction.
Assume that for all $\tau\ge 0$
$$
|X(\tau)|\le C
$$
for some potentially large constant $C$. Then
there exists a time $t_{0}$ such that
$$
t_{0}=\lim_{\tau\to \infty} x^{0}(\tau)
$$
Set $\tau_{0}$ be the value of the proper time $\tau$ for which
$t(\tau_{0}) = t_{0} -\de$ for some small constant $\de$.
Integrating the geodesic equation we obtain that for $\tau\ge \tau_{0}$
\beq
\label{eq:arrive}
\dot x^{0}(\tau) = \dot x^{0}(\tau_{0}) +
\int_{\tau_{0}}^{\tau} |\Gamma| |\dot x^{0}|^{2}\,d\tau' \le
\dot x^{0}(\tau_{0}) + \ve \int_{t_{0}-\de}^{t}\dot x^{0}\, dt'
\le \dot x^{0}(\tau_{0}) + \ve \de \sup_{\tau_{0}\le \tau'\le\tau}
\dot x^{0}(\tau)
\ee
Thus for any $\tau\ge \tau_{0}$
\beq
\label{eq:clonull}
\dot x^{0}(\tau) \le 2 \dot x^{0}(\tau_{0})
\ee
Choosing a sequence of times $\tau_{0}\to \infty$ such that
$\dot x^{0}(\tau_{0})\to 0$ (such a sequence must exist, otherwise
$x^{0}(\tau)\to \infty$) we infer from \eqref{eq:clonull} that
$$
\dot x^{0}(\tau) \to 0
$$
as $\tau\to \infty$.
We can then choose small constant $c, \de$ such that
$t(\tau_{0}) = t_{0}-\de$ and
$$
\dot x^{0}(\tau_{0})=c, \qquad \dot x^{0}(\tau) \le c
$$
for all $\tau\ge \tau_{0}$.
Returning to \eqref{eq:arrive} we see that
$$
|\dot x^{0}(\tau) - c |\le \ve \de c
$$
Thus
$$
\dot x^{0}(\tau) \ge \frac c2
$$
for all $\tau\ge \tau_{0}$ and we obtained contradiction.
\end{proof}

\end{document}